\documentclass[12pt,a4paper,reqno]{amsart}

\usepackage[utf8]{inputenc} \usepackage{amsmath,amssymb,amsfonts} \usepackage{mathrsfs} \usepackage{enumitem}
\usepackage[table]{xcolor}
\usepackage[numbers,sort&compress]{natbib}
\usepackage[pdfencoding=auto,colorlinks]{hyperref}

\newtheorem{thm}{Theorem}[section] \newtheorem{cor}[thm]{Corollary} \newtheorem{lem}[thm]{Lemma} \newtheorem{prop}[thm]{Proposition} \theoremstyle{definition} \newtheorem{defn}[thm]{Definition} \newtheorem{rem}[thm]{Remark} \newtheorem{ex}[thm]{Example}

\numberwithin{equation}{section} \allowdisplaybreaks

\usepackage{geometry}
\renewcommand{\emptyset}{\varnothing}
\newcommand{\one}{\mathbf 1}
\newcommand{\rmref}[1]{{\rm\ref{#1}}}
\newcommand{\ol}{\overline}
\newcommand{\wt}{\widetilde}

\newcommand{\rank}{\operatorname{rank}}

\newcommand{\R}{\mathbb R} %
\newcommand{\C}{\mathbb C} %
\newcommand{\Q}{\mathbb Q} %
\newcommand{\N}{\mathbb N} %

\newcommand{\calA}{\mathcal A} \newcommand{\calB}{\mathcal B} \newcommand{\calD}{\mathcal D} \newcommand{\calE}{\mathcal E} \newcommand{\calF}{\mathcal F} \newcommand{\calH}{\mathcal H} \newcommand{\calK}{\mathcal K} \newcommand{\calM}{\mathcal M} \newcommand{\calN}{\mathcal N} \newcommand{\calO}{\mathcal O} \newcommand{\calX}{\mathcal X}

 \newcommand{\veps}{\varepsilon} 

\newcommand{\linspan}{\operatorname{span}} \newcommand{\dom}{\operatorname{dom}} \newcommand{\ran}{\operatorname{ran}}

 \newcommand{\mc}[1]{\mathcal#1} \DeclareMathOperator*{\esssup}{ess\,sup\,}

\title[Koopman operators: invariance, compactness, and self-adjointness]{Koopman operators on function spaces: invariance, compactness, and self-adjointness}

\author[M. Colbrook]{Matthew Colbrook} \address{\hspace*{-.43cm}{\rm\bf M. Colbrook:}\newline Department of Applied Mathematics and Theoretical Physics, University of Cambridge, Cambridge, United Kingdom} \email{\href{mailto:m.colbrook@damtp.cam.ac.uk}{m.colbrook@damtp.cam.ac.uk}} \urladdr{\url{https://www.damtp.cam.ac.uk/user/mjc249/home.html}}

\author[R. Moldenhauer]{Robert Moldenhauer} \address{\hspace*{-.43cm}{\rm\bf R. Moldenhauer:}\newline Control and Signal Processing Group, Department of Electrical and Electronic Engineering, University of Melbourne, Melbourne, Australia\newline Centre de Recherche en Automatique de Nancy (CRAN), Université de Lorraine, Nancy, France} \email{\href{mailto:moldenhauer.r@student.unimelb.edu.au}{moldenhauer.r@student.unimelb.edu.au}}

\author[F.M. Philipp]{Friedrich M. Philipp} \address{\hspace*{-.43cm}{\rm\bf F.M. Philipp:}\newline Institute of Mathematics, Technische Universität Ilmenau, Ilmenau, Germany} \email{\href{mailto:friedrich.philipp@tu-ilmenau.de}{friedrich.philipp@tu-ilmenau.de}} \urladdr{\url{https://www.tu-ilmenau.de/obc/team/friedrich-philipp}}

\author[M. Schaller]{Manuel Schaller} \address{\hspace*{-.43cm}{\rm\bf M. Schaller:}\newline Chemnitz University of Technology, Chemnitz, Germany} \email{\href{mailto:manuel.schaller@mathematik.tu-chemnitz.de}{manuel.schaller@mathematik.tu-chemnitz.de}} \urladdr{\url{https://www.tu-chemnitz.de/mathematik/jpnuma/people/schaller.php}}

\author[K. Worthmann]{Karl Worthmann} \address{\hspace*{-.43cm}{\rm\bf K. Worthmann:}\newline Optimization-Based Control Group, Institute of Mathematics, Technische Universität Ilmenau, Ilmenau, Germany} \email{\href{mailto:karl.worthmann@tu-ilmenau.de}{karl.worthmann@tu-ilmenau.de}} \urladdr{\url{https://www.tu-ilmenau.de/obc/team/karl-worthmann}}

\subjclass[2020]{Primary 47B33; Secondary 37A30, 46E30, 47B07, 47B32, 60J35}
\keywords{Koopman operator, composition operator, Markov kernel, invariant function space, compact operator, self-adjoint operator, reproducing kernel Hilbert space}

\begin{document}
\begin{abstract}
Koopman operator theory studies nonlinear dynamical systems through the linear evolution of observable functions in a typically infinite-dimensional function space. 
This linear formulation makes tools from operator and semigroup theory available, including spectral analysis, functional calculus, resolvent and perturbation methods, and Galerkin-type approximation with convergence guarantees.
However, the spectral properties of the resulting Koopman operator depend strongly on the chosen space of observables.
In this paper, we investigate key properties of the 
Koopman operator for both deterministic and stochastic dynamics: boundedness, compactness, and self-adjointness. 
We give characterizations for Koopman operators acting on bounded continuous function spaces, on $L^p$ spaces with $1\le p\le\infty$, and on reproducing kernel Hilbert spaces (RKHSs), and link them to properties of the transition kernel and the dynamics. 
While the resulting conditions for boundedness are broadly applicable, those characterizing compactness and self-adjointness turn out to be rather restrictive in the deterministic case, as illustrated by the example of the Gaussian RKHS. 
For stochastic dynamics, smoothing can yield compactness, while reversibility yields self-adjointness. Langevin semigroups illustrate how both properties can occur beyond the restrictive
deterministic setting.
\end{abstract}

\maketitle

\section{Introduction}\label{sec:introduction}

Koopman's construction begins by replacing the evolution of states with the evolution of observables. If $F:\calX\to\calX$ advances a state by one time step, an observable $f$ evolves according to
\begin{equation}\label{eq:intro_deterministic_koopman}
(\calK f)(x)=f(F(x)). 
\end{equation}
Although $F$ need not be linear, the map $f\mapsto\calK f$ is linear. For an introduction to Koopman operator theory, see the edited volume~\cite{MaurMeziSusu2020}, the book \cite[Chapter 11]{colbrook2026spectral}, the recent tutorial paper~\cite{mezic2026koopman}, and the review and overview articles~\cite{Brunton2022ModernKoopman,Colbrook2026}; for a survey of applications, see~\cite{Budisic2012}. 
Koopman introduced the Hilbert-space viewpoint for Hamiltonian systems in 1931 \cite{Koop31}. 
Koopman and von Neumann soon further developed this concept in their study of systems with continuous spectra \cite{KoopNeum32}. 
This viewpoint became one of the operator-theoretic foundations of ergodic theory, whose early development included von Neumann's mean ergodic theorem~\cite{vonNeumann1932quasi}; see also~\cite{Halmos1956,EisnerFarkasHaaseNagel2015}.

On the other hand, Formula~\eqref{eq:intro_deterministic_koopman} also belongs to the classical theory of composition operators. 
Beginning with Littlewood's inequalities \cite{Littlewood1925}, that theory has sought to read properties of the symbol $F$ from boundedness, compactness~\cite{%
SinghKumar1979compact}, normality~\cite{singh1975normal}, and spectral behavior~\cite{Ridge1973Spectrum} of $f\mapsto f\circ F$ on a prescribed function space, see, e.g., \cite{SinghManhas1993}; 
for the analytic setting we refer to the monograph~\cite{CowenMacCluer1995} and the survey~\cite{nordgren2006composition}. 
Both traditions inform the modern spectral analysis of nonlinear systems~\cite{rowley2009spectral,WillRowl15,mezic2020spectrum}, data-driven approximation \cite{WillKevr15,KutzBrunton2016,KordMezi18:EDMD:convergence,NuskPeit23}, and kernel methods~\cite{Klus2020KernelEDMD}, see also~\cite{mezic2021koopman,PhilScha24:kernel}. 
Numerical schemes and convergence guarantees for the resulting spectral computations are developed, e.g., in~\cite{Mezi22,ColbTown24,LlamLlamLatzKlus24,doi:10.1137/24M1662370,colbrook2023mpedmd,conradie2026trustworthy,bartel2026extended}; see the overview \cite{COLBROOK2024127}, whereas finite-data error bounds for extended dynamic mode decomposition (EDMD) and kernel EDMD are provided in~\cite{PhilScha23,PhilScha24,KohnPhil24}, see also~\cite{PhilScha23,bold2025kernel} for extensions to systems with inputs; see the overview paper~\cite{StraWort26}.
Beyond spectral analysis, Koopman surrogate models are widely used for prediction, estimation, and control~\cite{proctor2018generalizing,KordMezi18:MPC,otto2021koopman,StraWort26}, and in the analysis of molecular dynamics~\cite{noenusk13}. In all of these applications, the choice of the observable space is essential.
The spectral viewpoint also extends naturally to random dynamical systems~\cite{CrnjaricZic2020,colbrook2024beyond}. We treat stochastic dynamics throughout using Markov kernels. In this setting, the Koopman operator is dual to the transfer, or Perron--Frobenius, operator~\cite{KlusKoltaiSchuette2016,SchlKord22} and can be viewed as a conditional expectation operator~\cite{Mollenhauer2020}. Sharp classifications of the computational difficulty of approximating
the spectra of Koopman and transfer operators from data have been
established in \cite{boulle2025convergent,colbrook2026adversarial}.

Formula~\eqref{eq:intro_deterministic_koopman} need not, however, define an operator on a prescribed observable space: composition can fail even to descend to equivalence classes modulo sets of measure zero, see Example~\ref{ex:equivalence_classes:counterexample}. 
In this work we address the following fundamental question: when does the Koopman rule define an operator on a prescribed function space, and when is that operator bounded, compact, or self-adjoint? 
Each of these properties has direct consequences for the analysis of the operator and for its data-driven approximation: boundedness places $\mathcal K$ within the usual theory of operators on Banach spaces, compactness makes the nonzero spectrum discrete with finite multiplicities, and self-adjointness brings the spectral theorem and, for Markov dynamics, an exact probabilistic meaning, see also~\cite{Shapiro1993} for precursor results on invariance, boundedness, and compactness of composition operators on Hardy and Bergman-type spaces.

This question is not merely formal. Koopman invariance is the standing hypothesis behind
error analysis for EDMD and kernel EDMD~\cite{PhilScha24:kernel,KohnPhil24,hertel2025koopman}, and RKHS
invariance in particular is known to impose strong rigidity on the underlying dynamics, forcing
affine $F$ for the Gaussian kernel~\cite{Gonzalez2024} and analogous restrictions for several other kernel families~\cite{Ikeda22}. 
For stochastic dynamics, these rigidity restrictions need not persist. For instance, sufficiently regular stochastic differential equations yield smoothing maps into suitably chosen RKHSs~\cite[Section~4.3]{hertel2026verifiable}; see also~\cite{Philipp2025} for Gaussian RKHSs and SDEs with constant matrix coefficients. A complementary viewpoint focusing on the spectral side is developed in~\cite{DasGiannakis2020}, and the kernel Banach-space setting is treated in~\cite{IkedaIshikawaSchlosser2022}.
Boundedness of composition operators on individual spaces, $L^p$, $C_b$, and RKHSs, has
likewise been studied separately in the classical composition-operator literature, see, e.g., \cite{SinghManhas1993}, and \cite{Shapiro1993} for the analytic setting.
What has been missing is a treatment that (i) jointly covers deterministic \emph{and} stochastic
(Markov) dynamics, and (ii) precisely characterizes invariance, boundedness, compactness, and self-adjointness --- rather than assuming them, as is common practice in the data-driven Koopman literature.
We largely close this gap by giving necessary and/or sufficient conditions for well-definedness, invariance,
boundedness, compactness, and self-adjointness of the Koopman operator on all three classes of
spaces, for both deterministic maps and Markov kernels. 
In particular, we provide a rather complete list of characterizations and sufficient conditions, identifying the relevant prior results explicitly, so that this paper may also serve as a reference work.

Beyond providing a unified treatment of existing results, this work makes several contributions to the functional-analytic theory of Koopman operators, particularly regarding compactness and self-adjointness. For stochastic dynamics, we characterize compactness on $L^\infty(\mu)$ in terms of total boundedness of the family of transition measures in the total variation norm outside a common $\mu$-null set (Theorem~\ref{t:compact:stoch:Linf}), and establish a characterization of compactness on $L^1(\mu)$ solely in terms of the underlying Markov kernel (Theorem~\ref{t:compact2:Lp}). A contribution of independent interest, going beyond Koopman operator theory, is a characterization of compact embeddings between arbitrary reproducing kernel Hilbert spaces in terms of the eigenvalues of matrices constructed from finite Gram matrices (Theorem~\ref{t:RKHS_comp_emb}). This result requires neither a topological nor a measure-theoretic structure on the underlying set and yields kernel-based characterizations of compactness of both deterministic and stochastic Koopman operators on RKHSs (Theorems~\ref{t:RKHS_comp_det} and \ref{t:RKHS_comp}). Furthermore, we derive necessary conditions for compactness in the deterministic RKHS setting (Proposition~\ref{p:nec}), revealing geometric restrictions on the underlying dynamics as well as restrictions on the asymptotic behavior of functions in the RKHS. In particular, these conditions rule out compactness on Gaussian and suitable Sobolev RKHSs over unbounded domains. Regarding self-adjointness, we establish a kernel-based criterion for deterministic Koopman operators on arbitrary RKHSs (Theorem~\ref{t:selfadjprep}). For a broad class of translation-invariant kernels, this leads to a structural characterization of self-adjointness in terms of kernel preservation and involutivity of the underlying dynamics (Corollary~\ref{c:rkhs:sa}), paralleling the classical characterization in the $L^2(\mu)$ setting.

In Table~\ref{tab:overview}, we provide an overview of the main theorems in this paper. For general Markov kernels and intermediate exponents $p\in(1,\infty)$, characterizations of
both Koopman invariance and compactness in the $L^p(\mu)$ setting remain open. In this setting, the present paper establishes only sufficient conditions.

\begin{table}[htb]
    \centering
    \begin{tabular}{|l|l|l|l|}
        \hline
        & Invariance & Compactness & Self-adjointness \\
        \hline
        $C_b(\calX)$    & Theorems~\ref{t:det_invariance_Cb} / \ref{t:inv_C_b}             & Theorems~\ref{t:C_det_compact} / \ref{t:Cb_compact}            & \hspace{1.4cm}---                                                        \\
        \hline
        $L^\infty(\mu)$ & Theorems~\ref{t:det:invariance:Linfty} / \ref{t:stoch:invariance:Linfty} & Theorems~\ref{t:Linfty_comp_det} / \ref{t:compact:stoch:Linf}  & \hspace{1.4cm}---                                                \\
        \hline
        $L^1(\mu)$      & Theorems~\ref{t:det:invariance:Lp} / \ref{t:invariance:stoch:Lp} & Theorems~\ref{t:compact_Lp_det} / \ref{t:compact2:Lp}          & \hspace{1.4cm}---                                                        \\
        \hline
        $L^p(\mu)$      & Theorems~\ref{t:det:invariance:Lp} / \colorbox{gray!30}{\ref{t:invariance:stoch:Lp}} & Theorems~\ref{t:compact_Lp_det} / \colorbox{gray!30}{\ref{t:Lp_comp_interpolation}}& \hspace{1.4cm}---                                                        \\
        \hline\hline
        $L^2(\mu)$      & Theorems~\ref{t:det:invariance:Lp} / \colorbox{gray!30}{\ref{t:invariance:stoch:Lp}} & Theorems~\ref{t:compact_Lp_det} / \colorbox{gray!30}{\ref{t:Lp_comp_interpolation}}& Theorems~\ref{t:det:selfadjL2} / \ref{t:sa_stoch}           \\
        \hline
        RKHS            & Theorems~\ref{t:RKHS_det_inv} / \ref{t:RKHS_inv}                 & Theorems~\ref{t:RKHS_comp_det} / \ref{t:RKHS_comp}             & Theorems~\ref{t:selfadjprep} / \ref{t:selfadj:stoch:rkhs}   \\
        \hline 
\end{tabular}

    \vspace{.5cm}
    \caption{Characterizations and sufficient conditions (in gray color) obtained in this work. Each entry lists the deterministic result first and the stochastic result second. A dash indicates that the property is not meaningful on the respective space.}
    \label{tab:overview}
\end{table}

Section~2 constructs the Koopman operator, first pointwise and then on equivalence classes; Section~3 treats invariance and boundedness using the closed graph theorem and positivity arguments; Section~4 develops compactness; and Section~5 characterizes self-adjointness on $L^2$ and on RKHSs. In Sections 3--5, the deterministic case is treated first, and the general Markov-kernel statement is specialized to it whenever the two arguments coincide.

\section{The Koopman operator} \label{sec:definition}
The deterministic formula \eqref{eq:intro_deterministic_koopman} defines a map on $\calM(\calX)$, whereas the integral in \eqref{eq:informal} need not be well defined at every point. To obtain an operator on $L^p(\mu)$, one must establish both integrability and independence of the chosen representative. We begin with the pointwise rule and then determine when it descends modulo $\mu$-almost-everywhere equality.

\subsection{Pointwise Koopman operators}
Let $\calX$ be a nonempty Polish space, let $\calB_\calX$ denote its Borel $\sigma$-algebra, and let $\calM(\calX)$ be the space of complex-valued Borel measurable functions on $\calX$. Measurability between topological spaces will always refer to their Borel $\sigma$-algebras.

\medskip \noindent\textbf{The deterministic case.} Consider the system
\begin{align}\label{eq:deterministic_system}
x^+ = F(x),
\end{align}
generated by a measurable map $F:\mc X\to\mc X$. Its \emph{Koopman operator} acts on observables by composition:
\begin{align}\label{eq:det:Koopman_operator}
\mc K:\mc M(\mc X)\to\mc M(\mc X),\qquad f\mapsto f\circ F.
\end{align}
The relevant question is usually whether this composition operator preserves a prescribed subspace $\mc F\subset\mc M(\mc X)$.

\begin{defn}\label{def:det:Koopman:invariance} A normed space $\calF\subset\calM(\calX)$ is \emph{Koopman invariant} for the deterministic system \eqref{eq:deterministic_system} if $\calK f\in\calF$ for every $f\in\calF$.
\end{defn}

\noindent\textbf{The stochastic case.} In the stochastic setting we consider, the dynamics are generated by a \emph{Markov kernel} on $\calX$, which is a map $\rho:\calX\times\calB_\calX\to[0,1]$ with the following properties: \begin{itemize} \item[(a)] $\rho_x := \rho(x,\cdot)$ is a probability measure for each $x\in\calX$. \item[(b)] For each $\Delta\in\calB_\calX$, we have $\rho(\,\cdot\,,\Delta)\in\calM(\calX)$.
\end{itemize}
Conditional on the current state $x\in\calX$, the next state belongs to $\Delta\in\calB_\calX$ with probability $\rho_x(\Delta)$.

Deterministic dynamics form the degenerate case
\begin{align}\label{eq:deterministic_special_case_of_stochastic}
\rho_x = \delta_{F(x)},\qquad x\in\mc X.
\end{align}
Hence, deterministic dynamics constitute a particular instance of the stochastic setting.

\begin{rem}\label{remark:Markov_kernel_correspondence} Equivalently, the dynamics may be represented as $x^+=G(x,v)$, where $G:\mc X\times\mc V\to\mc X$ is measurable and the independent random inputs have common law $\eta$ on the standard Borel space $\mc V$. The associated Markov kernel is
\begin{align}
\rho_x(\Delta) = \eta\bigl(\{v:G(x,v)\in\Delta\}\bigr), \qquad \Delta\in\mc B_\calX.\label{eq:Markov_kernel_correspondence}
\end{align}
Conversely, the kernel-randomization lemma \cite[Lemma 4.22]{Kallenberg2021} represents every Markov kernel on the Polish space $\mc X$ in this form with $\mc V=[0,1]$ and $\eta$ equal to Lebesgue measure. We use the kernel formulation throughout.
\end{rem}

For $f\in\mc M(\mc X)$, the set on which $f$ is $\rho_x$-integrable is
\begin{equation}\label{eq:Delta_f}
\Delta_f:=\left\{x\in\mc X:\int |f|\,d\rho_x<\infty\right\}.
\end{equation}
The zero-extended integral is
\begin{equation}\label{eq:zero_extended_K}
(\calK_0f)(x):= \begin{cases} \displaystyle\int f\,d\rho_x,&x\in\Delta_f,\\[2mm] 0,&x\notin\Delta_f.
\end{cases}
\end{equation}
Lemma~\ref{l:measurable} gives $\Delta_f\in\calB_\calX$ and $\calK_0f\in\mc M(\mc X)$. The zero extension need not be linear on all of $\mc M(\mc X)$, but it is linear on any vector space consisting of functions that belong to $L^1(\rho_x)$ for every $x$. On such a space the subscript is suppressed:
\begin{align}\label{eq:informal}
(\calK f)(x):=(\calK_0f)(x)=\int f\,d\rho_x =\mathbb E_{\rho_x}[f],\qquad x\in\calX.
\end{align}
In probability theory, the operator $\calK$ is commonly referred to as the \emph{Markov transition operator} and has been studied extensively in the classical theory of Markov processes; see, e.g., \cite{Doob53}, \cite{Feller53}, and \cite{Dynkin65}. In the modern Koopman literature, the Markov-kernel formulation of stochastic Koopman operators goes back at least to \cite{CrnjaricZic2020}.

\begin{defn}\label{def:stoch:Koopman:invariance} A linear space $\calF\subset\calM(\calX)$ is \emph{Koopman invariant} for $\rho$ if, for every $f\in\calF$, one has $\Delta_f=\calX$ and $\calK f\in\calF$.
\end{defn}

For the deterministic kernel \eqref{eq:deterministic_special_case_of_stochastic}, every $f\in\mc M(\mc X)$ belongs to $L^1(\rho_x)$, since $\int |f|\,d\rho_x=|f(F(x))|<\infty$. Definition~\ref{def:stoch:Koopman:invariance} therefore reduces to Definition~\ref{def:det:Koopman:invariance}.

Positivity converts invariance into boundedness on Banach lattices. Its pointwise form appears below; Theorem~\ref{t:schaefer_ec} gives the corresponding assertion on equivalence classes. Since complex-valued functions carry no compatible order, the assertion is formulated on a real Banach lattice and then passed to its complexification. The underlying result is the standard theorem that every positive, everywhere-defined operator between Banach lattices is bounded \cite[Theorem~4.3]{aliprantis_positive_2006}.

\begin{thm}\label{t:schaefer}
Let $\calF_{\mathbb R}$ be a Banach lattice of real-valued measurable functions on $\calX$ under the pointwise order, and let $\calF$ be its complexification, equipped with a Banach space norm equivalent to
\begin{align*}
f\longmapsto \|\operatorname{Re}f\|_{\calF_{\mathbb R}} +\|\operatorname{Im}f\|_{\calF_{\mathbb R}}.
\end{align*}
If $\calF$ is Koopman invariant, then $\calK:\calF\to\calF$ is bounded.
\end{thm}
\begin{proof} The restriction of $\calK$ to $\calF_{\mathbb R}$ is positive. If it were unbounded, there would be $g_n\in\calF_{\mathbb R}$ with $\|g_n\|=1$ and $\|\calK g_n\|\ge4^n$. For $f_n=|g_n|$, one has $\|f_n\|=1$, and positivity gives $|\calK g_n|\le\calK f_n$, hence $\|\calK f_n\|\ge\|\calK g_n\|$. The series $f=\sum_{n\ge1}2^{-n}f_n$ converges in the Banach lattice. Since its positive cone is closed, $f\ge2^{-n}f_n$ for every $n$. Positivity and monotonicity of the lattice norm now yield
\begin{align*}
\|\calK f\|\ge2^{-n}\|\calK f_n\|\ge2^n,
\end{align*}
a contradiction. Hence $\calK$ is bounded on the real part. Applying this bound separately to the real and imaginary parts, and then using norm equivalence, proves boundedness on $\calF$.
\end{proof}

\subsection{Koopman operators on equivalence classes}

The $L^p$ theory concerns equivalence classes modulo almost-everywhere equality, rather than pointwise functions. Let $\mu$ be a $\sigma$-finite Borel measure on $\mc X$ with $\mu(\calX)>0$. For $f,g\in\mc M(\mc X)$, let $f\sim g$ mean that $f=g$ $\mu$-almost everywhere, and let $[f]$ be the equivalence class of $f$. For $\calF\subset\calM(\calX)$, the corresponding family of classes is
\begin{align*}
[\calF] := \{[f] : f\in\calF\}.
\end{align*}
The issue is whether $\calK_0f$ and $\calK_0g$ determine the same class whenever $f\sim g$. We resolve this first for the general Markov kernel; the deterministic case then follows by specialization.

At the level of equivalence classes, $\rho_x$-integrability is needed only for $\mu$-almost every $x$. This leads to the linear space
\begin{equation}\label{eq:defnDmu}
\calD(\mu) := \big\{f\in\calM(\calX) : \mu(\calX\backslash\Delta_f)=0\big\},
\end{equation}
where $\Delta_f$ is defined in \eqref{eq:Delta_f}. If $f,g\in\calD(\mu)$ and $a,b\in\mathbb C$, then $af+bg$ is $\rho_x$-integrable for all $x$ in the full-measure set $\Delta_f\cap\Delta_g$ and
\begin{equation}\label{eq:K0_linearity_ae}
\calK_0(af+bg)=a\calK_0f+b\calK_0g\qquad\mu\text{-a.e.}
\end{equation}
Thus $\calD(\mu)$ is linear, and \eqref{eq:K0_linearity_ae} provides the almost-everywhere linearity of $\calK_0$ used below. In the deterministic case, $\Delta_f=\mc X$ for every $f\in\mc M(\mc X)$, and hence $\calD(\mu)=\calM(\calX)$.

Representative independence requires $\calK_0g\sim\calK_0f$ whenever $g\sim f$. The maximal domain on which this holds is
\begin{equation}\label{eq:wt_calD}
\wt\calD(\mu) := \big\{[f] : f\in\calD(\mu),\; g\sim f\,\Rightarrow\,(g\in\calD(\mu)\,\text{ and }\, \calK_0g \sim \calK_0f)\big\}.
\end{equation}
Proposition~\ref{p:equiv_classes} shows that, as soon as this set is nonempty, it is a linear subspace of $[\calM(\calX)]$. In the deterministic case, \eqref{eq:wt_calD} reduces to
\begin{align*}
\wt\calD(\mu) = \big\{[f] : f\in\calM(\calX),\; g\sim f\,\Rightarrow\,g\circ F \sim f\circ F\big\}.
\end{align*}
Example~\ref{ex:equivalence_classes:counterexample} shows that these conditions are substantive: even $[0]$ need not belong to $\wt\calD(\mu)$.

For $[f]\in\wt\calD(\mu)$, the induced operator is therefore unambiguous:
\begin{align}\label{eq:Koopman_equivalence_classes}
\calK[f] := [\calK_0 f].
\end{align}
Equation~\eqref{eq:K0_linearity_ae} shows that this induced operator is linear.

Representative independence is governed by the one-step measure
\begin{align}\label{e:nu}
\nu(\Delta) := \int\rho_x(\Delta)\,d\mu(x),\qquad\Delta\in\calB_\calX.
\end{align}
The countable additivity of the measures $\rho_x$, together with monotone convergence, proves that $\nu$ is countably additive; the value $+\infty$ is allowed. Lemma~\ref{l:nu} provides the corresponding Tonelli--Fubini formula and its precise integrability consequences. For the deterministic kernel, $\rho_x=\delta_{F(x)}$, $\nu$ is the push-forward measure
\begin{align}\label{e:nu_det}
\nu(\Delta) = \left(\mu\circ F^{-1}\right)(\Delta),\qquad\Delta\in\mc B_{\mc X}.
\end{align}
Although $\mu$ is $\sigma$-finite, $\nu$ need not be: on $\calX=\N_0$ with counting measure and $F(n)=0$, one has $\nu(\Delta)=\infty$ when $0\in\Delta$ and $\nu(\Delta)=0$ otherwise. Thus the Radon--Nikodym theorem may be invoked only after the $\sigma$-finiteness of $\nu$ has been established. As usual, $\nu\ll\mu$ means that every $\mu$-null Borel set is $\nu$-null.

\begin{prop}\label{p:equiv_classes} The following statements are equivalent: \begin{enumerate} \item[{\rm (i)}] The maximal domain $\wt\calD(\mu)$ is nonempty, i.e., $\wt\calD(\mu) \neq\emptyset$; \item[{\rm (ii)}] $[0]\,\in\,\wt\calD(\mu)$; \item[{\rm (iii)}] $[\calD(\mu)] = \wt\calD(\mu)$; \item[{\rm (iv)}] $\nu$ is absolutely continuous with respect to $\mu$, i.e., $\nu\ll\mu$.
\end{enumerate}
\end{prop}
\begin{proof} (i)$\Rightarrow$(ii). We first show that a nonempty maximal domain contains the zero class. If $[f]\in\wt\calD(\mu)$ and $g\sim-f$, then $-g\sim f$; hence $-g\in\calD(\mu)$ and $\calK_0(-g)\sim\calK_0f$. Since $\calD(\mu)$ is a vector space and $\calK_0$ is linear almost everywhere on it, this proves $[-f]\in\wt\calD(\mu)$.

Likewise, let $[f],[h]\in\wt\calD(\mu)$ and let $g\sim f+h$. Then $g-h\sim f$, so $g-h\in\calD(\mu)$ and $\calK_0(g-h)\sim\calK_0f$. Since $h\in\calD(\mu)$, it follows that $g\in\calD(\mu)$ and
\begin{align*}
\calK_0g \sim \calK_0(g-h)+\calK_0h \sim \calK_0f+\calK_0h \sim \calK_0(f+h).
\end{align*}
Thus, $[f+h]\in\wt\calD(\mu)$. Consequently, any $[f]\in\wt\calD(\mu)$ gives $[0]=[f]+[-f]\in\wt\calD(\mu)$.

(ii)$\Rightarrow$(iv). Let $\Delta\in\calB_\calX$ satisfy $\mu(\Delta)=0$ and take $g=\one_\Delta$. Since $g\sim0$, condition (ii) gives $\calK_0g\sim0$. But $(\calK_0g)(x)=\rho_x(\Delta)$ for every $x$, and therefore $\nu(\Delta)=0$.

(iv)$\Rightarrow$(iii). Let $f\in\calD(\mu)$ and $g\sim f$. There is a Borel set $N$ with $\mu(N)=0$ such that $f=g$ on $N^c$. Since $\nu(N)=0$, the nonnegative measurable function $x\mapsto\rho_x(N)$ vanishes for $\mu$-almost every $x$. For all such $x$,
\begin{align*}
\int|g|\,d\rho_x=\int|f|\,d\rho_x, \qquad \int g\,d\rho_x=\int f\,d\rho_x
\end{align*}
whenever the integrals are finite. Hence $g\in\calD(\mu)$ and $\calK_0g\sim\calK_0f$, proving $[f]\in\wt\calD(\mu)$. The reverse inclusion $\wt\calD(\mu)\subset[\calD(\mu)]$ is part of the definition.

(iii)$\Rightarrow$(i) follows because $0\in\calD(\mu)$.
\end{proof}

\begin{rem} The condition $\nu\ll\mu$ does \emph{not} imply $\rho_x\ll\mu$ for $\mu$-almost every $x\in\calX$. It says instead that, for each fixed $\mu$-null set $\Delta\in\calB_\calX$, one has $\rho_x(\Delta)=0$ for $\mu$-almost every $x$. The exceptional set of points $x$ may depend on $\Delta$. For example, take $\calX=(0,1)$, $\mu$ equal to Lebesgue measure, and $\rho_x=\delta_x$. Then $\nu=\mu\ll\mu$, although $\rho_x\not\ll\mu$ for every $x\in\calX$.
\end{rem}

\begin{rem}[Invariant measure] The measure $\mu$ is called \emph{invariant} or \emph{stationary} if $\nu=\mu$. For an invariant measure, Jensen's inequality and Tonelli's theorem give
\begin{align*}
\|\calK f\|_{L^p(\mu)}^p \le \int\calK(|f|^p)\,d\mu =\int|f|^p\,d\mu\qquad(f\in L^p(\mu),\ 1\le p<\infty),
\end{align*}
and
\begin{align*}
\|\calK f\|_{L^\infty(\mu)}\le\|f\|_{L^\infty(\mu)}\qquad(f\in L^\infty(\mu)).
\end{align*}
Thus, $\calK$ is a contraction on every $L^p(\mu)$, including $p=\infty$. We do not assume stationarity in general.
\end{rem}

For deterministic dynamics, a familiar change-of-variables hypothesis ensures the equivalent conditions of Proposition~\ref{p:equiv_classes}.

\begin{prop}\label{prop:det:special_case_absolute_continuous} If $\mc X\subset\R^n$ is open, $F:\mc X\to\mc X$ is a diffeomorphism, and $\mu$ is equivalent to Lebesgue measure $\lambda$, that is, $\mu\ll\lambda$ and $\lambda\ll\mu$, then $\nu\ll\mu$.
\end{prop}
\begin{proof} Since $\mu\ll\lambda$, the Radon--Nikodym theorem gives a density $h:=\frac{d\mu}{d\lambda}\in\calM(\calX)$. Let $\Delta\in\mc B_\calX$ satisfy $\mu(\Delta)=0$. The converse absolute continuity $\lambda\ll\mu$ gives $\lambda(\Delta)=0$, and the change-of-variables formula yields
\begin{align*}
\nu(\Delta) &=\mu(F^{-1}(\Delta)) =\int_{F^{-1}(\Delta)}h(x)\,dx\\ &=\int_\Delta h(F^{-1}(y))\cdot|\det D(F^{-1})(y)|\,dy=0.
\end{align*}
The last expression vanishes because $\lambda(\Delta)=0$.
\end{proof}

Under the equivalent conditions of Proposition~\ref{p:equiv_classes}, the pointwise rule \eqref{eq:informal} induces a linear operator on every class arising from $\calD(\mu)$. This gives the appropriate notion of invariance for quotient spaces.

\begin{defn}\label{def:Koopman-invariance:equivalence-class} A normed space $\wt\calF$ of equivalence classes {\rm (}w.r.t.\ $\mu$-a.e.\ equality{\rm )} is called \emph{Koopman invariant} if $\wt\calF\subset\wt\calD(\mu)$ and $\calK [f]\in\wt\calF$ for every $[f]\in\wt\calF$.
\end{defn}

On quotient spaces, brackets are henceforth suppressed when no ambiguity can arise: $f$ and $\calK f$ denote the classes $[f]$ and $\calK[f]$, respectively, and a Borel representative is understood whenever a pointwise expression occurs. Whenever $\mu$ is fixed, $\|\cdot\|_p$ abbreviates $\|\cdot\|_{L^p(\mu)}$, with the analogous convention for operator norms.

\begin{thm}\label{t:schaefer_ec}
Let $\wt\calF_{\mathbb R}$ be a Banach lattice of equivalence classes {\rm (}w.r.t.\ $\mu$-a.e.\ equality{\,\rm )} of real-valued functions, ordered almost everywhere, and let $\wt\calF$ be its complexification with a Banach norm equivalent to $\|\operatorname{Re}[f]\|_{\wt\calF_{\mathbb R}}+ \|\operatorname{Im}[f]\|_{\wt\calF_{\mathbb R}}$. If $\wt\calF$ is Koopman invariant, then $\calK:\wt\calF\to\wt\calF$ is bounded.
\end{thm}
\begin{proof}
If $[f]\ge0$ in $\wt\calF_{\mathbb R}$, choose a nonnegative measurable representative $f$. Then $\calK_0f\ge0$ almost everywhere, and representative independence makes the induced Koopman operator positive. The series argument from Theorem~\ref{t:schaefer}, now with inequalities understood almost everywhere, proves boundedness on the real part. Boundedness on $\wt\calF$ then
follows from complex linearity and the assumed
norm equivalence.
\end{proof}

\section{Koopman invariance and boundedness}\label{sec:invariance}
Koopman invariance is the key ingredient to prove %
boundedness of the Koopman operator and to establish %
bounds on the approximation error in Koopman learning; %
see, for example, \cite{KohnPhil24,hertel2025koopman}. Its characterization depends on the chosen space of observables.
We develop sufficient (and mostly also necessary) conditions on the dynamics for Koopman invariance,
first for deterministic maps and then for general Markov kernels.

Invariance of $C_b(\calX)$ will turn out to be equivalent to continuity of $F$ and, for a Markov kernel, to its weak continuity (Theorems~\ref{t:det_invariance_Cb} and \ref{t:inv_C_b}). 
Absolute continuity of $\nu$ w.r.t.\ $\mu$ removes the ambiguity caused by null sets and characterizes %
invariance of $L^\infty(\mu)$, where $\calK$ has norm one (Theorems~\ref{t:det:invariance:Linfty} and \ref{t:stoch:invariance:Linfty}). 
For deterministic dynamics, invariance at a finite exponent is equivalent to the stronger domination condition $\nu\le C\mu$ (Theorem~\ref{t:det:invariance:Lp}). For Markov kernels, this equivalence holds at $p=1$; invariance at a larger exponent need not imply domination, but it propagates to all larger exponents (Theorem~\ref{t:invariance:stoch:Lp}).
On an RKHS, transported kernel sections replace densities, and invariance is characterized by the kernel order $k_F\preceq c^2k$ or $k_\rho\preceq c^2k$; by the closed graph theorem, it then implies boundedness (Theorems~\ref{t:RKHS_det_inv} and \ref{t:RKHS_inv}). 
Finally, we show that Sobolev regularity of the transition density yields the smoothing property $\calK:L^2(\Omega)\to H^\ell(\Omega)$ on a bounded Lipschitz domain $\Omega\subset\R^n$ and hence invariance of suitable Mat\'ern and Wendland RKHSs (Theorem~\ref{t:sob_stoch}).

\subsection{Deterministic systems}\label{s:invariance:det} Throughout this subsection, $F$ generates the deterministic system \eqref{eq:deterministic_system}. We use \eqref{eq:det:Koopman_operator} on functions and \eqref{eq:Koopman_equivalence_classes} on equivalence classes.

\subsubsection{Invariance on \texorpdfstring{$C_b(\calX)$}{Cb(X)}} \label{s:invariance:det:Cb} Let $C_b(\calX)$ denote the Banach space of bounded continuous complex-valued functions with the supremum norm. The estimate $\|f\circ F\|_\infty\le\|f\|_\infty$ shows that composition is a contraction whenever $C_b(\calX)$ is Koopman invariant. Continuity is therefore the only obstruction.

\begin{thm}\label{t:det_invariance_Cb}
The space $C_b(\calX)$ is Koopman invariant for the deterministic system \eqref{eq:deterministic_system} if and only if $F$ is continuous. In this case, $\mc K:C_b(\mc X)\to C_b(\mc X)$ is bounded with norm $\|\mc K\|_{C_b(\mc X)\to C_b(\mc X)}=1$.
\end{thm}
\begin{proof} Since $\calK f=f\circ F$, continuity of $F$ immediately implies invariance. Conversely, suppose that $f\circ F$ is continuous for every $f\in C_b(\calX)$. Let $d$ be a metric inducing the topology of $\calX$\footnote{Such a metric exists because $\calX$ is Polish.}, and let $B_r(z):=\{x\in\calX:d(x,z)<r\}$. Fix $x_0\in\calX$ and $\veps>0$. The bounded continuous function $f(z)=\min\{d(z,F(x_0)),\veps\}$ vanishes at $F(x_0)$. Continuity of $f\circ F$ therefore gives $\delta>0$ such that $f(F(x))<\veps$ whenever $x\in B_\delta(x_0)$. By the definition of $f$, this means $d(F(x),F(x_0))<\veps$. Thus $F$ is continuous at $x_0$, and hence on $\calX$.

For $f\in C_b(\calX)$, the estimate $|f(F(x))|\leq\|f\|_\infty$ gives $\|\calK f\|_\infty\le\|f\|_\infty$. Since $\mc K\one_\calX=\one_\calX$, equality holds in the operator-norm bound.
\end{proof}

\begin{cor}\label{c:dense_in_Cb}
If $\calF$ is a linear function space which is Koopman invariant for $F$ and dense in $C_b(\calX)$, then $C_b(\calX)$ is Koopman invariant. In particular, $F$ is continuous.
\end{cor}
\begin{proof}
Let $f\in C_b(\calX)$. Then there exists $(f_n)\subset\calF$ such that $f_n\to f$ in $C_b(\calX)$. Hence, $\|f_n\circ F - f\circ F\|_\infty\le \|f_n - f\|_\infty\to 0$ as $n\to\infty$, which implies $f\circ F\in C_b(\calX)$.
\end{proof}

\subsubsection{Invariance on \texorpdfstring{$L^\infty(\mu)$}{L-infinity(mu)}}\label{s:invariance:det:Linf}
For a $\sigma$-finite Borel measure $\mu$ on $\calX$, consider the space of equivalence classes
\begin{align*}
L^\infty(\mu) := \left\{[f]\in[\mc M(\mc X)]: \mu\text{-}\esssup_{x\in\mc X} |f(x)|<\infty\right\}.
\end{align*}
Equipped with $\|[f]\|_{\infty}:=\mu\text{-}\esssup_{x\in\mc X}|f(x)|$, this is a Banach space. On $L^\infty(\mu)$, representative independence is the only obstruction to invariance; once it holds, the Koopman operator is a contraction.

\begin{thm}\label{t:det:invariance:Linfty} The following statements are equivalent for the deterministic system \eqref{eq:deterministic_system}: \begin{enumerate} \item[{\rm (i)}] One, and hence all, of the conditions in Proposition \rmref{p:equiv_classes} holds. \item[{\rm (ii)}] $L^\infty(\mu)\subset\wt\calD(\mu)$. \item[{\rm (iii)}] $L^\infty(\mu)$ is Koopman invariant.
\end{enumerate}
If one of {\rm (i)--(iii)} holds, then $\calK : L^\infty(\mu)\to L^\infty(\mu)$ is bounded with norm $\|\mc K\|_{L^\infty(\mu)\to L^\infty(\mu)}=1$.
\end{thm}

This is the deterministic specialization of Theorem~\ref{t:stoch:invariance:Linfty}, obtained by taking $\rho_x=\delta_{F(x)}$.

\subsubsection{Invariance on \texorpdfstring{$L^p(\mu)$}{Lp(mu)}, \texorpdfstring{$1\le p<\infty$}{1 <= p < infinity}} \label{s:invariance:det:Lp}

At finite exponents, representative independence no longer suffices: the transported mass must also be controlled. For $p\in[1,\infty)$, consider
\begin{align*}
L^p(\mu) := \left\{[f]\in[\mc M(\mc X)]:\int_{\mc X}|f|^p\,d\mu<\infty\right\}.
\end{align*}
With norm $\|[f]\|_{L^p(\mu)}:=\left(\int_{\mc X}|f|^p\,d\mu\right)^{1/p}$, this is a Banach space for every $p\in[1,\infty)$. Absolute continuity of the push-forward measure $\nu=\mu\circ F^{-1}$ no longer suffices for invariance. The obstruction is the possible unboundedness of its density, as Example~\ref{ex:invariance:det:Linfty} will show. The appropriate stronger condition is domination by $\mu$: there exists $C<\infty$ such that
\begin{equation}\label{e:measure_domination}
\nu(\Delta)\le C\mu(\Delta) \qquad(\Delta\in\calB_\calX).
\end{equation}
We abbreviate this inequality as $\nu\le C\mu$. By Lemma~\ref{l:measure_domination}, it is equivalent to $\nu\ll\mu$ with $d\nu/d\mu\in L^\infty(\mu)$.

This density condition characterizes invariance simultaneously for every finite $p$ and also determines the exact operator norm. The equivalence of (i)--(iv) is classical; see, for example, \cite[Theorem 2.1.1]{SinghManhas1993}. We include the proof.

\begin{thm}\label{t:det:invariance:Lp}
The following statements are equivalent for the deterministic system \eqref{eq:deterministic_system}: \begin{enumerate}[label={\rm(\roman*)}] \item For some $p\in[1,\infty)$, $L^p(\mu)$ is Koopman invariant, \item For all $p\in[1,\infty)$, $L^p(\mu)$ is Koopman invariant, \item For all $p\in[1,\infty)$, $L^p(\mu)$ is Koopman invariant, and $\mc K:L^p(\mu)\to L^p(\mu)$ is bounded, \item $\nu\le C\mu$ for some $C<\infty$.
\end{enumerate}
When these conditions hold, the optimal domination constant is
\begin{equation}\label{eq:Cstar}
C_*:=\inf\{C>0:\nu(\Delta)\le C\mu(\Delta) \text{ for every }\Delta\in\calB_\calX\}.
\end{equation}
Then $\nu$ is $\sigma$-finite, $g=d\nu/d\mu\in L^\infty(\mu)$, and
\begin{equation}\label{eq:det_exact_Lp_norm}
\|\calK\|_{L^p(\mu)\to L^p(\mu)} =C_*^{1/p}=\|g\|_{L^\infty(\mu)}^{1/p} \qquad(1\le p<\infty).
\end{equation}
If $\mc X\subset\R^n$ is open, $F:\mc X\to\mc X$ is a diffeomorphism and $\mu$ is the Lebesgue measure, then the equivalence extends to \begin{enumerate} \item[\rm{(v)}] $\inf_{x\in\mc X}|\det DF(x)|>0$.
\end{enumerate}
\end{thm}
\begin{proof} (i)$\Rightarrow$(iv). Suppose first that $L^p(\mu)$ is invariant for some $1\le p<\infty$. Theorem~\ref{t:schaefer_ec} makes the induced operator bounded. Let $M=\|\calK\|_{L^p\to L^p}$. If $\mu(\Delta)<\infty$, then
\begin{align*}
\nu(\Delta)=\int\one_\Delta\circ F\,d\mu =\|\calK\one_\Delta\|_p^p \le M^p\mu(\Delta).
\end{align*}
Here, $M>0$: if $M=0$, a countable finite-$\mu$ cover of $\calX$ would give $\nu(\calX)=0$, contrary to $\nu(\calX)=\mu(\calX)>0$. For $\mu(\Delta)=\infty$ the inequality is automatic, so (iv) follows.

(iv)$\Rightarrow$(iii). Conversely, if $\nu\le C\mu$, then a countable finite-$\mu$ cover is also a finite-$\nu$ cover. Thus $\nu$ is $\sigma$-finite, and the Radon--Nikodym theorem gives $d\nu=g\,d\mu$ with $g\in L^\infty(\mu)$. Moreover $C_*=\|g\|_\infty$ by Lemma~\ref{l:measure_domination}. Since $\nu\ll\mu$, the induced operator is well defined by Proposition~\ref{p:equiv_classes}. For every $1\le p<\infty$,
\begin{align}\label{e:p_iso}
\|\calK f\|_p^p &=\int|f\circ F|^p\,d\mu =\int|f|^p\,d\nu =\int|f|^p g\,d\mu \le\|g\|_\infty\|f\|_p^p.
\end{align}
Thus (iii) holds. The implications (iii)$\Rightarrow$(ii)$\Rightarrow$(i) are immediate.

It remains to identify the norm. Identity \eqref{e:p_iso} gives the upper bound. Let $a=\|g\|_\infty$. If $a=0$ there is nothing to prove. Given $0<\varepsilon<a$, the set $A_\varepsilon=\{g>a-\varepsilon\}$ has positive $\mu$-measure. By $\sigma$-finiteness there is a measurable $E_\varepsilon\subset A_\varepsilon$ with $0<\mu(E_\varepsilon)<\infty$. For $f_\varepsilon=\mu(E_\varepsilon)^{-1/p}\one_{E_\varepsilon}$,
\begin{align*}
\|f_\varepsilon\|_p=1,\qquad \|\calK f_\varepsilon\|_p^p =\int |f_\varepsilon|^p g\,d\mu \ge a-\varepsilon.
\end{align*}
Letting $\varepsilon\downarrow0$ proves the reverse inequality and hence \eqref{eq:det_exact_Lp_norm}.

(iv)$\Leftrightarrow$(v). Under the additional hypotheses, change of variables gives
\begin{align*}
g(y)=|\det D(F^{-1})(y)| =|\det DF(F^{-1}(y))|^{-1}.
\end{align*}
The continuous positive function $x\mapsto|\det DF(x)|$ has the same essential infimum as its infimum on the open set $\calX$: a strict pointwise inequality persists on a nonempty open set, which has positive Lebesgue measure. Hence $g\in L^\infty$ is equivalent to (v).
\end{proof}

The following example isolates the gap between absolute continuity and the bounded-density condition.

\begin{ex}\label{ex:invariance:det:Linfty} Let $\mc X=(0,1)$, let $F(x)=x^2$, and let $\mu$ be Lebesgue measure. The system $x^+=F(x)$ is the exact time-one discretization of $\dot x=\ln(2)x\ln(x)$. The push-forward measure in \eqref{e:nu_det} satisfies $\frac{d\nu}{d\mu}(x)=\frac1{2\sqrt x}$. Thus $\nu\ll\mu$, but $d\nu/d\mu\notin L^\infty(\mu)$. Accordingly, Theorem~\ref{t:det:invariance:Lp} rules out invariance of $L^p(\mu)$ for every $p\in[1,\infty)$. More explicitly, $f(x)=x^{-1/(2p)}$ belongs to $L^p(\mu)$, whereas $(\mc Kf)(x)=x^{-1/p}$ does not. By contrast, $L^\infty(\mu)$ is invariant by Theorem~\ref{t:det:invariance:Linfty}.
\end{ex}

\subsubsection{Invariance on reproducing kernel Hilbert spaces}\label{s:invariance:det:RKHS}
The preceding $L^p$ criteria depend only on the transported measure $\nu$. An RKHS retains pointwise geometry, and its invariance is instead encoded by an order relation between kernels. All Hilbert spaces in this paper are complex, with inner products linear in the first variable. Let $k:\calX\times\calX\to\mathbb C$ be a measurable Hermitian positive-definite kernel: $k(y,x)=\overline{k(x,y)}$ and, for all $N\in\N$, $c_1,\ldots,c_N\in\mathbb C$, and $x_1,\ldots,x_N\in\calX$,
\begin{align*}
\sum_{i,j=1}^N c_i\overline{c_j}\,k(x_i,x_j)\ge0.
\end{align*}
Let $\phi_x(y):=k(x,y)$. Its RKHS $\calH$ is the completion of $\operatorname{span}\{\phi_x:x\in\calX\}$ with
\begin{align*}
\langle\phi_x,\phi_y\rangle_\calH=k(x,y),\qquad f(x)=\langle f,\phi_x\rangle_\calH, \qquad \|\phi_x\|_\calH^2=k(x,x).
\end{align*}
We refer to \cite{pr} for these standard facts. The deterministic theory requires no global boundedness assumption on $k$. Finite linear combinations of the functions $\phi_x$ are Borel measurable, while convergence in $\calH$ implies pointwise convergence; hence every $f\in\calH$ is Borel measurable.

A general RKHS need not be a Banach lattice. For example, Gaussian RKHSs consist of analytic functions \cite{Minh10}, whereas taking absolute values does not preserve analyticity. We instead use the closed graph theorem. The maximal RKHS domain is
\begin{equation}\label{e:maxdomain}
\calD_\calH := \{f\in\calH : \calK f\in\calH\}.
\end{equation}
The Koopman operator on this domain is
\begin{align*}
\calK_\calH f := \calK f = f\circ F, \qquad f\in\calD_\calH.
\end{align*}
Continuity of point evaluation shows that this operator is closed; see also \cite[Lemma 1]{Gonzalez2024} and \cite[Lemma 3.1]{KohnPhil24}.

\begin{lem}\label{l:RKHS_invar_det} The operator $\calK_\calH$ is closed. In particular, if $\calH$ is Koopman invariant, equivalently $\calD_\calH=\calH$, then $\calK_\calH:\calH\to\calH$ is bounded.
\end{lem}
\begin{proof} Let $f_n\in\calD_\calH$, $f\in\calH$, and $g\in\calH$ satisfy $f_n\to f$ and $\calK_\calH f_n\to g$ in $\calH$. Continuity of point evaluation gives, for every $x\in\calX$,
\begin{align*}
g(x)=\lim_{n\to\infty}f_n(F(x))=f(F(x)).
\end{align*}
Thus $f\circ F=g\in\calH$, so $f\in\calD_\calH$ and $\calK_\calH f=g$. If $\calD_\calH=\calH$, the closed graph theorem gives boundedness.
\end{proof}
For positive-definite kernels $k_1$ and $k_2$ on $\calX\times\calX$, with associated RKHSs $\calH_1$ and $\calH_2$, let $k_1\preceq k_2$ mean that $k_2-k_1$ is positive definite. Aronszajn's inclusion theorem \cite[Theorem 5.1]{pr} gives the following equivalences: \begin{itemize} \item $\calH_1$ is continuously embedded in $\calH_2$; \item $\calH_1\subset\calH_2$; \item $k_1\preceq c^2k_2$ with some $c>0$.
\end{itemize}
By \cite[Theorem 5.7]{pr} (see also Lemma~\ref{l:Hrho} for the stochastic analog), the pullback kernel $k_F(x,y)=k(F(x),F(y))$ induces the RKHS
\begin{align}\label{e:H_F}
\calH_F = \{f\circ F : f\in\calH\} \quad\text{and}\quad \|g\|_{\calH_F} = \inf\{\|f\|_\calH : f\in\calH,\,g=f\circ F\},
\end{align}
Aronszajn's theorem therefore turns Koopman invariance into a comparison of kernels (see also \cite[Theorem 5.10]{pr}).

\begin{thm}\label{t:RKHS_det_inv} The following statements are equivalent: \begin{enumerate} \item[{\rm (i)}] $\calH$ is Koopman invariant; \item[{\rm (ii)}] $\calH_F\subset\calH$; \item[{\rm (iii)}] $\calH_F$ is continuously embedded in $\calH$; \item[{\rm (iv)}] $k_F \preceq c^2k$ for some $c>0$.
\end{enumerate}
\end{thm}
\begin{proof} By \eqref{e:H_F}, conditions (i) and (ii) are identical. Aronszajn's inclusion theorem, applied to the kernels $k_F$ and $k$, gives the equivalence of (ii), (iii), and (iv).
\end{proof}

\begin{ex} Let $\mathbb D$ be the unit disc in $\C$, and let $F:\mathbb D\to\mathbb D$ be analytic. Then Littlewood's inequalities for analytic self-maps imply Koopman invariance of the Hardy space $H^2(\mathbb D)$~\cite{Littlewood1925}.
\end{ex}

When $\calX\subset\mathbb R^n$ is a bounded Lipschitz domain, the RKHS $\calH$ corresponding to a Wendland or Mat\'ern kernel and the fractional Sobolev space $H^{\sigma(n,s)}(\calX)$ coincide as sets and have equivalent norms, where the Sobolev order $\sigma(n,s)\ge1$ depends on the kernel smoothness $s\in\N$ and the dimension $n\in\N$ (see \cite[Appendices B and C]{hertel2025koopman}):
\begin{equation}\label{eq:sigma:Wendland:Matern}
\sigma(n,s) = \begin{cases}
\frac{n+1}{2}+s & \qquad\text{Wendland kernels} \\
\frac n2+s & \qquad\text{Mat\'ern kernels}
\end{cases}
\end{equation}%
Sobolev composition estimates can therefore be used to verify the domain condition in Lemma~\ref{l:RKHS_invar_det}, as in \cite{KohnPhil24}.

\begin{thm}\label{t:sob_det} 
Let $\calX\subset\mathbb R^n$ be a bounded Lipschitz domain, let $m\in\N$ satisfy $m>n/2$, and let $F:\calX\to\calX$ be a $C^1$-diffeomorphism such that $F\in C_b^m(\calX;\mathbb R^n)$ and
\begin{align*}
\inf_{x\in\calX}|\det DF(x)|>0.
\end{align*}
Here, $C_b^m$ means that all derivatives of order at most $m$ are bounded. 
If, for some smoothness parameter $s\in\N$, the Wendland or Mat\'ern RKHS $\calH$ is norm-equivalent to $H^{\sigma(n,s)}(\calX)$ with
\begin{align*}
\frac n2<\sigma(n,s)\le m,
\end{align*}
where $\sigma(n,s)$ is defined by~\eqref{eq:sigma:Wendland:Matern}, then $\calD_\calH=\calH$ and $\calK_\calH:\calH\to\calH$ is bounded.
\end{thm}
\begin{proof}
Under the stated hypotheses, composition with $F$ defines a bounded linear operator from $H^{\sigma(n,s)}(\calX)$ into itself, see \cite[Theorem 4.2]{KohnPhil24}. The claim now follows as $\calH$ and $H^{\sigma(n,s)}(\calX)$ coincide as sets and their norms are equivalent.
\end{proof}

Koopman invariance is a central hypothesis in error analysis for kernel EDMD, see \cite{KohnPhil24,bold2025kernel}. However, depending on the RKHS, it may impose strong rigidity on the dynamics. In particular, invariance of a Gaussian RKHS over any subset of $\R^n$ with nonempty interior forces the dynamics to be affine, see Appendix~\ref{a:GaussianRKHS} and  \cite{Gonzalez2024,Ikeda22}. Analogous restrictions occur for polynomial kernels $k(x,y)=(1+x^\top y)^m$ on $\mathbb R^d$ ($m\in\N$), Fock kernels $k(z,w)=e^{\alpha\overline z w}$ on $\mathbb C$, exponential dot-product kernels $k(x,y)=e^{\alpha x^\top y}$ on $\mathbb R^d$, Paley--Wiener kernels $k(x,y)=\operatorname{sinc}(\alpha(x-y))$ on $\mathbb R$, and other kernels treated in \cite{Ikeda22}. Here $\alpha>0$, $\operatorname{sinc}(t)=\sin(t)/t$ for $t\ne0$, and $\operatorname{sinc}(0)=1$. We refer to \cite{Gonzalez2024} and the references therein for details. \subsection{Markov kernels}\label{s:invariance:stoch} For a Markov kernel $\rho$, the same sequence of spaces reveals precisely where averaging changes the deterministic picture. The Koopman operator is given by \eqref{eq:informal} on functions and by \eqref{eq:Koopman_equivalence_classes} on equivalence classes.

\subsubsection{Invariance on \texorpdfstring{$C_b(\calX)$}{Cb(X)}} \label{s:invariance:stoch:Cb} For deterministic systems, Theorem~\ref{t:det_invariance_Cb} identifies invariance of $C_b(\mc X)$ with continuity of $F$. For a Markov kernel, the corresponding condition is weak continuity: the map $x\mapsto\int f\,d\rho_x$ is continuous for every $f\in C_b(\mc X)$. This is also called the \emph{weak Feller property}; see \cite[Definition C.3]{hernandez-lerma_discrete-time_1996}. In the locally compact setting, compare \cite[Section 6.1.1]{meyn_markov_2012}. The Portmanteau theorem gives several equivalent formulations.

\begin{thm}\label{t:inv_C_b} The following statements are equivalent for a Markov kernel $\rho$: \begin{enumerate}[label=\rm(\roman*)] \item $C_b(\mc X)$ is Koopman invariant; \item The Markov kernel $\rho$ is weakly continuous, i.e., $\mc Kf\in C_b(\mc X)$ for all $f\in C_b(\mc X)$; \item For every bounded real-valued lower semicontinuous function $f$, $\mc Kf$ is lower semicontinuous; \item For every open set $\calO\subset\mc X$, the function $\mc K\one_{\calO}$ is lower semicontinuous; \item For every $x_0\in\mc X$ and $\Delta\in\calB_\calX$ satisfying $\rho_{x_0}(\partial\Delta)=0$, one has $\rho_x(\Delta)\to\rho_{x_0}(\Delta)$ as $x\to x_0$.
\end{enumerate}
If one of {\rm(i)--(v)} holds, then $\mc K:C_b(\mc X)\to C_b(\mc X)$ is bounded with norm $\|\mc K\|_{C_b(\mc X)\to C_b(\mc X)}=1$.
\end{thm}
\begin{proof} (i)$\Leftrightarrow$(ii). One implication is part of the definition of Koopman invariance. Conversely, a function $f\in C_b(\mc X)$ belongs to $L^1(\rho_x)$ for every $x$, because $f$ is bounded. Thus $\Delta_f=\mc X$, and weak continuity gives invariance.

(ii)$\Leftrightarrow$(iii)$\Leftrightarrow$(iv)$\Leftrightarrow$(v). For each $x_0$, condition (ii) says precisely that $\rho_x\Rightarrow\rho_{x_0}$ as $x\to x_0$. The four displayed formulations are the bounded-lower-semicontinuous, open-set, and continuity-set formulations of the Portmanteau theorem \cite[Theorem 2.1]{Bill99}; see also \cite[Proposition C.4]{hernandez-lerma_discrete-time_1996}.

Finally, under any of (i)--(v),
\begin{align*}
|(\mc Kf)(x)|\leq\int |f|\,d\rho_x\leq\|f\|_\infty \qquad(f\in C_b(\mc X),\ x\in\mc X).
\end{align*}
Hence $\|\mc K\|\le1$, while $\mc K\one_\calX=\one_\calX$ gives equality.
\end{proof}

\subsubsection{Invariance on \texorpdfstring{$L^\infty(\mu)$}{L-infinity(mu)}} \label{s:invariance:stoch:Linf}
In Proposition \rmref{p:equiv_classes}, we established several equivalent conditions for $\calK$ to be well-defined on equivalence classes (representative independence). As it turns out, representative independence is both necessary and sufficient for the Koopman invariance of $L^\infty(\mu)$.

\begin{thm}\label{t:stoch:invariance:Linfty} The following statements are equivalent for a Markov kernel $\rho$:
\begin{enumerate}
\item[{\rm (i)}] One, and hence all, of the conditions in Proposition \rmref{p:equiv_classes} holds.
\item[{\rm (ii)}] $L^\infty(\mu)\subset\wt\calD(\mu)$.
\item[{\rm (iii)}] $L^\infty(\mu)$ is Koopman invariant.
\end{enumerate}
If one of {\rm (i)--(iii)} holds, then $\calK : L^\infty(\mu)\to L^\infty(\mu)$ is bounded with norm $\|\mc K\|_{L^\infty(\mu)\to L^\infty(\mu)}=1$.
\end{thm}
\begin{proof}
(iii)$\Rightarrow$(ii) follows from Definition~\ref{def:Koopman-invariance:equivalence-class}. If (ii) holds, then $[0]\in\wt\calD(\mu)$, so Proposition~\ref{p:equiv_classes} gives (i).

(i)$\Rightarrow$(iii). Let $[f]\in L^\infty(\mu)$. There is a Borel set $\Delta$ with $\mu(\Delta^c)=0$ such that $|f(x)|\le\|f\|_{L^\infty}$ on $\Delta$. Condition (i) gives $\nu(\Delta^c)=0$, and hence $\rho_x(\Delta^c)=0$ for $\mu$-almost every $x$. For such $x$,
\begin{align*}
\int|f|\,d\rho_x =\int_\Delta|f|\,d\rho_x \le\|f\|_{L^\infty}\rho_x(\Delta) =\|f\|_{L^\infty}.
\end{align*}
Thus $f\in\calD(\mu)$ and, by Proposition~\ref{p:equiv_classes}, $[f]\in\wt\calD(\mu)$. The same estimate gives $\|\calK[f]\|_{L^\infty}\le\|[f]\|_{L^\infty}$, proving invariance and $\|\mc K\|\le1$. Since $\mc K\one_\calX=\one_\calX$, the norm is one.
\end{proof}

\subsubsection{Invariance on \texorpdfstring{$L^p(\mu)$}{Lp(mu)}, \texorpdfstring{$1\le p<\infty$}{1 <= p < infinity}} \label{s:invariance:stoch:Lp}

For finite $p$, stochastic averaging changes the picture. Deterministic invariance for one $L^p$ forces invariance for every finite exponent, whereas a Markov kernel may preserve $L^q$ without preserving any smaller $L^p$. Invariance nevertheless propagates toward larger exponents, with a quantitative norm bound.

\begin{thm}\label{t:invariance:stoch:Lp}
The following statements are equivalent for a Markov kernel $\rho$ and $q\in[1,\infty)$: \begin{enumerate}[label={\rm(\roman*)}] \item $L^q(\mu)$ is Koopman invariant, \item For all $p \in [q,\infty]$, $L^p(\mu)$ is Koopman invariant, \item For all $p \in [q,\infty]$, $L^p(\mu)$ is Koopman invariant, and $\calK : L^p(\mu)\to L^p(\mu)$ is bounded.
\end{enumerate}
Under these conditions, if $M_r:=\|\calK\|_{L^r(\mu)\to L^r(\mu)}$, then
\begin{equation}\label{eq:stoch_Lq_Lp_interpolation}
M_p\le M_q^{q/p}\quad(q\le p<\infty), \qquad M_\infty=1.
\end{equation}
If $q=1$, the equivalence extends to \begin{enumerate} \item[\rm (iv)] $\nu\le C\mu$ for some $C<\infty$.
\end{enumerate}
In this case, $\nu$ is $\sigma$-finite, $g=d\nu/d\mu\in L^\infty(\mu)$, and the $L^1$ norm is exact:
\begin{equation}\label{eq:stoch_exact_L1_norm}
\|\calK\|_{L^1(\mu)\to L^1(\mu)} =C_*=\|g\|_{L^\infty(\mu)},
\end{equation}
where
\begin{align*}
C_*:=\inf\{C>0:\nu(\Delta)\le C\mu(\Delta) \text{ for every }\Delta\in\calB_\calX\}
\end{align*}
is the best domination constant {\rm (}the same quantity as in \eqref{eq:Cstar}{\rm )}.
\end{thm}
\begin{proof} (i)$\Rightarrow$(iii). Invariance of $L^q(\mu)$ first implies $[0]\in\wt\calD(\mu)$, hence $\nu\ll\mu$ and $[\calD(\mu)]=\wt\calD(\mu)$ by Proposition~\ref{p:equiv_classes}. Theorem~\ref{t:schaefer_ec} gives $M_q<\infty$, and Theorem~\ref{t:stoch:invariance:Linfty} gives invariance on $L^\infty$ and $M_\infty=1$.

Fix $q\le p<\infty$, let $r=p/q\ge1$, and take $f\in L^p(\mu)$. Then $|f|^r\in L^q(\mu)$. For almost every $x$, the function $|f|^r$ is $\rho_x$-integrable; because $\rho_x$ is a probability measure, this also implies $f\in L^1(\rho_x)$. Jensen's inequality yields
\begin{align*}
|\calK f|^r\le\calK(|f|^r)\qquad\mu\text{-a.e.}
\end{align*}
(with the case $r=1$ included). Therefore
\begin{align*}
\|\calK f\|_p^p \le\|\calK(|f|^r)\|_q^q \le M_q^q\||f|^r\|_q^q =M_q^q\|f\|_p^p.
\end{align*}
Hence $L^p(\mu)$ is invariant and $M_p\le M_q^{q/p}$. This proves (iii) and \eqref{eq:stoch_Lq_Lp_interpolation}. The implications (iii)$\Rightarrow$(ii)$\Rightarrow$(i) are immediate.

It remains to prove the additional equivalence when $q=1$. Assume first that (i) holds, and let $M_1=\|\calK\|_{L^1\to L^1}$. For every Borel $\Delta$ of finite $\mu$-measure, Tonelli's theorem gives
\begin{align*}
\nu(\Delta)=\int\calK\one_\Delta\,d\mu \le M_1\mu(\Delta).
\end{align*}
Moreover $M_1>0$: otherwise a countable finite-$\mu$ cover of $\calX$ would imply $\nu(\calX)=0$, contradicting $\nu(\calX)=\mu(\calX)>0$. The domination is automatic when $\mu(\Delta)=\infty$, so (iv) follows.

Conversely, assume (iv). Domination by $C\mu$ makes $\nu$ $\sigma$-finite, and the Radon--Nikodym theorem gives $d\nu=g\,d\mu$ with $g\in L^\infty(\mu)$. For $f\in L^1(\mu)$, Tonelli's theorem gives
\begin{equation}\label{e:proof:t:stoch:invariance:Lp}
\int\!\int |f(y)|\,d\rho_x(y)d\mu(x) =\int|f|\,d\nu =\int|f|g\,d\mu<\infty.
\end{equation}
Hence $f\in\calD(\mu)$, $[f]\in\wt\calD(\mu)$, and $\calK[f]\in L^1(\mu)$, proving (i).

To obtain the exact norm, take a nonnegative $f\in L^1(\mu)$. Then
\begin{align*}
\|\calK f\|_1=\int\calK f\,d\mu=\int f g\,d\mu.
\end{align*}
The preceding estimate bounds the norm of $\calK$ from above by $\|g\|_\infty$. For the reverse bound, let $a=\|g\|_\infty$. If $a>0$, then for every $0<\varepsilon<a$, $\sigma$-finiteness supplies a measurable set $E\subset\{g>a-\varepsilon\}$ with $0<\mu(E)<\infty$. Taking $f=\mu(E)^{-1}\one_E$ gives $\|f\|_1=1$ and $\|\calK f\|_1=\int fg\,d\mu\ge a-\varepsilon$. Hence $\|\calK\|=\|g\|_\infty$. Lemma~\ref{l:measure_domination} identifies this number with $C_*$, proving \eqref{eq:stoch_exact_L1_norm}.
\end{proof}

\begin{cor}\label{c:Lp_inv}
If $\nu\le C\mu$ for some $C>0$, then every $L^p(\mu)$, $1\le p\le\infty$, is Koopman invariant and
\begin{align*}
\|\calK\|_{L^p\to L^p}\le C_*^{1/p}\qquad(1\le p<\infty), \qquad \|\calK\|_{L^\infty\to L^\infty}=1.
\end{align*}
In particular, if $\mu$ is invariant, then $C_*=1$, and $\calK$ is a contraction on every $L^p(\mu)$.
\end{cor}

Propagation toward larger exponents cannot in general be reversed: the following genuinely stochastic example preserves $L^2(\mu)$ but not $L^1(\mu)$.

\begin{ex} Let $\calX=(0,1)$ and let $\mu$ be Lebesgue measure. Consider the transition kernel
\begin{align*}
d\rho_x(y)=h(y)\,dy, \qquad h(y)=\frac34y^{-1/4}.
\end{align*}
The normalization $\int_0^1h(y)\,dy=1$ makes $\rho_x$ a probability measure. Since it is independent of $x$ and $\mu(\calX)=1$, one has $\nu=\rho_x$ and $d\nu/d\mu=h$. For $f\in L^2(\mu)$,
\begin{align*}
\calK f=\int_0^1 f(y)h(y)\,dy, \qquad \|\calK f\|_2\le\|h\|_2\|f\|_2.
\end{align*}
Thus $h\in L^2(\mu)$ implies that $L^2(\mu)$ is Koopman invariant. By contrast, $h\notin L^\infty(\mu)$, so Theorem~\ref{t:invariance:stoch:Lp} shows that $L^1(\mu)$ is not Koopman invariant.
\end{ex}

If $\mu$ is finite, $L^\infty(\mu)$ is Koopman invariant if and only if $\frac{d\nu}{d\mu}\in L^1(\mu)$, whereas $L^1(\mu)$ is Koopman invariant if and only if $\frac{d\nu}{d\mu}\in L^\infty(\mu)$. One may therefore be inclined to conjecture that $L^p(\mu)$ is Koopman invariant if and only if $\frac{d\nu}{d\mu}\in L^q(\mu)$, where $\frac 1p + \frac 1q = 1$. However, as the following proposition shows, the latter condition is in fact equivalent to an exact two-space criterion, including both endpoints.

\begin{prop}\label{p:stoch:Lp_to_L1}
Let $1\le p\le\infty$, and let $q$ be its conjugate exponent. The induced Koopman operator is well defined and bounded from $L^p(\mu)$ to $L^1(\mu)$ if and only if $d\nu = g\,d\mu$ with $g\in L^q(\mu)$.
In that event,
\begin{equation}\label{eq:exact_Lp_to_L1_norm}
\|\calK\|_{L^p(\mu)\to L^1(\mu)}=\|g\|_{L^q(\mu)}.
\end{equation}
\end{prop}
\begin{proof} Suppose first that $g=d\nu/d\mu$ exists and belongs to $L^q(\mu)$. For $f\in L^p(\mu)$, Tonelli's theorem and H\"older's inequality give
\begin{align*}
\int\!\int |f(y)|\,d\rho_x(y)d\mu(x) = \int|f|\,d\nu = \int|f|g\,d\mu \le\|f\|_p\|g\|_q.
\end{align*}
Thus $f\in\calD(\mu)$, and $\nu\ll\mu$ implies $[f]\in\wt\calD(\mu)$. Moreover, $\|\calK f\|_1\le\|f\|_p\|g\|_q$, so the induced map is bounded.

Conversely, suppose that the induced map is bounded with norm $M$. Then $L^p(\mu)\subset\wt\calD(\mu)$, so $\nu\ll\mu$ by Proposition~\ref{p:equiv_classes}. Choose $E_n\uparrow\calX$ with $\mu(E_n)<\infty$. Since
\begin{align*}
\nu(E_n)=\|\calK\one_{E_n}\|_1 \le M\|\one_{E_n}\|_p<\infty
\end{align*}
(with $\|\one_{E_n}\|_\infty\le1$), $\nu$ is $\sigma$-finite. Hence $g=d\nu/d\mu$ exists by Radon--Nikodym. For every nonnegative $f\in L^p(\mu)$, Tonelli gives
\begin{align*}
\int fg\,d\mu=\int f\,d\nu =\int\calK f\,d\mu=\|\calK f\|_1 \le M\|f\|_p.
\end{align*}
If $1<p<\infty$, $L^p$--$L^q$ duality yields $g\in L^q(\mu)$ and $\|g\|_q\le M$. If $p=\infty$, take $f=\one_\calX$ in the preceding inequality to obtain $g\in L^1(\mu)$ and $\|g\|_1\le M$. Finally, let $p=1$. If $\mu(\{g>M+\varepsilon\})>0$ for some $\varepsilon>0$, then $\sigma$-finiteness supplies a measurable subset $E$ of this set with $0<\mu(E)<\infty$. Taking $f=\one_E$ gives
\begin{align*}
\int_E g\,d\mu\le M\mu(E),
\end{align*}
contrary to the definition of $E$. Thus $g\in L^\infty(\mu)$ and $\|g\|_\infty\le M$. Together with the upper bound in the first paragraph, these three cases prove \eqref{eq:exact_Lp_to_L1_norm}.
\end{proof}

Returning to invariant spaces, stationarity yields simultaneous invariance throughout the $L^p$ scale. Overdamped Langevin dynamics provide a model case that will recur in the sections on compactness and self-adjointness, below.

\begin{ex}[Overdamped Langevin dynamics]\label{ex:olan}
Let $\calX=\mathbb R^n$, let $D,\beta>0$, and let $V\in C^2(\mathbb R^n)$ have globally Lipschitz gradient. Moreover, suppose that
\begin{align*}
Z := \int_{\R^n}e^{-\beta V(x)}\,dx<\infty,
\end{align*}
and define the Gibbs probability measure
\begin{align*}
d\mu(x)=Z^{-1}e^{-\beta V(x)}\,dx.
\end{align*}
Consider the stochastic differential equation (SDE)
\begin{align}\label{eq:langevin}
dX_t=-D\nabla V(X_t)\,dt +\sqrt{\frac{2D}{\beta}}\,dW_t,\qquad X_0=x,
\end{align}
where $(W_t)_{t\ge0}$ is standard $n$-dimensional Brownian motion. Since $\nabla V$ is globally Lipschitz, it satisfies the linear growth condition $\|\nabla V(x)\|\le \|\nabla V(0)\| + L\|x\|$ for some $L>0$. Hence, by \cite[Theorem 5.2.1]{Oksendal2000}, the SDE \eqref{eq:langevin} admits a unique global strong solution $(X_t^x)_{t\ge0}$ for every $x\in\mathbb R^n$. Its transition probabilities
\begin{align*}
\rho_t(x,\Delta):=\mathbb P(X_t^x\in\Delta), \qquad t\ge0,
\end{align*}
form Markov kernels.

Consider the differential operator $\mathcal L_0$ associated with the diffusion, defined on $C_c^2(\mathbb R^n)$,
\begin{align*}
\mathcal L_0 f =\frac{D}{\beta}\Delta f-D\nabla V\cdot\nabla f =\frac{D}{\beta}e^{\beta V}\operatorname{div}\bigl(e^{-\beta V}\nabla f\bigr).
\end{align*}
For every $f\in C_c^2(\mathbb R^n)$, integration by parts gives
\begin{align*}
\int\mathcal L_0 f\,d\mu =\frac{D}{\beta Z}\int\operatorname{div}\bigl(e^{-\beta V}\nabla f\bigr)\,dx = 0.
\end{align*}
As the coefficients of \eqref{eq:langevin} are globally Lipschitz, the martingale problem for $(\mathcal L_0,C_c^\infty(\mathbb R^n))$ is well posed. Thus, Echeverr\'ia's invariant-measure criterion \cite{Echeverria1982} implies that $\mu$ is invariant, i.e.,
\begin{align*}
\int\rho_t(x,\Delta)\,d\mu(x)=\mu(\Delta) \qquad(\Delta\in\calB_{\mathbb R^n},\ t\ge0).
\end{align*}
Hence, by Corollary \ref{c:Lp_inv}, the associated Koopman operators $\calK_t : L^p(\mu)\to L^p(\mu)$ are well-defined contractions for every $p\in [1,\infty)$. The Markov property yields $\mathcal K_{t+s}=\mathcal K_t\mathcal K_s$ and $\mathcal K_0=I$. Finally, the continuity of the sample paths implies that, for every $f\in C_c^\infty(\mathbb R^n)$,
\begin{align*}
\lim_{t\downarrow0}\mathcal K_t f(x)=f(x)
\qquad(x\in\mathbb R^n).
\end{align*}
Since $|\mathcal K_t f(x)-f(x)|\le2\|f\|_\infty$, dominated convergence gives
\begin{align*}
\lim_{t\downarrow0}
\|\mathcal K_t f-f\|_{L^p(\mu)}=0.
\end{align*}
As $C_c^\infty(\mathbb R^n)$ is dense in $L^p(\mu)$ and the operators $\mathcal K_t$ are contractions, this convergence extends to every $f\in L^p(\mu)$. Thus, $(\mathcal K_t)_{t\ge0}$ is a strongly continuous semigroup of positive contractions on $L^p(\mu)$.

Specifically for $p=2$, let $\mathcal L$ denote the infinitesimal generator of the contraction semigroup $(\calK_t)_{t\ge0}$ on $L^2(\mu)$. By It\^o's formula \cite[Theorem 4.2.1]{Oksendal2000} and strong continuity, $C_c^2(\mathbb R^n)\subset D(\mathcal L)$, and $\mathcal Lf = \mathcal L_0f$ for $f\in C_c^2(\R^n)$. Define the unitary operator $U : L^2(\mu)\to L^2(\mathbb R^n)$ by $Uf=Z^{-1/2}e^{-\beta V/2}f$ and the Schr\"odinger operator $H_0 = -\frac D\beta\Delta+W$ on $C_c^\infty(\R^n)$, where
\begin{align*}
W(x) = \frac{D\beta}{4}\|\nabla V(x)\|^2 - \frac D2\Delta V(x).
\end{align*}
A direct computation shows that
\[
H_0 = -U\mathcal L_0 U^{-1}|_{C_c^\infty(\R^n)}\,\subset\,-U\mathcal LU^{-1}.
\]
By the Lumer--Phillips theorem, $-\mathcal L$ is maximally accretive, and hence so is $-U\mathcal LU^{-1}$. Since $\nabla V$ is globally Lipschitz, $W$ is continuous and bounded below. Hence, $H_0$ is essentially self-adjoint \cite[Theorem X.28]{ReedSimon1975}. Moreover, integration by parts shows that $H_0$ is nonnegative. Thus, its self-adjoint closure $H$ is maximally accretive. But as $H\subset -U\mathcal LU^{-1}$, we conclude that
\begin{align*}
H = -U\mathcal LU^{-1},\qquad\calK_t=U^{-1}e^{-tH}U.
\end{align*}
In particular, both $\mathcal L$ and all $\calK_t$, $t\ge 0$, are self-adjoint operators.
\end{ex}

\subsubsection{Invariance on reproducing kernel Hilbert spaces}\label{s:invariance:stoch:RKHS}
The one-step measure $\nu$ governs the preceding $L^p$ results. On an RKHS, the averaged kernel features carry the corresponding information. We use the conventions of Subsection~\ref{s:invariance:det:RKHS}; for the stochastic RKHS results in this subsection, we impose the additional bounded-diagonal hypothesis
\begin{equation}\label{eq:kappa_bounded_diagonal}
\kappa:=\sup_{y\in\calX}\sqrt{k(y,y)}<\infty.
\end{equation}
Then $|f(y)|\le\kappa\|f\|_\calH$ for every $f\in\calH$, and the kernel Cauchy--Schwarz inequality gives $|k(x,y)|\le\kappa^2$. Thus $\calK f$ is a bounded Borel function for every $f\in\calH$, with $\|\calK f\|_\infty\le\kappa\|f\|_\calH$.

With the maximal domain $\calD_\calH:=\{f\in\calH:\calK f\in\calH\}$, the stochastic Koopman operator is
\begin{align*}
\calK_\calH : \calH\supset\calD_\calH\to\calH, \qquad (\calK_\calH f)(x) := (\calK f)(x) = \int f\,d\rho_x.
\end{align*}
The pointwise bound above replaces the composition argument in the deterministic case and again yields closedness, see also \cite[Proposition 2.4]{hertel2025koopman}.

\begin{lem}\label{l:RKHS_invar_stochastic} The operator $\calK_\calH$ is closed. In particular, if $\calH$ is Koopman invariant, i.e., $\calD_\calH=\calH$, then $\calK_\calH:\calH\to\calH$ is bounded.
\end{lem}
\begin{proof} Let $f_n\in\calD_\calH$, $f\in\calH$, and $g\in\calH$ satisfy $f_n\to f$ and $\calK_\calH f_n\to g$ in $\calH$. For every $x$,
\begin{align*}
|\calK f_n(x)-\calK f(x)| \le\int|f_n-f|\,d\rho_x \le\kappa\|f_n-f\|_\calH\longrightarrow0.
\end{align*}
Continuity of point evaluation also gives $\calK_\calH f_n(x)\to g(x)$. Hence $g(x)=\calK f(x)$ for every $x$, so $\calK f=g\in\calH$ and $f\in\calD_\calH$. If $\calD_\calH=\calH$, the closed graph theorem gives boundedness.
\end{proof}

The bounded-diagonal hypothesis also allows each stochastic evaluation functional to be represented in $\calH$. For each $x\in\calX$, the linear functional
\begin{align*}
\Lambda_x(f):=\int f(y)\,d\rho_x(y),\qquad f\in\calH,
\end{align*}
satisfies $|\Lambda_x(f)|\le\kappa\|f\|_\calH$. The Riesz representation theorem therefore gives a unique $\Psi_x\in\calH$ such that
\begin{equation}\label{eq:Psi_Riesz}
\int f(y)\,d\rho_x(y)=\langle f,\Psi_x\rangle_\calH \quad(f\in\calH), \qquad \|\Psi_x\|_\calH\le\kappa.
\end{equation}
Thus, $\Psi_x$ is the Pettis (or weak) integral of $y\mapsto\phi_y$ against $\rho_x$.

The functions $\Psi_x$ define the Hermitian positive-definite kernel
\begin{align}\label{e:krho}
k_\rho(x,y) :=\langle\Psi_x,\Psi_y\rangle_\calH =\int\!\int k(u,v)\,d\rho_x(u)\,d\rho_y(v).
\end{align}
Indeed, the reproducing property and \eqref{eq:Psi_Riesz} give
\begin{align*}
\Psi_x(v)=\int k(u,v)\,d\rho_x(u),
\end{align*}
and a second application of \eqref{eq:Psi_Riesz}, followed by Fubini's theorem, gives the displayed formula. The double integral is absolutely convergent by \eqref{eq:kappa_bounded_diagonal}. Let $\calH_\rho$ be the RKHS induced by $k_\rho$. The range-space construction in Lemma~\ref{l:Hrho} gives
\begin{align}\label{e:H_rho}
\calH_\rho = \{\calK f : f\in\calH\} \quad\text{and}\quad \|g\|_{\calH_\rho} =\inf\{\|f\|_\calH:f\in\calH,\ g=\calK f\}.
\end{align}
As in the deterministic case, invariance is therefore equivalent to comparing $k_\rho$ with the original kernel $k$.

\begin{thm}\label{t:RKHS_inv} The following statements are equivalent: \begin{enumerate} \item[{\rm (i)}] $\calH$ is Koopman invariant; \item[{\rm (ii)}] $\calH_\rho\subset\calH$; \item[{\rm (iii)}] $\calH_\rho$ is continuously embedded in $\calH$; \item[{\rm (iv)}] $k_\rho\preceq c^2k$ for some $c>0$.
\end{enumerate}
\end{thm}
\begin{proof} By \eqref{e:H_rho}, conditions (i) and (ii) are identical. Aronszajn's inclusion theorem, applied to $k_\rho$ and $k$, gives the equivalence of (ii), (iii), and (iv).
\end{proof}

For Sobolev RKHSs associated with Wendland and Mat\'ern kernels, regularity of the transition density gives a verifiable %
sufficient condition. The following result establishes an $L^2$-to-$H^\ell$ smoothing property and derives its implications for RKHSs of continuous functions. It strengthens \cite[Theorem 2.7]{hertel2025koopman}, which gives boundedness of $\calK:H^s(\Omega)\to H^s(\Omega)$ for $s\in(0,\ell]$ under the same condition on~$p$. The smoothing property itself is the case $\calH=H^\ell(\Omega)$ of \cite[Theorem 3.3(a)]{hertel2026verifiable}, which even shows that $\calK:L^2(\Omega)\to H^\ell(\Omega)$ is in fact Hilbert--Schmidt.

\begin{thm}\label{t:sob_stoch} 
Let $\calX=\Omega\subset\mathbb R^n$ be a bounded Lipschitz domain, and let $\ell\in\N$. Assume $d\rho_x(y)=p(x,y)\,dy$, where $p:\Omega\times\Omega\to[0,\infty)$ is measurable, and its weak derivatives in the first variable satisfy
\begin{equation}\label{eq:stoch_sob_density}
D_1^\alpha p\in L^2(\Omega\times\Omega) \qquad(|\alpha|\le\ell).
\end{equation}
Then the %
Koopman operator is well defined on equivalence classes, defines a bounded operator
\begin{align*}
\calK:L^2(\Omega)\to H^\ell(\Omega)
\end{align*}
and, consequently, a bounded operator $\calK : H^r(\Omega)\to H^r(\Omega)$ for every $0\le r\le\ell$, with all mappings understood in the sense of equivalence classes.

For the pointwise RKHS conclusion, suppose in addition that $\ell>n/2$ and that the specified map $x\mapsto p(x,\cdot)$ agrees at every $x$ with the continuous $L^2(\Omega)$-valued representative of $p\in H^\ell(\Omega;L^2(\Omega))$. 
Then every Wendland or Mat\'ern RKHS $\calH$ that, for some smoothness parameter $s\in\N$, is norm-equivalent to $H^{\sigma(n,s)}(\Omega)$, as a space of its continuous representatives, with $\sigma(n,s)$ defined by~\eqref{eq:sigma:Wendland:Matern} and satisfying
\begin{align*}
\frac n2<\sigma(n,s)\le\ell
\end{align*}
is Koopman invariant, and $\calK_\calH:\calH\to\calH$ is bounded.
\end{thm}
\begin{proof} 
Let $\mu$ be the Lebesgue measure on $\Omega$. Since every $\rho_x$ has a density, $\nu(\Delta)=\int_\Omega\rho_x(\Delta)\,dx=0$ holds whenever $\mu(\Delta)=0$. Hence $\nu\ll\mu$, and the Koopman rule is independent of the representative by Proposition~\ref{p:equiv_classes}. Condition \eqref{eq:stoch_sob_density} is precisely the assertion that $p\in H^\ell(\Omega;L^2(\Omega))$, where $L^2(\Omega\times\Omega)$ is identified with $L^2(\Omega;L^2(\Omega))$ via $x\mapsto p(x,\cdot)$. Let $f\in L^2(\Omega)$. Testing the distributional derivative of $x\mapsto\int_\Omega p(x,y)f(y)\,dy$ against a function in $C_c^\infty(\Omega)$ and applying Fubini's theorem, which is justified since $p,D_1^\alpha p\in L^2(\Omega\times\Omega)$ and $f\in L^2(\Omega)$, gives
\begin{align*}
D^\alpha(\calK f)(x) =\int_\Omega D_1^\alpha p(x,y)f(y)\,dy, \qquad |\alpha|\le\ell.
\end{align*}
Cauchy--Schwarz therefore yields
\begin{align*}
\|D^\alpha(\calK f)\|_{L^2(\Omega)} \le\|D_1^\alpha p\|_{L^2(\Omega\times\Omega)}\|f\|_{L^2(\Omega)}.
\end{align*}
Summing over $|\alpha|\le\ell$ proves the smoothing property. The embeddings $H^r(\Omega)\hookrightarrow L^2(\Omega)$ and $H^\ell(\Omega)\hookrightarrow H^r(\Omega)$, valid for $0\le r\le\ell$, prove the equivalence-class assertions. Under the additional compatibility hypothesis, $p(x,\cdot)\in L^2(\Omega)$ for every $x$, so the pointwise integral defining $\calK f(x)$ exists for every $f\in L^2(\Omega)$ and every $x$. For such $f$, this pointwise function agrees almost everywhere with the $H^\ell$ element just constructed. Both are continuous when $\ell>n/2$ (the first by continuity of $x\mapsto p(x,\cdot)$ in $L^2$, the second by Sobolev embedding), so they agree everywhere. Taking $r=\sigma(n,s)$ in the two continuous embeddings proves pointwise Koopman invariance of the stated RKHSs. Since the norms of $\calH$ and $H^{\sigma(n,s)}(\Omega)$ are equivalent, the smoothing estimate and the same two embeddings give $\|\calK_\calH f\|_\calH\le C\|f\|_\calH$ for all $f\in\calH$, with $C$ independent of~$f$. Hence $\calK_\calH:\calH\to\calH$ is bounded.
\end{proof}

The deterministic and stochastic RKHS criteria have the same range-space form but arise through different mechanisms: composition can force geometric rigidity, whereas averaging can create Sobolev regularity. Koopman invariance of these RKHSs also enters error estimates for kernel approximations of the Koopman operator; see \cite{hertel2025koopman}.

\section{Compactness of the Koopman operator}
\label{sec:compactness} 

The preceding section established %
boundedness once the chosen space is Koopman invariant. %
In this section, we investigate compactness of the Koopman operator. %
Recall that a bounded operator between Banach spaces is \emph{compact} if it maps the closed unit ball to a relatively compact set. For a compact operator on a complex infinite-dimensional Banach space, every nonzero spectral value is an isolated eigenvalue of finite algebraic multiplicity, and zero is the only possible accumulation point of the spectrum.

In this section, we characterize the compactness of the Koopman operator on the function spaces already considered in Section \ref{sec:invariance}. Let us briefly outline the main results: %
A deterministic composition operator on $C_b(\calX)$ with continuous $F$ will turn out to be compact if and only if $F(\calX)$ is finite and, under $\nu\ll\mu$, on $L^\infty(\mu)$ if and only if $F$ has essentially finite range (Theorems~\ref{t:C_det_compact} and \ref{t:Linfty_comp_det}). On $L^p(\mu)$, $1\le p<\infty$, compactness is instead governed by the nonatomic part of $\mu$ together with a decay condition on the transported atomic masses (Theorem~\ref{t:compact_Lp_det}). For stochastic Koopman operators for which $C_b(\calX)$ is Koopman invariant, compactness on $C_b(\calX)$ is equivalent to total boundedness of the set of the transition laws $\rho_x$, implying total-variation continuity of $x\mapsto\rho_x$, which is also sufficient for compact $\calX$ (Theorem~\ref{t:Cb_compact}). On $L^\infty(\mu)$, compactness is equivalent to total boundedness of $\{\rho_x:x\in\calX\backslash N\}$ for some $\mu$-null set $N$ (Theorem~\ref{t:compact:stoch:Linf}). %
Under domination, that is, $\nu(\Delta)\le C\mu(\Delta)$ for all $\Delta\in\calB_\calX$, we will give an exact $L^1$-criterion of Kolmogorov--Riesz type, based on uniform integrability, tightness, and local mean oscillation, from which compactness propagates to every $L^p(\mu)$, $1\le p<\infty$ (Lemma~\ref{l:compact1:Lp} and Theorem~\ref{t:compact2:Lp}). On an RKHS, both regimes are captured by compactness of the canonical inclusion of the range space, which we will characterize by uniform eigenvalue decay in finite Gram-matrix comparisons (Theorems~\ref{t:RKHS_comp_det} and \ref{t:RKHS_comp}). 
For a translation-invariant kernel satisfying the strict-separation hypothesis of Proposition~\ref{p:nec}(i) and deterministic systems, compactness forces $F(\calX)$ to be bounded.

\subsection{Deterministic systems}

\subsubsection{Compactness on \texorpdfstring{$C_b(\calX)$}{Cb(X)}}
Here, compactness detects the size of the image of $F$ exactly; no compactness assumption on $\calX$ is needed.

\begin{thm}[{\cite[Corollary 2.2]{SinghSummers1987}}]\label{t:C_det_compact}
Let $F:\calX\to\calX$ be continuous. Then the Koopman operator $\calK : C_b(\calX)\to C_b(\calX)$, $\calK f = f\circ F$, is compact if and only if $F(\calX)$ is finite. In particular, if $\calK$ is compact, then $F$ is constant on any connected component of $\calX$.
\end{thm}
\begin{proof} If $F(\calX)=\{y_1,\ldots,y_m\}$, then every $\calK f$ is constant on each fiber $F^{-1}(\{y_j\})$. Thus $\calK$ has rank at most $m$.

Conversely, suppose that $F(\calX)$ is infinite and choose $x_n\in\calX$ so that the points $F(x_n)$ are pairwise distinct. Evaluation at $x$ defines $\delta_x\in C_b(\calX)^*$ with norm one, and
\begin{align*}
\calK^*\delta_x=\delta_{F(x)}.
\end{align*}
For distinct $y,z\in\calX$, the metric-space Urysohn lemma gives $f\in C_b(\calX)$ with $\|f\|_\infty\leq1$, $f(y)=1$, and $f(z)=-1$. Consequently,
\begin{align*}
\|\delta_y-\delta_z\|_{C_b(\calX)^*}=2.
\end{align*}
The sequence $(\calK^*\delta_{x_n})_n$ therefore has no convergent subsequence. Hence, $\calK^*$ cannot be compact. By Schauder's theorem, the same holds for $\calK$.
\end{proof}

\subsubsection{Compactness on \texorpdfstring{$L^\infty$}{L-infinity}}

Passing from pointwise functions to equivalence classes changes ``finite'' to ``essentially finite''. Recall that $\nu=\mu\circ F^{-1}\ll\mu$ is precisely the condition that composition by $F$ be well defined on $L^\infty(\mu)$.

\begin{thm}\label{t:Linfty_comp_det} Assume that $\nu\ll\mu$. Then
\begin{align*}
\calK:L^\infty(\mu)\longrightarrow L^\infty(\mu),\qquad \calK f=f\circ F,
\end{align*}
is compact if and only if there is a measurable $\mu$-null set $N$ such that $F(\calX\backslash N)$ is finite.
\end{thm}

\begin{proof} This is the point-mass specialization of the general $L^\infty$ criterion for Markov kernels, Theorem~\ref{t:compact:stoch:Linf}. Indeed, $\rho_x=\delta_{F(x)}$, and
\begin{align*}
\|\delta_y-\delta_z\|_{TV}= \begin{cases} 0,&y=z,\\ 2,&y\neq z.
\end{cases}
\end{align*}
A family of such Dirac masses is therefore totally bounded in total variation precisely when it is finite. Condition~{\rm(iii)} of Theorem~\ref{t:compact:stoch:Linf} now gives the assertion.
\end{proof}

\subsubsection{Compactness on \texorpdfstring{$L^p(\mu)$}{Lp(mu)}, \texorpdfstring{$1\leq p<\infty$}{1 <= p < infinity}}

The obstruction on $L^p$ is measure-theoretic rather than topological. Throughout this subsection, assume that $\mu$ is $\sigma$-finite. The set $\calX_{\mathrm a} = \{x\in\calX : \mu(\{x\})>0\}$ of points of positive $\mu$-mass is at most countable. Thus, for some finite or countable index set $I$,
\begin{align*}
\calX_{\mathrm a}=\{x_i:i\in I\},\qquad \mu_i=\mu(\{x_i\})>0,\qquad \calX_{\mathrm{na}}=\calX\backslash\calX_{\mathrm a}.
\end{align*}
The restriction of $\mu$ to $\calX_{\mathrm{na}}$ is nonatomic. Indeed, every $\mu$-atom has finite measure by $\sigma$-finiteness. If $A$ is such an atom, then the normalized measure $\mu(\,\cdot\,\cap A)/\mu(A)$ is a $\{0,1\}$-valued probability measure on the standard Borel space $\calX$, and is therefore a Dirac measure. Thus $A$ agrees modulo a null set with a singleton of positive $\mu$-measure.

Let $\nu=\mu\circ F^{-1}$ and $\nu_i=\nu(\{x_i\})$. If $\nu\ll\mu$, then $\mu(\{F(x_i)\})>0$ for every $i$. Otherwise the $\mu$-null singleton $\{F(x_i)\}$ would have a preimage of measure at least $\mu_i$. Consequently, there is a map $\pi_F:I\to I$ such that $F(x_i)=x_{\pi_F(i)}$, and
\begin{align}\label{e:numu_j}
\nu_j = \mu\bigl(F^{-1}(\{x_j\})\cap\calX_{\mathrm a}\bigr) +\mu\bigl(F^{-1}(\{x_j\})\cap\calX_{\mathrm{na}}\bigr) = \sum_{i\in\pi_F^{-1}(j)}\mu_i +\mu\bigl(F^{-1}(\{x_j\})\cap\calX_{\mathrm{na}}\bigr).
\end{align}
Note that $\mu_i < \infty$ for all $i\in I$ as $\mu$ is $\sigma$-finite. Moreover, if $\calK : L^p(\mu)\to L^p(\mu)$ is well defined for some $p\in [1,\infty)$, then also $\nu_i < \infty$ for all $i\in I$. Indeed, $\mu_i < \infty$ implies $\one_{\{x_i\}}\in L^p(\mu)$, and hence $\one_{F^{-1}(\{x_i\})} = \calK\one_{\{x_i\}}\in L^p(\mu)$. Thus, $\nu_i = \mu(F^{-1}(\{x_i\})) = \int|\one_{F^{-1}(\{x_i\})}|^p\,d\mu < \infty$.

Recall that, for any finite or countable index set $I$, the space $c_0(I)$ consists of those $a\in\ell^\infty(I)$ satisfying
\begin{align*}
\{i\in I:|a_i|\geq\varepsilon\}\text{ is finite for every } \varepsilon>0.
\end{align*}
In particular, $c_0(I)=\mathbb C^I$ when $I$ is finite.

\begin{thm}[{\cite[Theorems 1 and 3]{Takagi1992}}]\label{t:compact_Lp_det} Let $1\leq p<\infty$. The composition operator $\calK:L^p(\mu)\to L^p(\mu)$ is well defined and compact if and only if
\begin{align*}
\nu(\calX_{\mathrm{na}})=0 \qquad\text{and}\qquad \left(\frac{\nu_i}{\mu_i}\right)_{i\in I}\in c_0(I).
\end{align*}
The final condition requires every ratio $\nu_i/\mu_i$ to be finite, also when $I$ itself is finite.
\end{thm}
\begin{proof}
Suppose first that $\calK$ is compact. Boundedness gives $\nu\ll\mu$ and $g=d\nu/d\mu\in L^\infty(\mu)$, see Theorem \ref{t:det:invariance:Lp}. If $\nu(\calX_{\mathrm{na}})>0$, choose $\varepsilon>0$ and a nonatomic measurable set $E\subset\calX_{\mathrm{na}}$ such that $0<\mu(E)<\infty$ and $g\geq\varepsilon$ on $E$. Split $E$ into pairwise disjoint measurable sets $E_n$ of positive measure. The normalized indicators
\begin{align*}
f_n=\mu(E_n)^{-1/p}\one_{E_n}
\end{align*}
have norm one, whereas, for $n\neq m$,
\begin{align*}
\|\calK f_n-\calK f_m\|_p^p =\int |f_n-f_m|^p g\,d\mu \geq 2\varepsilon.
\end{align*}
The resulting separated sequence contradicts compactness. Hence $\nu(\calX_{\mathrm{na}})=0$.

On each atom, $g(x_i)=\nu_i/\mu_i < \infty$. If these ratios do not belong to $c_0(I)$, choose $\varepsilon>0$ and distinct indices $i_n$ with $\nu_{i_n}/\mu_{i_n}\geq\varepsilon$. The functions $f_n=\mu_{i_n}^{-1/p}\one_{\{x_{i_n}\}}$ have unit norm, their images have pairwise disjoint supports, and
\begin{align*}
\|\calK f_n-\calK f_m\|_p^p =\frac{\nu_{i_n}}{\mu_{i_n}}+\frac{\nu_{i_m}}{\mu_{i_m}} \geq2\varepsilon.
\end{align*}
Compactness again rules out this separated sequence.

Conversely, assume the two displayed conditions, and set $a_i=\nu_i/\mu_i$. Then $F(x)\in\calX_{\mathrm a}$ for $\mu$-almost every $x$. Let $E_i=F^{-1}(\{x_i\})$ and $I_+=\{i:\nu_i>0\}$. Note that $\calK = WDT$, where $T:L^p(\mu)\to\ell^p(I)$, $D:\ell^p(I)\to\ell^p(I_+)$, and $W:\ell^p(I_+)\to L^p(\mu)$ are defined by
\[
(Tf)_i=\mu_i^{1/p}f(x_i),\qquad (Dz)_i=a_i^{1/p}z_i, \qquad\text{and}\qquad Wz=\sum_{i\in I_+}z_i\,\nu_i^{-1/p}\one_{E_i},
\]
respectively. All three operators are bounded with $D$ even being compact because $a\in c_0(I)$. Hence, $\calK$ is compact.
\end{proof}

\begin{cor} Let $1\leq p<\infty$. \begin{enumerate} \item[{\rm (a)}] If $\mu$ is nonzero and nonatomic, then $\calK:L^p(\mu)\to L^p(\mu)$ is not compact. \item[{\rm (b)}] If $0<\inf_{i\in I}\mu_i\leq\sup_{i\in I}\mu_i<\infty$ and $\calK$ is compact, then $F(\calX_{\mathrm a})$ is finite.
\end{enumerate}
\end{cor}
\begin{proof}
For~{\rm(a)}, the theorem would force $\nu(\calX)=0$, whereas $\nu(\calX)=\mu(\calX)\neq0$. For~{\rm(b)}, let $c=\inf_i\mu_i$ and $C=\sup_i\mu_i$. If $j\in\pi_F(I)$, then \eqref{e:numu_j} gives $\nu_j/\mu_j\geq c/C$. The $c_0(I)$ condition therefore forces $\pi_F(I)$ to be finite.
\end{proof}

\subsubsection{Compactness on reproducing kernel Hilbert spaces}
The range-space construction from Section~\ref{s:invariance:det:RKHS} factors composition as a coisometry onto $\calH_F$ followed by the canonical inclusion $\calH_F\hookrightarrow\calH$. Thus all compactness lies in this inclusion. Theorem~\ref{t:RKHS_comp_emb} turns that abstract condition into uniform spectral decay of finite Gram matrices. For a kernel $k_0$ and an ordered finite family $X=(x_1,\ldots,x_m)$ in $\calX$, the associated matrix is
\begin{align*}
K_{k_0,X}=(k_0(x_i,x_j))_{i,j=1}^m.
\end{align*}
For a matrix $A$, let $A^\dagger$ denote its Moore--Penrose pseudoinverse. For a positive semidefinite $m\times m$ matrix $A$, let $(\lambda_n(A))$ be its eigenvalues in nonincreasing order, counted with multiplicity and completed by zeros beyond $m$.

\begin{thm}\label{t:RKHS_comp_det}
The following statements are equivalent:
\begin{enumerate}
\item[{\rm (i)}] The space $\calH$ is Koopman invariant and $\calK_\calH:\calH\to\calH$ is compact;
\item[{\rm (ii)}] $\calH_F$ is compactly embedded in $\calH$;
\item[{\rm (iii)}] $k_F\preceq c^2k$ for some $c>0$ and
\begin{align*}
\sup_{\substack{X=(x_1,\ldots,x_m)\\m<\infty}} \lambda_n\!\left( (K_{k,X}^{\dagger})^{1/2} K_{k_F,X} (K_{k,X}^{\dagger})^{1/2}\right) \longrightarrow0\qquad(n\to\infty).
\end{align*}
\end{enumerate}
Here, $K_{k_F,X}=(k(F(x_i),F(x_j)))_{i,j=1}^m$.
\end{thm}

The domination clause in {\rm(iii)} is the boundedness component of the criterion and follows here from Theorem~\ref{t:RKHS_det_inv}; the spectral decay is its compactness component. The proof of Theorem \ref{t:RKHS_comp_det} is essentially the same as that of Theorem \ref{t:RKHS_comp} for the stochastic setting and is therefore omitted.

For translation-invariant kernels, compactness of the Koopman operator has immediate geometric consequences.

\begin{prop}\label{p:nec}
If $\calK_\calH$ is compact, the following statements hold:
\begin{enumerate}
\item[{\rm (i)}] If $\calX\subset\mathbb R^n$ and $k(x,y)=\Phi(x-y)$ with a positive definite function $\Phi : \R^n\to\R$ satisfying
\begin{align*}
\sup_{\|u\|\geq r}\operatorname{Re}\Phi(u)<\Phi(0) \qquad\text{for every }r>0,
\end{align*}
then $F(\calX)$ is bounded.
\item[{\rm (ii)}] If $\inf_{x\in\calX}k(x,x)>0$, then every sequence $(x_n)\subset\calX$ admits an $f\in\calH$ for which $\limsup_{n\to\infty}|f(x_n)|>0$.
\end{enumerate}
\end{prop}
\begin{proof} For~{\rm(i)}, the feature vectors $\phi_x=k(x,\cdot)$ have constant norm $\sqrt{\Phi(0)}$, and $\calK_\calH^*\phi_x=\phi_{F(x)}$. Given any sequence $(x_n)\subset\calX$, compactness of $\calK_\calH^*$ supplies a subsequence along which $(\phi_{F(x_n)})$ is Cauchy. Since
\begin{align*}
\|\phi_{F(x)}-\phi_{F(y)}\|_\calH^2 =2\bigl(\Phi(0)-\operatorname{Re}\Phi(F(x)-F(y))\bigr),
\end{align*}
the strict-separation hypothesis forces $(F(x_n))$ along that subsequence to be Cauchy. Hence, every sequence in $F(\calX)$ has a convergent subsequence in $\mathbb R^n$, which means that $F(\calX)$ is relatively compact in $\mathbb R^n$, and thus bounded.

For~{\rm(ii)}, suppose instead that $f(x_n)\to0$ for every $f\in\calH$. Then $\phi_{x_n}\rightharpoonup0$; uniform boundedness also shows that this sequence is norm bounded. A compact operator maps weakly null sequences to norm-null sequences, and therefore $\phi_{F(x_n)}=\calK_\calH^*\phi_{x_n}\to0$ in norm. But
\begin{align*}
\|\phi_{F(x_n)}\|_\calH^2 =k(F(x_n),F(x_n))\geq\inf_{x\in\calX}k(x,x)>0,
\end{align*}
which is a contradiction.
\end{proof}

Functions in Gaussian and Sobolev RKHSs tend to zero at infinity. In the Sobolev case, provided that $\Omega$ admits a bounded extension operator $E:H^s(\Omega)\to H^s(\R^n)$, $s>\frac n2$ (e.g., $\Omega=\R^n$ or $\Omega$ with uniform Lipschitz boundary), this follows from Sobolev's embedding theorem $H^s(\Omega)\subset C_0(\ol\Omega)$, applied to $Ef$, whereas in the Gaussian case, it is a consequence of $\lim_{|x|\to\infty}g(x) = 0$ for all $g\in\linspan\{k_y : y\in\R^n\}$. Moreover, $\inf_x k(x,x)>0$ in both cases: for the Gaussian kernel, $k(x,x)=1$, and in the Sobolev case, the restrictions to $\Omega$ of the translates $\psi(\cdot-x)$ of a fixed $\psi\in C_c^\infty(\R^n)$ with $\psi(0)=1$ have uniformly bounded norms and take the value $1$ at $x$. Hence, from Proposition \ref{p:nec} (ii), we immediately draw the following conclusion.

\begin{cor}
No deterministic Koopman operator on the Gaussian RKHS $\calH_\sigma(\calX)$ for unbounded $\calX\subset\R^n$ is compact. Moreover, no deterministic Koopman operator on the Sobolev space $H^s(\Omega)$, $s>\frac n2$, for unbounded $\Omega\subset\R^n$ admitting a bounded extension operator $H^s(\Omega)\to H^s(\R^n)$ {\rm (}e.g., $\Omega$ with uniform Lipschitz boundary{\rm )} is compact.
\end{cor}

For $\calX=\R^n$, the Gaussian case of this corollary is also contained in~\cite{Ikeda22}, where it is shown that bounded composition operators on RKHSs associated with a large class of analytic positive definite functions, including the Gaussian kernel, cannot be compact.

\subsection{Markov kernels}
In the stochastic setting, a new mechanism appears. In the deterministic case the transition laws are Dirac masses, and distinct Dirac masses lie at total-variation distance two; total-variation precompactness therefore collapses to finite range. For a nondeterministic Markov kernel, by contrast, the transition laws may form a compact family while ranging over infinitely many states. Averaging and smoothing thus provide an alternative source of compactness.

\subsubsection{Compactness on \texorpdfstring{$C_b(\calX)$}{Cb(X)}}
On $C_b(\calX)$, the relevant topology on transition laws is total variation. For a finite signed Borel measure $\tau$ on the Polish space $\calX$, its total-variation norm is, see \cite[Theorem IV.6.2]{DunfSchw58},
\begin{align*}
\|\tau\|_{TV} = |\tau|(\calX) = \sup_{\|f\|_\infty\leq1,\ f\in C_b(\calX)} \left|\int f\,d\tau\right|.
\end{align*}
For two probability measures $\mu_1$ and $\mu_2$, we have $\|\mu_1-\mu_2\|_{TV} =2\sup_{A\in\calB_\calX}|\mu_1(A)-\mu_2(A)|$, see \cite[Proposition 4.2]{Levin17}.

\begin{thm}\label{t:Cb_compact}
Let $\calX$ be a Polish space and assume that
$C_b(\calX)$ is Koopman invariant. Then the following
statements are equivalent:
\begin{enumerate}
\item[{\rm (i)}]
The Koopman operator
$\calK:C_b(\calX)\to C_b(\calX)$ is compact;
\item[{\rm (ii)}]
The set $\{\rho_x:x\in\calX\}$ is totally
bounded in total variation.
\end{enumerate}
If {\rm (i)} or {\rm (ii)} holds,
\begin{enumerate}
\item[{\rm (iii)}] the map $x\mapsto\rho_x$ is continuous in total variation.
\end{enumerate}
If $\calX$ is compact, {\rm (i)}--{\rm (iii)} are equivalent.
\end{thm}
\begin{proof}
For each $x\in\calX$, identify $\rho_x$ with the bounded linear functional on $C_b(\calX)$ given by
\begin{align*}
\rho_x(f) = \int f\,d\rho_x.
\end{align*}
Since $\calX$ is Polish, the total variation norm of a finite signed Borel measure agrees with its norm as a functional on $C_b(\calX)$.

{\rm (i)} $\Rightarrow$ {\rm (ii)}. Suppose that $\calK$ is compact. By Schauder's theorem, its adjoint $\calK^*$ is compact as well. Since for $x\in\calX$ we have $\calK^*\delta_x = \rho_x$, and the family $(\delta_x)_{x\in\calX}$ is bounded in $C_b(\calX)^*$, the set $\{\rho_x:x\in\calX\}$ is relatively compact in total variation, and hence totally bounded.

{\rm (ii)} $\Rightarrow$ {\rm (i)}. Suppose that $\{\rho_x:x\in\calX\}$ is totally bounded in total variation. Then
\begin{align*}
M:= \overline{\{\rho_x:x\in\calX\}}^{\|\cdot\|_{\mathrm{TV}}}
\end{align*}
is a compact metric space consisting of probability measures. 

We first show (iii). Let $x_n\to x$ in $\calX$. Suppose that $\rho_{x_n}$ does not converge to  $\rho_x$ in total variation. Then there exist $\varepsilon>0$ and a subsequence, again denoted by $(x_n)$, such that
\begin{align*}
\|\rho_{x_n}-\rho_x\|_{\mathrm{TV}}\geq\varepsilon
\end{align*}
for all $n$. By compactness of $M$, there exists a further subsequence such that
\begin{align*}
\rho_{x_{n_k}}\to\eta\quad\text{in total variation}
\end{align*}
for some $\eta\in M$. For every $f\in C_b(\calX)$, Koopman invariance implies $\calK f\in C_b(\calX)$. Consequently,
\begin{align*}
\int f\,d\eta = \lim_{k\to\infty}\int f\,d\rho_{x_{n_k}} = \lim_{k\to\infty}(\calK f)(x_{n_k}) = (\calK f)(x) = \int f\,d\rho_x.
\end{align*}
Since bounded continuous functions determine finite Borel measures on Polish spaces, we obtain $\eta=\rho_x$, a contradiction. Thus, (iii) follows.

To prove (i), we factor $\calK$ through $C(M)$ by defining $A : C_b(\calX)\to C(M)$ and $T : C(M)\to C_b(\calX)$ by
\begin{align*}
Af(\eta) = \int f\,d\eta
\qquad\text{and}
\qquad (Tg)(x)=g(\rho_x).
\end{align*}
Note that both $A$ and $T$ are bounded operators and that $\calK = TA$.

The image of the unit ball under $A$ is bounded in $C(M)$ and $1$-Lipschitz on the compact metric space $M$ and thus relatively compact by the Arzel\`a--Ascoli theorem. Hence, $A$ is compact and so is $\calK$.

It remains to prove that (iii) implies (ii) when $\calX$ is compact. But this is trivial as the continuous map $x\mapsto\rho_x$ maps the compact set $\calX$ to the compact set $\{\rho_x : x\in\calX\}$ which, consequently, is totally bounded.
\end{proof}

\begin{rem} Total-variation continuity is stronger than the strong Feller property. Indeed, the latter asks that $x\mapsto\int f\,d\rho_x$ be continuous for every bounded measurable $f$, whereas the continuity of $x\mapsto\rho_x$ in TV-norm makes this continuity locally uniform over the unit ball of bounded measurable functions.
\end{rem}

Theorem \ref{t:C_det_compact} shows that the Koopman operator for deterministic systems $x^+ = F(x)$ with infinite image $F(\calX)$ is never compact on $C_b(\calX)$, even when $\calX$ is compact. The next example, however, demonstrates that there is a large class of \emph{stochastic} systems whose Koopman operators are all compact on $C(\calX)$.

\begin{ex}
Suppose that $\calX$ is compact and $d\rho_x(y) = h(x,y)\,d\mu(y)$ with a continuous function $h:\calX^2\to[0,\infty)$ and a finite measure $\mu$. Uniform continuity then gives
\begin{align*}
\|\rho_{x'}-\rho_x\|_{TV} \leq\int|h(x',y)-h(x,y)|\,d\mu(y)\longrightarrow0\qquad\text{as $x'\to x$},
\end{align*}
so $\calK$ is compact on $C(\calX)$.
\end{ex}

\subsubsection{Compactness on \texorpdfstring{$L^\infty$}{L-infinity}}
Passing from $C_b(\calX)$ to $L^\infty(\mu)$ removes pointwise continuity from view. What remains is the essential total-variation geometry of the transition laws: compactness is equivalent to total boundedness after the deletion of a null set. The partition formulation below makes this criterion intrinsic to equivalence classes. For a measurable $\Delta\subset\calX$ and a nonnegative function $d$ on $\calX^2$, we use the convention
\begin{align*}
\mathop{\mu\text{-}\mathrm{ess\,sup}}_{x,y\in\Delta}d(x,y) :=\inf_{\mu(N)=0}\ \sup_{x,y\in\Delta\backslash N}d(x,y),
\end{align*}
where the supremum over the empty set is zero.

\begin{thm}\label{t:compact:stoch:Linf}
Assume that $L^\infty(\mu)$ is Koopman invariant. Then the following statements are equivalent: \begin{enumerate}
\item[{\rm (i)}] $\calK:L^\infty(\mu)\to L^\infty(\mu)$ is compact;
\item[{\rm (ii)}] for every $\varepsilon>0$ there is a finite measurable partition $\pi$ of $\calX$ such that, for every $\Delta\in\pi$,
\begin{align}
\mathop{\mu\text{-}\mathrm{ess\,sup}}_{x,y\in\Delta} \|\rho_x-\rho_y\|_{TV} \le \varepsilon; \label{e:Linfty_comp}
\end{align}
\item[{\rm (iii)}] there is a measurable $\mu$-null set $N$ such that $\{\rho_x:x\in\calX\backslash N\}$ is totally bounded in total variation.
\end{enumerate}
\end{thm}
\begin{proof}
We start with an essential oscillation result from \cite{Cherkas70}. For this, let $B$ be the unit ball of $L^\infty(\mu)$. According to \cite{Cherkas70}, (i) holds if and only if
\begin{enumerate}
\item[{\rm (i')}] For any $\veps>0$ there exists a finite measurable partition $\pi$ of $\calX$ such that, for each $\Delta\in\pi$,
\begin{align}\label{e:Linfty_operator_osc}
\sup_{f\in B}\ \mathop{\mu\text{-}\mathrm{ess\,sup}}_{x,y\in\Delta} |(\calK f)(x)-(\calK f)(y)|\leq\veps.
\end{align}
\end{enumerate}
The passage from essential oscillation to total variation requires one exceptional null set for all test functions. The countable norming class from Theorem~\ref{t:TV_with_calE} provides precisely this reduction. Let $\calE$ be the class from that theorem. Suppose first that~{\rm(i')} holds. Fix $\veps>0$, a corresponding partition $\pi$ and $\Delta\in\pi$. For each $E\in\calE$, apply \eqref{e:Linfty_operator_osc} to $f_E = 2\one_E-\one_{\calX}\in B$, and let $N_\Delta$ denote the countable union of the corresponding null sets. Simultaneously for all $E\in\calE$, we then have
\begin{align*}
2|(\rho_x-\rho_y)(E)|\leq\veps \qquad(x,y\in\Delta\backslash N_\Delta).
\end{align*}
Theorem~\ref{t:TV_with_calE} gives $\|\rho_x-\rho_y\|_{TV}\leq\varepsilon$, and proves~{\rm(ii)}.

Conversely, assume~{\rm(ii)}, fix $\veps>0$, a corresponding partition $\pi$ and $\Delta\in\pi$. Given a class $[f]$ in the unit ball, we may assume that $|f|\le 1$ everywhere. Koopman invariance implies $\nu\ll\mu$, so $\calK f$ represents $\calK[f]$. For $x,y\in\Delta\backslash N_\Delta$,
\begin{align*}
|(\calK f)(x)-(\calK f)(y)| \leq\|\rho_x-\rho_y\|_{TV}\le\varepsilon.
\end{align*}
Thus (i') holds.

It remains to compare~{\rm(ii)} and~{\rm(iii)}. For any $k\in\N$, apply~{\rm(ii)} with $\varepsilon=1/k$ and take the countable union $N$ of the finitely many exceptional null sets arising for each $k$. Thereby, include possible null cells $\Delta$. Given $\veps>0$, choose $k$ such that $1/k < \veps$, and select one point from each non-null cell, after deleting $N$. These points then yield a finite $\veps$-net of $\{\rho_x : x\in\calX\backslash N\}$, proving~{\rm(iii)}. Conversely, Corollary~\ref{c:meas_map_meas} makes $x\mapsto\rho_x$ measurable into its separable total-variation range. Pull back a finite cover of the range by open total-variation balls of radius $\varepsilon/3$, form a measurable partition from the preimages by successive set differences, and adjoin the exceptional null set to one cell. This partition satisfies~{\rm(ii)}.
\end{proof}

\subsubsection{Compactness on \texorpdfstring{$L^p(\mu)$}{Lp(mu)}, \texorpdfstring{$1\leq p<\infty$}{1 <= p < infinity}}
Throughout this subsection, it is assumed that $\nu\leq C\mu$ for some $C>0$. This assumption supplies boundedness of the Koopman operator on every $L^p(\mu)$, $p\in [1,\infty]$. Compactness at one finite exponent then propagates through the entire open interpolation scale.

\begin{thm}\label{t:Lp_comp_interpolation}
Assume $\nu\leq C\mu$ for some $C>0$. If $\calK$ is compact on $L^{p_0}(\mu)$ for some $p_0\in[1,\infty)$, then it is compact on $L^p(\mu)$ for every $p\in (1,\infty)$.
\end{thm}
\begin{proof} If $1<p<p_0$, interpolate the compact map on $L^{p_0}$ with the bounded map on $L^1$. If $p_0<p<\infty$, interpolate it with the bounded map on $L^\infty$. Krasnosel'ski\u{\i}'s compact interpolation theorem gives each conclusion; see \cite[Corollary 3.2]{Cobos98}.
\end{proof}

Compact interpolation leaves the $L^1$ endpoint unresolved. Normalized indicators provide a necessary test at every finite exponent and, at $p=1$, a complete one. For $p\in [1,\infty)$, consider the family
\begin{align*}
\calA_p = \left\{\frac{\rho(\,\cdot\,,A)}{\mu(A)^{1/p}}: A\in\calB_\calX,\ 0<\mu(A)<\infty\right\}\subset L^p(\mu).
\end{align*}

\begin{lem}\label{l:compact1:Lp}
Let $p\in [1,\infty)$. Then the following statements hold:
\begin{enumerate}
\item[{\rm (i)}] If $\calK:L^p(\mu)\to L^p(\mu)$ is compact, then $\calA_p$ is relatively compact in $L^p(\mu)$.
\item[{\rm (ii)}] If $\calA_1$ is relatively compact in $L^1(\mu)$, then $\calK$ is compact on $L^p(\mu)$ for every $p\in [1,\infty)$.
\end{enumerate}
In particular, $\calK$ is compact on $L^1(\mu)$ if and only if $\calA_1$ is relatively compact in $L^1(\mu)$.
\end{lem}
\begin{proof} If $0<\mu(A)<\infty$, then $\mu(A)^{-1/p}\one_A$ has unit $L^p$ norm and $\calK$ maps it to the corresponding element of $\calA_p$. This proves~{\rm(i)}.

For~{\rm(ii)}, let $f=\sum_{j=1}^n a_j\one_{A_j}$, where the $A_j$ are pairwise disjoint and have positive finite measure. Then
\begin{align*}
\calK f=\sum_{j=1}^n a_j\mu(A_j) \frac{\rho(\,\cdot\,,A_j)}{\mu(A_j)}, \qquad \sum_{j=1}^n|a_j|\mu(A_j)=\|f\|_1.
\end{align*}
Simple functions with support of finite measure are dense in $L^1(\mu)$. It follows that $\calK(B_{L^1})$ lies in the closed absolutely convex hull of $\calA_1$. This hull is compact whenever $\calA_1$ is relatively compact \cite[Theorem 4.8.9(b)]{Narici11}. Hence $\calK$ is compact on $L^1$, and Theorem~\ref{t:Lp_comp_interpolation} extends compactness to $p\in (1,\infty)$.
\end{proof}

Relative compactness in $L^1$ has three distinct obstructions: concentration on sets of small measure, escape to infinity, and local mean oscillation. The next result excludes all three on an arbitrary $\sigma$-finite measure space.

The measure-theoretic compactness theorem behind this trichotomy is proved for arbitrary subsets of $L^1$ in Appendix~\ref{a:L1_compactness}. Its application to the normalized kernel sections $\calA_1$ gives the following exact criterion.

\begin{thm}\label{t:compact2:Lp} 
The operator $\calK:L^1(\mu)\to L^1(\mu)$ is compact if and only if the following three conditions hold: 
\begin{enumerate} 
\item[{\rm (i)}] 
\begin{align*}
\lim_{\delta\searrow 0}\;\sup_{\substack{0<\mu(A)<\infty\\\mu(B)\leq\delta}} \frac1{\mu(A)}\int_B\rho_x(A)\,d\mu(x)=0; 
\end{align*}
\item[{\rm (ii)}] for every $\varepsilon>0$ there is $E\in\calB_\calX$ with $\mu(E)<\infty$ such that
\begin{align*}
\sup_{0<\mu(A)<\infty} \frac1{\mu(A)} \int_{\calX\backslash E}\rho_x(A)\,d\mu(x)<\varepsilon;
\end{align*}
\item[{\rm (iii)}] for every~$\varepsilon>0$ and every $\Delta \in \calB_\calX$ with $0<\mu(\Delta)<\infty$, there are measurable sets $\Delta_1,\ldots,\Delta_n\subset\Delta$ of positive measure such that, for every $A\in\calB_\calX$ with $0<\mu(A)<\infty$, some $j\in\{1,\ldots,n\}$ satisfies
\begin{align}
\int_{\Delta_j}\int_{\Delta_j} |\rho_x(A)-\rho_z(A)|\,d\mu(z)\,d\mu(x) <\varepsilon\,\mu(A)\mu(\Delta_j)^2. \label{e:kernel_bocce_oscillation}
\end{align}
\end{enumerate}
Whenever these conditions hold, $\calK$ is compact on every $L^p(\mu)$ with $1\leq p<\infty$. 
\end{thm}

\begin{proof} 
For $A\in\calB_\calX$ with $0<\mu(A)<\infty$, define 
\begin{align*}
    s_A:=\frac{\rho(\,\cdot\,,A)}{\mu(A)}. 
\end{align*}
Apply Theorem~\ref{t:Lp_comp_family} to the family
\begin{align*}
\calA_1=\left\{s_A:0<\mu(A)<\infty\right\}.
\end{align*}
It is bounded because
\begin{align*}
\|s_A\|_1=\frac{\nu(A)}{\mu(A)}\leq C.
\end{align*}
Conditions~{\rm(i)}, {\rm(ii)}, and {\rm(iii)} are, respectively, the uniform-integrability, tightness, and local mean-oscillation conditions in that theorem. The normalization by $\mu(A)$ is essential in each case. Lemma~\ref{l:compact1:Lp} now gives both the equivalence and the final assertion. 
\end{proof}

A more direct sufficient condition is available when the transition kernel has a jointly Borel-measurable density $h:\calX\times\calX\to[0,\infty)$ with respect to $\mu$:
\begin{align*}
\rho_x(A)=\int_A h(x,y)\,d\mu(y),\qquad \int h(x,y)\,d\mu(y)=1
\end{align*}
for $x\in\calX$ and $A\in\calB_\calX$. In this notation, $\nu\leq C\mu$ means
\begin{align*}
\mathop{\mathrm{ess\,sup}}_{y\in\calX} \int h(x,y)\,d\mu(x)\leq C.
\end{align*}

\begin{thm}[{\cite[Theorem 11.6]{Joergens82}}]
Let $p\in (1,\infty)$ and $q = p/(p-1)$. If
\begin{align}\label{e:mixed_kernel_compact}
\int\left(\int|h(x,y)|^q\,d\mu(y)\right)^{p/q} d\mu(x)<\infty,
\end{align}
then $\calK$ is compact on $L^p(\mu)$ and, consequently, on every $L^r(\mu)$ with $1<r<\infty$. 
\end{thm}
\begin{proof}
First of all, condition \eqref{e:mixed_kernel_compact} implies $\nu\ll\mu$, hence $\calK$ is well defined on equivalence classes. Moreover, \eqref{e:mixed_kernel_compact} places $x\mapsto h(x,\cdot)$ in $L^p(\mu;L^q(\mu))$. Indeed, $L^q(\mu)$ is separable because $\calX$ is Polish, $\mu$ is $\sigma$-finite, and $q<\infty$; joint measurability therefore makes this map strongly measurable. Finite sums $h_N(x,y)=\sum_{j=1}^{m_N}a_j(x)b_j(y)$ with $a_j\in L^p(\mu)$ and $b_j\in L^q(\mu)$ are dense in this Bochner space. Define
\begin{align*}
(\calK_Nf)(x):=\int h_N(x,y)f(y)\,d\mu(y).
\end{align*}
These operators have finite rank, while pointwise H\"older's inequality followed by the $L^p$ norm gives
\begin{align*}
\|\calK-\calK_N\|_{L^p\to L^p} \leq\|h-h_N\|_{L^p(\mu;L^q(\mu))}\longrightarrow0.
\end{align*}
Thus, $\calK$ is an operator-norm limit of finite-rank operators on $L^p$. Theorem~\ref{t:Lp_comp_interpolation} gives the remaining exponents.
\end{proof}

The preceding criterion extracts compactness from integrability of the transition density. Diffusion semigroups also admit a spectral route: discreteness of the $L^2$ spectrum gives compactness at $p=2$, after which Theorem~\ref{t:Lp_comp_interpolation} covers the open interpolation scale. The overdamped Langevin dynamics from Example~\ref{ex:olan} illustrate this second mechanism.

\begin{ex}[Overdamped Langevin dynamics]\label{ex:olan_compact}
Recall the function $W$, the self-adjoint Schr\"odinger operator $H=-\frac D\beta\Delta+W$, and the unitary operator $U:L^2(\mu)\to L^2(\R^n)$ from Example \ref{ex:olan}, and suppose additionally
that
\begin{align*}
\lim_{\|x\|\to\infty}\|\nabla V(x)\| = \infty.
\end{align*}
Since $\Delta V$ is bounded, this condition implies $W(x)\to\infty$ as $\|x\|\to\infty$. Consequently, $H$ has compact resolvent \cite{Simon2009}, and hence $e^{-tH}$ is compact for every $t>0$ by the spectral theorem. By unitary equivalence, $\calK_t=U^{-1}e^{-tH}U$ is compact on $L^2(\mu)$. Since $\calK_t$ is a contraction on $L^p(\mu)$ for every $p\in[1,\infty)$, Theorem \ref{t:Lp_comp_interpolation} implies that $\calK_t$ is compact on $L^p(\mu)$ for every $p\in(1,\infty)$.
\end{ex}

\subsubsection{Compactness on reproducing kernel Hilbert spaces}
On an RKHS, stochastic averaging is encoded by the range kernel $k_\rho$. The deterministic characterization is %
therefore maintained if %
$\calH_F$ is replaced by $\calH_\rho$. 
We remain in the bounded-diagonal Hermitian setting of Subsection~\ref{s:invariance:stoch:RKHS}.

\begin{thm}\label{t:RKHS_comp} 
The following statements are equivalent: 
\begin{enumerate} 
\item[{\rm (i)}]   The space $\calH$ is Koopman invariant and $\calK_\calH:\calH\to\calH$ is compact; 
\item[{\rm (ii)}]  $\calH_\rho$ is compactly embedded in $\calH$;
\item[{\rm (iii)}] $k_\rho\preceq c^2k$ for some $c>0$ and
\begin{align*}
\sup_{\substack{X=(x_1,\ldots,x_m)\\m<\infty}} \lambda_n\!\left( (K_{k,X}^{\dagger})^{1/2} K_{k_\rho,X} (K_{k,X}^{\dagger})^{1/2}\right) \longrightarrow0\qquad(n\to\infty).
\end{align*}
\end{enumerate}
\end{thm}

The domination clause in {\rm(iii)} is the boundedness component of the criterion and follows here from Theorem~\ref{t:RKHS_inv}; the spectral decay is its compactness component.

\begin{proof}
Since all three conditions (i)--(iii) include Koopman invariance of $\calH$, we may assume that $\calK_\calH : \calH\to\calH$ is well defined and bounded. Let $S:\calH\to\calH_\rho$ be defined by $Sf = \calK f$ for $f\in\calH$, equip $\calH_\rho$ with the range norm in \eqref{e:H_rho}, and let $J:\calH_\rho\hookrightarrow\calH$ be the inclusion. The restriction of $S$ to $(\ker S)^\perp$ is an isometric isomorphism onto $\calH_\rho$; equivalently,
\begin{align*}
SS^*=I_{\calH_\rho}.
\end{align*}
The identities $\calK_\calH = JS$ and $J = \calK_\calH S^*$ show that one of the operators $\calK_\calH$ and $J$ is compact precisely when the other is. This proves the equivalence of {\rm(i)} and {\rm(ii)}. Theorem~\ref{t:RKHS_comp_emb} yields the equivalence of {\rm(ii)} and {\rm(iii)}.
\end{proof}

Theorem~\ref{t:RKHS_comp} concerns the Koopman operator $\calK_\calH:\calH\to\calH$ on a single invariant RKHS. For diffusions, parabolic smoothing yields a related statement in which observables and their images lie in different spaces~\cite{hertel2026verifiable}. Consider the stochastic differential equation%
\begin{align}\label{e:sde_hs}
\mathrm{d}X_t=b(X_t)\,\mathrm{d}t + \sigma(X_t)\,\mathrm{d}W_t,\qquad t \geq%
0,\qquad X_0=x\in\calX\subseteq\R^n,
\end{align}
with $n$-dimensional Brownian motion $(W_t)_{t\geq0}$, drift $b$, and diffusion coefficient $\sigma$.
Here, observables $g$ are defined on $\R^n$, whereas their images $\calK_tg$ are evaluated on a domain $\calX\subseteq\R^n$ only. Under the assumptions of Theorem~\ref{t:sde_HS} below, the transition kernels associated with the SDE~\eqref{e:sde_hs} admit densities $\eta_t(x,\cdot)$ with respect to Lebesgue measure. That is, the time-$t$ transition kernel $\rho^t_x(\Delta):=\mathbb P(X_t\in\Delta\mid X_0=x)$, which is a Markov kernel on $\R^n$ in the sense of Section~\ref{sec:definition}, satisfies $\mathrm{d}\rho^t_x(y)=\eta_t(x,y)\,\mathrm{d}y$, so that
\begin{align}\label{e:sde_koopman}
(\calK_tg)(x)=\mathbb E\left[g(X_t)\mid X_0=x\right] = \int_{\R^n}g(y)\,\eta_t(x,y)\,\mathrm{d}y, \qquad x\in\calX,
\end{align}
for bounded measurable $g$ and $t>0$. Thus $\calK_t$ is the Koopman operator of $\rho^t$, evaluated at $x\in\calX$ only, and $\eta_t$ plays the role of the density $p$ in a two-space version of Theorem~\ref{t:sob_stoch} for this kernel.
The following result shows that, under uniform ellipticity and regularity of the coefficients~$b$ and~$\sigma$, $\calK_t$ is Hilbert--Schmidt, and thus compact, as an operator from $L^2(\R^n)$ into an RKHS~$\calH$ of functions on~$\calX$. Here, all $L^2$ spaces are taken with respect to Lebesgue measure, $H^\ell$ denotes the Sobolev space of order $\ell\in\N$, and $\mathrm{HS}$ and $\mathrm{Tr}$ denote the classes of Hilbert--Schmidt and trace-class operators. We write $\calH\simeq H^\ell(\calX)$ if $\calH$ and $H^\ell(\calX)$ coincide as sets of functions and have equivalent norms, where the elements of $H^\ell(\calX)$ are identified with their continuous representatives; this requires $\ell>n/2$. Finally, $\eta_t\in L^2(\R^n;\calH)$ means that the map $y\mapsto\eta_t(\cdot,y)$, whose values are functions of $x\in\calX$, is square-integrable with values in~$\calH$, that is, $\|\eta_t\|_{L^2(\R^n;\calH)}^2=\int_{\R^n}\|\eta_t(\cdot,y)\|_\calH^2\,\mathrm{d}y<\infty$.

\begin{thm}[{\cite[Theorem 4.11]{hertel2026verifiable}}]\label{t:sde_HS} 
Let $t>0$ and let $\ell\in\N$ with $\ell>n/2$. Suppose that
\begin{enumerate} 
\item[{\rm(i)}] $\sigma$ is uniformly elliptic, i.e., 
$\exists\,c_\sigma>0 : c_\sigma\|\xi\|_2^2\leq\|\sigma(x)^\top\xi\|_2^2$ for all $x,\xi\in\R^n$ and
\item[{\rm(ii)}] $b\in C_b^\ell(\R^n,\R^n)$ and $\sigma\in C_b^\ell(\R^n,\R^{n\times n})$, 
where $C_b^\ell$ denotes the $\ell$-times continuously differentiable functions with bounded derivatives up to order $\ell$. 
\end{enumerate}
Then, we have 
\begin{align*}
\eta_t\in L^2(\R^n;\calH)
\qquad\text{and}\qquad
\calK_t\in\mathrm{HS}(L^2(\R^n),\calH)
\end{align*}
with $\|\calK_t\|_{\mathrm{HS}(L^2(\R^n),\calH)} =\|\eta_t\|_{L^2(\R^n;\calH)}$ if either of the following holds:
\begin{enumerate}
\item[{\rm(A)}] $\calX\subset\R^n$ is a bounded Lipschitz domain and $\calH\simeq H^\ell(\calX)$; then $\calK_t\in\mathrm{HS}(L^2(\R^n),L^2(\calX))$ and even $\calK_t\in\mathrm{Tr}(L^2(\R^n),L^2(\calX))$;
\item[{\rm(B)}] $\calX=\R^n$ and $\calH=\calH_{k_\tau}$ for the Gaussian-reweighted kernel
\begin{align*}
k_\tau(x,x')=m_\tau(x)^{-1}k_0(x,x')m_\tau(x')^{-1}
\end{align*}
with $m_\tau(x)=\exp(-\tfrac\tau2 x^\top\Gamma^{-1}x)$, $\tau>0$, $\Gamma\in\R^{n\times n}$ symmetric positive definite, and $\calH_{k_0}\simeq H^\ell(\R^n)$. 
\end{enumerate}
\end{thm}

Theorem~\ref{t:sde_HS} is not an instance of Theorem~\ref{t:RKHS_comp}, since $\calK_t$ maps $L^2(\R^n)$ into an RKHS of functions on~$\calX$. Its norm identity quantifies the Hilbert--Schmidt norm of this two-space operator. In case~{\rm(A)}, the condition $\eta_t\in L^2(\R^n;\calH)$ is~\eqref{eq:stoch_sob_density} for $p=\eta_t$, with $\Omega\times\Omega$ replaced by $\calX\times\R^n$, so Theorem~\ref{t:sde_HS} verifies the density condition of Theorem~\ref{t:sob_stoch} for uniformly elliptic diffusions. Since the diffusion is not confined to~$\calX$, however, this yields no invariance of a single space.

The underlying mechanism is parabolic smoothing in the backward variable~$x$. Uniform ellipticity and the regularity of the coefficients yield Gaussian bounds on $\partial_x^\alpha\eta_t(x,y)$ for $|\alpha|\le\ell$, and integration over~$y$ gives $\eta_t\in L^2(\R^n;\calH)$. If $\calX$ is bounded, as in case~{\rm(A)}, the resulting integral over $\calX\times\R^n$ is finite. On $\calX=\R^n$, however, each slice $\eta_t(\cdot,y)$ decays in~$x$, but its norm need not decay in~$y$.\footnote{For instance, for $b=0$ and $\sigma=\sqrt2\,I$, the transition density is the heat kernel $\eta_t(x,y)=g_t(y-x)$ with $g_t(z)=(4\pi t)^{-n/2}\exp(-\|z\|^2/(4t))$, so that $\|\eta_t(\cdot,y)\|_{H^\ell(\R^n)}=\|g_t\|_{H^\ell(\R^n)}$ is independent of~$y$ and hence $\eta_t\notin L^2(\R^n;H^\ell(\R^n))$.} This is the reason for the reweighted kernel in case~{\rm(B)}, whose norm $\|f\|_{\calH_{k_\tau}}=\|m_\tau f\|_{\calH_{k_0}}$ damps the slices $\eta_t(\cdot,y)$ for large $\|y\|$ and thereby restores square-integrability in~$y$. The price is that the diagonal $k_\tau(x,x)=m_\tau(x)^{-2}k_0(x,x)$ is unbounded, so that case~{\rm(B)} lies outside the bounded-diagonal setting of Subsection~\ref{s:invariance:stoch:RKHS}.

\section{Self-adjointness}\label{sec:selfadjointness}
In this section, we address self-adjointness of the Koopman operator. This property is valuable in numerical analysis because the spectrum is real and, for every $f\in\dom\calK\backslash\{0\}$ and $\lambda\in\C$, the spectral theorem gives
\[
 \operatorname{dist}(\lambda,\sigma(\calK)) \le \frac{\|(\calK - \lambda I)f\|}{\|f\|}.
\]
For a nonnormal operator, a small normalized residual instead provides pseudospectral information and does not in general certify a small distance to the spectrum; see~\cite{Colb23,ColbTown24} for rigorous residual-based spectral analysis. Self-adjointness also gives variational principles for eigenvalues in the appropriate isolated or discrete spectral regime. Such principles underpin the variational approach to slow processes in reversible molecular dynamics~\cite{noenusk13,nuske2017variational}. Eigenvector and eigenspace accuracy nevertheless depends on spectral separation. In essence, self-adjointness expresses a symmetry of the dynamics, whose form depends on the Hilbert space.

We will show that, for deterministic dynamics, the composition rule defines a bounded self-adjoint operator on $L^2(\mu)$ if and only if $\mu\circ F^{-1}=\mu$ and $F^2=\operatorname{id}$ $\mu$-almost everywhere, in which case $\calK$ is a self-adjoint involution with spectrum contained in $\{-1,1\}$ (Theorem~\ref{t:det:selfadjL2} and Corollary~\ref{c:det:selfadj:spectral}). For stochastic dynamics, the corresponding condition is symmetry of the joint measure $Q(\mathrm{d}x,\mathrm{d}y)=\mu(\mathrm{d}x)\rho_x(\mathrm{d}y)$ of two successive states, which forces stationarity, makes $\calK$ a self-adjoint contraction, and yields a Dirichlet identity together with intrinsic descriptions of the eigenspaces corresponding to $1$ and $-1$ (Theorem~\ref{t:sa_stoch}). On an invariant RKHS and, in the stochastic case, under the bounded-diagonal hypothesis, self-adjointness will be characterized by symmetry of the transported kernel features, that is, by $k(x,F(y))=k(F(x),y)$ and $\Psi_x(y)=\overline{\Psi_y(x)}$, respectively, where $\Psi_x\in\calH$ is determined by $\langle f,\Psi_x\rangle_\calH=\int f\,d\rho_x$ for $f\in\calH$ (Theorems~\ref{t:selfadjprep} and \ref{t:selfadj:stoch:rkhs}).

For RKHSs we retain the convention of Subsection~\ref{s:invariance:det:RKHS}: the kernel is Hermitian, $k(x,y)=\overline{k(y,x)}$, and positive definite, and the canonical features are $\phi_x(y) = k(x,y)$.%

\subsection{Deterministic systems}

\subsubsection{Self-adjointness on the space \texorpdfstring{$L^2$}{L2}}
On $L^2(\mu)$, self-adjointness is equivalent to invariance of the measure together with involutivity of the dynamics. Neither boundedness nor nonsingularity is assumed separately; the two stated conditions imply that the composition operator is an isometry. For bounded composition operators on $L^2$ of a $\sigma$-finite measure, this characterization is classical; see~\cite[Theorem~3]{singh1975normal}. Theorem~\ref{t:det:selfadjL2} restates it in the present setting, and we include the short proof for completeness.

\begin{thm}\label{t:det:selfadjL2}
The formula $\calK f = f\circ F$ defines a bounded self-adjoint operator on $L^2(\mu)$ if and only if \begin{enumerate}
\item[\rm (i)] $\mu\circ F^{-1}=\mu$;
\item[\rm (ii)] $F^2=\operatorname{id}$ $\mu$-almost everywhere.
\end{enumerate}
\end{thm}
\begin{proof} Assume first that \rm{(i)} and \rm{(ii)} hold. Invariance (i) makes composition independent of the representative and yields
\begin{align*}
\|\calK f\|_2^2=\int |f|^2\,d(\mu\circ F^{-1})=\|f\|_2^2, \qquad f\in L^2(\mu),
\end{align*}
that is, $\calK^*\calK = I$. Condition \rm{(ii)} gives $\calK^2=I$. Hence, $\calK^* = \calK^*\calK^2 = \calK$.

Conversely, suppose that the formula defines a bounded self-adjoint operator. Choose Borel sets $X_m\uparrow\calX$ with $\mu(X_m)<\infty$. Testing self-adjointness on indicators of finite-measure sets $A,B\in\calB_\calX$ gives
\begin{align}\label{e:det_detailed_balance}
\mu\bigl(A\cap F^{-1}(B)\bigr) &=\langle\calK\one_B,\one_A\rangle =\langle\one_B,\calK\one_A\rangle =\mu\bigl(B\cap F^{-1}(A)\bigr).
\end{align}
Fix $B$ of finite measure and take $A=X_m$ in \eqref{e:det_detailed_balance}. Monotone convergence then gives
\begin{align*}
\mu(F^{-1}(B)) =\lim_m\mu\bigl(X_m\cap F^{-1}(B)\bigr) =\lim_m\mu\bigl(B\cap F^{-1}(X_m)\bigr)=\mu(B).
\end{align*}
Apply this identity to $B\cap X_m$ and pass once more to the limit. The result holds for every Borel set $B$, so $\mu$ is invariant and $\calK$ is an isometry. Self-adjointness now implies $\calK^2=\calK^*\calK=I$.

It remains to recover the pointwise assertion about $F$ from this operator identity. Let $(U_j)_{j\ge1}$ be a countable base of the Polish space $\calX$. The countable family
\begin{align*}
\mathscr A:=\{U_j\cap X_m:j,m\ge1\}
\end{align*}
consists of finite-measure sets and separates the points of $\calX$. For each $A\in\mathscr A$, the equality $\calK^2\one_A=\one_A$ supplies a null set $N_A$ outside of which $\one_A(F^2(x))=\one_A(x)$. All these equalities hold simultaneously outside $N=\bigcup_{A\in\mathscr A}N_A$. Since $\mathscr A$ separates points, $F^2(x)=x$ for every $x\in\calX\backslash N$.
\end{proof}

\begin{cor}\label{c:det:selfadj:spectral}
If $\calK f = f\circ F$ defines a bounded self-adjoint operator on $L^2(\mu)$, it is an involution, $\sigma(\calK)\subset\{-1,1\}$,
\begin{align*}
L^2(\mu)=L^2_+(\mu)\oplus L^2_-(\mu),\qquad L^2_\pm(\mu):=\{f\in L^2(\mu):f\circ F=\pm f\ \mu\text{-a.e.}\},
\end{align*}
and the orthogonal projections onto these two subspaces are $P_\pm=(I\pm\calK)/2$.
\end{cor}
\begin{proof} For the self-adjoint involution $\calK$, the operators $P_\pm$ satisfy $P_\pm^*=P_\pm$, $P_\pm^2=P_\pm$, $P_+P_-=0$, and $P_++P_-=I$. Their ranges are precisely the indicated eigenspaces. Moreover, $(\calK-\lambda I)^{-1}=(\calK+\lambda I)/(1-\lambda^2)$ whenever $\lambda\notin\{-1,1\}$, which proves the spectral assertion.
\end{proof}

\subsubsection{Self-adjointness on reproducing kernel Hilbert spaces}
The relevant realization is the maximal closed composition operator $\calK_\calH:\calD_\calH\subset\calH\to\calH$ from Subsection~\ref{s:invariance:det:RKHS}. Here, self-adjointness can be read directly from the action of the adjoint on the canonical features.

\begin{thm}\label{t:selfadjprep}
If $\calK_\calH$ is self-adjoint, then
\begin{align}\label{e:rkhs:sa}
k(x,F(y))=k(F(x),y)\qquad(x,y\in\calX).
\end{align}
Conversely, if $\calH$ is Koopman invariant and \eqref{e:rkhs:sa} holds, then $\calK_\calH$ is self-adjoint.
\end{thm}
\begin{proof} We first derive an adjoint identity under the sole assumption that $\calD_\calH$ is dense. For $f\in\calD_\calH$, the reproducing property gives
\begin{align*}
\langle\calK_\calH f,\phi_x\rangle =f(F(x))=\langle f,\phi_{F(x)}\rangle.
\end{align*}
Hence
\begin{equation}\label{e:det_feature_adjoint}
\phi_x\in\dom(\calK_\calH^*)
\qquad\text{and}\qquad
\calK_\calH^*\phi_x=\phi_{F(x)}.
\end{equation}
If $\calK_\calH$ is self-adjoint, its domain is dense, so the preceding identity applies and shows that $\phi_x\in\calD_\calH$ and $\calK_\calH\phi_x=\phi_{F(x)}$. Evaluating this equality at $y$ yields
\begin{align*}
k(x,F(y)) = \phi_x(F(y)) = (\calK_\calH\phi_x)(y) = \phi_{F(x)}(y)=k(F(x),y).
\end{align*}
Conversely, assume Koopman invariance and \eqref{e:rkhs:sa}. Lemma~\ref{l:RKHS_invar_det} makes $\calK_\calH$ bounded and everywhere defined. Equations \eqref{e:rkhs:sa} and \eqref{e:det_feature_adjoint} give $\calK_\calH\phi_x=\calK_\calH^*\phi_x$ for every $x$. The canonical features have dense linear span; continuity therefore extends this equality to all of $\calH$.
\end{proof}

As before, let $k_F(x,y)=k(F(x),F(y))$. Under the strict separation hypothesis below, the RKHS criterion has the same two-part structure as Theorem~\ref{t:det:selfadjL2}: preservation of the ambient structure and involutivity. Here, preservation is expressed by $k_F=k$.

\begin{cor}\label{c:rkhs:sa}
Let $\calX\subset\R^n$ and let $k(x,y)=\Phi(x-y)$, where $\Phi:\R^n\to\C$ is positive definite, $\Phi(0)>0$, and
\begin{align*}
|\Phi(z)|<\Phi(0)\qquad(z\ne0).
\end{align*}
Then $\calK_\calH$ is self-adjoint if and only if \begin{enumerate}
\item[\rm (i)]  $k_F=k$;
\item[\rm (ii)] $F^2=\operatorname{id}$ pointwise on $\calX$.
\end{enumerate}
In that event, $\calK_\calH$ is a bounded self-adjoint involution, $\sigma(\calK_\calH)\subset\{-1,1\}$,
\begin{align*}
\calH = \calH_+\oplus\calH_-,\qquad \calH_\pm:=\{f\in\calH:f\circ F=\pm f\},
\end{align*}
and the corresponding orthogonal projections are $(I\pm\calK_\calH)/2$.
\end{cor}
\begin{proof} Suppose first that $\calK_\calH$ is self-adjoint. Substitution of $y=F(x)$ in \eqref{e:rkhs:sa}, together with translation invariance, yields
\begin{align*}
\Phi(x-F^2(x))=\Phi(0),
\end{align*}
and the strict inequality in the hypothesis forces $F^2(x)=x$ for every $x$. Replacing $y$ by $F(y)$ in \eqref{e:rkhs:sa} now gives
\begin{align*}
k_F(x,y)=k(F(x),F(y))=k(x,F^2(y))=k(x,y).
\end{align*}
Conversely, assume \rm{(i)} and \rm{(ii)}. The identity $k_F=k$ implies Koopman invariance by Theorem~\ref{t:RKHS_det_inv}, whereas
\begin{align*}
k(x,F(y)) = k(F(F(x)),F(y)) = k_F(F(x),y) = k(F(x),y)
\end{align*}
proves \eqref{e:rkhs:sa}. Theorem~\ref{t:selfadjprep} therefore gives self-adjointness. The remaining assertions follow as in Corollary~\ref{c:det:selfadj:spectral}.
\end{proof}

\begin{rem}
The separation hypothesis cannot be omitted. For example, take $\calX=\R$, $k(x,y)=1$, and $F(x)=x+1$. The kernel is translation invariant and $\calK_{\calH} = I$, but $F^2\ne\operatorname{id}$.
\end{rem}

\subsection{Markov kernels}
\subsubsection{Self-adjointness on the space \texorpdfstring{$L^2$}{L2}}
For a general Markov kernel, the transition mechanism need not be represented by a single map. The appropriate symmetry is instead carried by the joint measure of two successive states, namely the $\sigma$-finite Borel measure $Q$ on $\calX\times\calX$ given by
\begin{align}\label{e:joint_measure}
Q(dx,dy) := \rho_x(dy)\mu(dx), \qquad Q(E) = \int\!\int\one_E(x,y)\,d\rho_x(y)\,d\mu(x).
\end{align}
Its first marginal is $\mu$, and its second marginal is the measure $\nu$ defined in \eqref{e:nu}. Whenever the Koopman formula induces an operator on $L^2(\mu)$, its matrix coefficients on indicators of finite-measure Borel sets satisfy
\begin{align*}
Q(A\times B)=\langle\calK\one_B,\one_A\rangle.
\end{align*}
This identity singles out interchange of the initial and final states as the relevant symmetry. If $\vartheta(x,y)=(y,x)$, we call the kernel \emph{reversible with respect to $\mu$} when
\begin{align*}
Q=\vartheta_\#Q.
\end{align*}
Equivalently, reversibility means that the \emph{detailed-balance} identities
\begin{equation}\label{e:detailed_balance_rectangles}
\int_A\rho_x(B)\,d\mu(x) =\int_B\rho_x(A)\,d\mu(x) \qquad(A,B\in\calB_\calX)
\end{equation}
hold. This is the usual reversibility condition for Markov dynamics; see, for example, \cite{noenusk13}.

\begin{thm}\label{t:sa_stoch}
The following statements are equivalent.
\begin{enumerate}
\item[\rm (i)] The joint measure $Q$ in \eqref{e:joint_measure} is symmetric: $Q=\vartheta_\#Q$.
\item[\rm (ii)] $L^2(\mu)$ is Koopman invariant, and $\calK : L^2(\mu)\to L^2(\mu)$ is self-adjoint.
\end{enumerate}
Under these equivalent conditions, $\nu=\mu$, $\calK$ is a contraction, $\sigma(\calK)\subset[-1,1]$, and
\begin{align}\label{e:dirichlet_identity}
\langle(I-\calK)f,f\rangle_{L^2(\mu)} =\frac12\int_{\calX\times\calX}|f(x)-f(y)|^2\,dQ(x,y) \qquad(f\in L^2(\mu)).
\end{align}
Moreover, the eigenspaces corresponding to $1$ and $-1$ admit the intrinsic descriptions
\begin{align}
\begin{split}\label{e:stoch_extremal_eigenspaces}
\ker(I-\calK) &=\{f\in L^2(\mu):f(x)=f(y)\text{ for $Q$-a.e.\ }(x,y)\},\\ \ker(I+\calK) &=\{f\in L^2(\mu):f(x)=-f(y)\text{ for $Q$-a.e.\ }(x,y)\}.
\end{split}
\end{align}
\end{thm}
\begin{proof} Assume first that $Q$ is symmetric. Its marginals are then equal, so $\nu=\mu$. By Theorem \ref{t:invariance:stoch:Lp}, $L^2(\mu)$ is Koopman invariant, and $\calK : L^2(\mu)\to L^2(\mu)$ is a contraction. For $f,g\in L^2(\mu)$, the function $(x,y)\mapsto f(y)\overline{g(x)}$ is $Q$-integrable by the Cauchy--Schwarz inequality and the marginal identities. Symmetry of $Q$ now gives
\begin{align*}
\langle\calK f,g\rangle &=\int f(y)\overline{g(x)}\,dQ(x,y) =\int f(x)\overline{g(y)}\,dQ(x,y) =\langle f,\calK g\rangle.
\end{align*}
Hence, $\calK$ is self-adjoint.

Conversely, suppose that (ii) holds. Testing on indicators of sets $A,B$ of finite $\mu$-measure gives
\begin{align*}
Q(A\times B)=\langle\calK\one_B,\one_A\rangle =\langle\one_B,\calK\one_A\rangle=Q(B\times A).
\end{align*}
For arbitrary $A,B$, apply this equality to $A\cap X_m$ and $B\cap X_m$, where $X_m\uparrow\calX$ and $\mu(X_m)<\infty$, and pass to the limit by monotone convergence. Thus $Q$ and $\vartheta_\#Q$ agree on all measurable rectangles. Both are $\sigma$-finite, and the rectangles $X_m\times X_m$ form a common finite-measure exhaustion. Uniqueness on the product $\sigma$-algebra therefore gives $Q=\vartheta_\#Q$.

A self-adjoint contraction clearly has its spectrum in $[-1,1]$. To identify the two eigenspaces, use the marginal identities and expand the square:
\begin{align*}
\frac12\int|f(x)-f(y)|^2\,dQ &=\|f\|_2^2-\operatorname{Re}\int f(y)\overline{f(x)}\,dQ\\ &=\|f\|_2^2-\langle\calK f,f\rangle = \langle(I-\calK)f,f\rangle,
\end{align*}
where $\langle\calK f,f\rangle$ is real because $\calK$ is self-adjoint. The same calculation with a plus sign gives
\begin{align}\label{e:dirichlet_identity_plus}
\langle(I+\calK)f,f\rangle = \frac12\int|f(x)+f(y)|^2\,dQ.
\end{align}
Since $\sigma(\calK)\subset[-1,1]$, the operators $I-\calK$ and $I+\calK$ are non-negative. A non-negative operator $T$ satisfies $\langle Tf,f\rangle=0$ if and only if $f\in\ker T$. This proves the identities in \eqref{e:stoch_extremal_eigenspaces}.
\end{proof}

For a deterministic kernel $\rho_x=\delta_{F(x)}$, the joint measure is $Q=(\operatorname{id},F)_\#\mu$. Theorem~\ref{t:sa_stoch} then reduces to Theorem~\ref{t:det:selfadjL2}: symmetry of the graph measure is equivalent to invariance of $\mu$ and $F^2=\operatorname{id}$ almost everywhere.

\begin{ex}[Overdamped Langevin dynamics]\label{ex:olan:selfadj}
Consider again the setting of Example~\ref{ex:olan}, where it was shown that the associated Koopman operators $\calK_t$ form a strongly continuous contraction semigroup on each $L^p(\mu)$, $1\le p<\infty$, with invariant Gibbs measure $\mu$. Moreover, both the infinitesimal generator $\mathcal L$ and all $\calK_t$, $t\ge0$, are self-adjoint on $L^2(\mu)$. Equivalently, by Theorem~\ref{t:sa_stoch}, the measure $\mu(dx)\rho_t(x,dy)$ satisfies detailed balance.

To illustrate this property explicitly, consider the quadratic potential $V(x)=\tfrac12x^\top Ax$, where $A=A^\top>0$. Then $\mu=\calN(0,\Sigma)$, where $\Sigma=(\beta A)^{-1}$. The solution of
\eqref{eq:langevin} is given by
\begin{align*}
X_t^x = B_tx+\sqrt{\frac{2D}{\beta}} \int_0^t B_{t-s}\,dW_s,\qquad B_t=e^{-DtA}.
\end{align*}
Consequently, the transition probabilities are
\begin{align*}
\rho_t(x,\cdot) = \calN\bigl(B_tx,\Sigma-B_t\Sigma B_t\bigr).
\end{align*}
If $X\sim\calN(0,\Sigma)$ and, conditionally on $X=x$, the variable $Y$ has this transition law, then $(X,Y)$ is centered Gaussian with covariance matrix
\begin{align*}
\begin{pmatrix}
\Sigma & \Sigma B_t\\
B_t\Sigma & \Sigma
\end{pmatrix}.
\end{align*}
Since $B_t$ and $\Sigma$ are symmetric and commute, this covariance matrix is invariant under interchange of $X$ and $Y$. Hence, $\mu(dx)\rho_t(x,dy)$ is symmetric, and $\calK_t$ is a self-adjoint contraction by Theorem~\ref{t:sa_stoch}.
\end{ex}

\begin{rem}
A self-adjoint stochastic Koopman operator on $L^2(\mu)$ with invariant probability measure $\mu$ need not be nonnegative. For instance, if $\mu=\mathcal N(0,1)$ and $F(x)=-x$, then the Koopman operator $\calK f=f\circ F$ is self-adjoint and satisfies $\calK f=-f$ for $f(x)=x$.

In contrast, if $(\calK_t)_{t\ge0}$ is a semigroup of bounded self-adjoint operators on $L^2(\mu)$, then
\begin{align*}
\langle\calK_t f,f\rangle = \langle\calK_{t/2}^2f,f\rangle = \|\calK_{t/2}f\|^2 \ge 0
\end{align*}
for every $f\in L^2(\mu)$ and $t\ge0$. Thus, every operator in a self-adjoint semigroup is necessarily nonnegative.
\end{rem}

\subsubsection{Self-adjointness on reproducing kernel Hilbert spaces}
On an RKHS, self-adjointness is detected by the averaged canonical features. Under the bounded-diagonal Hermitian hypotheses of Subsection~\ref{s:invariance:stoch:RKHS}, the vectors $\Psi_x\in\calH$ are characterized by \eqref{eq:Psi_Riesz}; equivalently, they are Pettis integrals of the feature map.

\begin{lem}\label{l:K*psi}
If $\calD_\calH$ is dense in $\calH$, i.e., $\calK_\calH$ is densely defined, then $\phi_x\in\dom(\calK_\calH^*)$ with $\calK_\calH^*\phi_x=\Psi_x$. In particular, $\calK_\calH^*$ is also densely defined.
\end{lem}
\begin{proof} For all $f\in\calD_\calH$ we have
\begin{align*}
\langle\calK_\calH f,\phi_x\rangle_\calH =(\calK f)(x)=\int f(y)\,d\rho_x(y) =\langle f,\Psi_x\rangle_\calH.
\end{align*}
Thus $\phi_x\in\dom(\calK_\calH^*)$ and $\calK_\calH^*\phi_x=\Psi_x$. Since the span of the $\phi_x$ is dense, the adjoint is densely defined.
\end{proof}

\begin{thm}\label{t:selfadj:stoch:rkhs}
If $\calK_\calH$ is self-adjoint, then
\begin{align}\label{e:RKHS_self-adjoint}
\Psi_x(y)=\overline{\Psi_y(x)} \qquad(x,y\in\calX).
\end{align}
Conversely, if $\calH$ is Koopman invariant and \eqref{e:RKHS_self-adjoint} holds, then $\calK_\calH$ is self-adjoint.
\end{thm}
\begin{proof}
Suppose that $\calK_\calH$ is self-adjoint. Then its domain is dense, and Lemma~\ref{l:K*psi} gives
\begin{equation}\label{e:Kfeature_hermitian}
\phi_x\in\dom(\calK_\calH^*)
\qquad\text{and}\qquad
\calK_\calH^*\phi_x=\Psi_x.
\end{equation}
Consequently, $\calK_\calH\phi_x=\Psi_x$. The defining identity \eqref{eq:Psi_Riesz} and Hermitian symmetry give
\begin{align*}
\Psi_x(y) = (\calK_\calH\phi_x)(y) = \int\phi_x(z)\,d\rho_y(z) = \langle\phi_x,\Psi_y\rangle_\calH = \overline{\Psi_y(x)},
\end{align*}
which proves \eqref{e:RKHS_self-adjoint}.

Conversely, assume Koopman invariance and \eqref{e:RKHS_self-adjoint}. Lemmas~\ref{l:RKHS_invar_stochastic} and \ref{l:K*psi} make $\calK_\calH$ bounded and everywhere defined and give \eqref{e:Kfeature_hermitian}. For $x,y\in\calX$, the kernel identity yields
\begin{align*}
(\calK_\calH\phi_x)(y) = \int \phi_x(z)\,d\rho_y(z) = \overline{\Psi_y(x)}=\Psi_x(y).
\end{align*}
Together with \eqref{e:Kfeature_hermitian}, this implies $\calK_\calH\phi_x=\Psi_x=\calK_\calH^*\phi_x$. The canonical features have dense linear span in $\calH$, so continuity gives $\calK_\calH=\calK_\calH^*$.
\end{proof}

For a deterministic kernel $\rho_x=\delta_{F(x)}$, one has $\Psi_x=\phi_{F(x)}$, and \eqref{e:RKHS_self-adjoint} reduces to \eqref{e:rkhs:sa}. Thus the deterministic and stochastic RKHS criteria express the same symmetry of transported features.

\section*{Acknowledgments} We thank the Isaac Newton Institute for Mathematical Sciences, Cambridge, for its support and hospitality during the program \textit{Operator approaches to dynamics: new connections}, where part of this work was undertaken. This work was partially supported by EPSRC grant EP/Z000580/1. Karl Worthmann is thankful for support by the  Deutsche Forschungsgemeinschaft (DFG, German
Research Foundation) within the research unit ALeSCo – Active Learning in Systems and Control, project-ID 535860958.

\appendix
\renewcommand{\thesection}{\Alph{section}}

\section{Auxiliary measure and RKHS results}\label{a:aux} This appendix collects the measurability of kernel integrals, counterexamples to representative independence, and the range-space description of $\calH_\rho$.

\begin{lem}[Measurability of kernel integrals]\label{l:measurable}
For every Borel function $h:\calX\to[0,\infty]$, the extended-valued function
\begin{align}\label{e:I_h}
I_h(x):=\int h(y)\,d\rho_x(y),\qquad x\in\calX,
\end{align}
is Borel measurable. Consequently, for every $f\in\calM(\calX)$,
\begin{align*}
\Delta_f=\left\{x\in\calX:I_{|f|}(x)<\infty\right\}\in\calB_\calX,
\end{align*}
and the zero-extended function $\calK_0f$ defined in \eqref{eq:zero_extended_K} is Borel measurable.
\end{lem}
\begin{proof}
Choose nonnegative simple functions $h_n\nearrow h$. The defining properties of a Markov kernel make each $I_{h_n}$ Borel measurable, while monotone convergence gives $I_{h_n}\nearrow I_h$. Hence $I_h$ is Borel measurable. Taking $h=|f|$ then yields
\begin{align*}
\Delta_f=\bigcup_{m=1}^\infty\{I_{|f|}\le m\}\in\calB_\calX.
\end{align*}
For the zero extension, consider the truncations $f_n=f\one_{\{|f|\le n\}}$. Each $f_n$ is bounded and Borel measurable. Separating the real and imaginary parts into their positive and negative parts shows that $x\mapsto\int f_n\,d\rho_x$ is Borel measurable. If $x\in\Delta_f$, dominated convergence gives
\begin{align*}
\int f_n\,d\rho_x\longrightarrow\int f\,d\rho_x.
\end{align*}
Consequently,
\begin{align*}
\calK_0f =\lim_{n\to\infty}\one_{\Delta_f} \left(x\mapsto\int f_n\,d\rho_x\right)
\end{align*}
pointwise, proving the final assertion.
\end{proof}

Measurability alone does not ensure the compatibility of $f\mapsto\calK f$ with equality $\mu$-almost everywhere.

\begin{ex}\label{ex:equivalence_classes:counterexample}%
The first two examples are deterministic: they show that $g\sim 0$ need not imply $\calK g\sim 0$. The third is stochastic and exhibits the stronger obstruction that $g\sim0$ need not imply $g\in\calD(\mu)$.

\smallskip\noindent (a) Let $\calX=\R$, $\mu=\delta_0$, and $F(x) = x+1$ for $x\in\calX$. For $g=\one_{\{1\}}$, one has $g\sim0$ but $\calK g=\one_{\{0\}}\not\sim0$. Hence $[0]\notin\wt\calD(\mu)$.

\medskip\noindent (b) Let $\calX = (-\frac 12,\frac 12)$, let $\mu$ be Lebesgue measure, and let $F(x)=0$ for $x\in\calX$. For $g=\one_{\{0\}}$, one again has $g\sim0$ but $\calK g=\one\not\sim0$. Thus $[0]\notin\wt\calD(\mu)$.

\medskip\noindent (c) Let $\calX = (0,1)$, let $\mu$ be Lebesgue measure, and take $\rho_x = \sum_{i=1}^\infty2^{-i}\delta_{x_i}$ for every $x\in\calX$, where $\{x_i : i\in\N\}$ enumerates $\calX\cap\Q$. Consider $g_1=\sum_{i=1}^\infty\one_{\{x_i\}}$ and $g_2=\sum_{i=1}^\infty 2^i\one_{\{x_i\}}$. Both satisfy $g_1\sim g_2\sim0$. Although $g_1\in\calD(\mu)$, one has $(\calK g_1)(x)=1$ for every $x\in\calX$. Hence $\calK g_1\not\sim0$ and $[0]\notin\wt\calD(\mu)$. By contrast, $g_2\notin\calD(\mu)$ because $\Delta_{g_2}=\varnothing$.
\end{ex}

The remaining auxiliary result identifies the range RKHS used in the invariance and compactness criteria.

\begin{lem}[The range RKHS]\label{l:Hrho}
Assume the bounded-diagonal Hermitian RKHS hypotheses of Subsection~\rmref{s:invariance:stoch:RKHS}. Let $S:\calH\to\calM(\calX)$ be the pointwise Markov operator, $(Sf)(x)=\int f\,d\rho_x$, and let $k_\rho$ be defined by \eqref{e:krho}. Then the RKHS induced by $k_\rho$ is
\begin{align*}
\calH_\rho=S\calH,
\end{align*}
with range norm
\begin{align*}
\|g\|_{\calH_\rho} =\inf\{\|f\|_\calH:f\in\calH,\ Sf=g\}.
\end{align*}
Moreover, if $M=\overline{\linspan}\{\Psi_x:x\in\calX\}$, then $M=(\ker S)^\perp$, and $S|_M : M\longrightarrow\calH_\rho$ is an isometric isomorphism.
\end{lem}
\begin{proof} Let $\phi_x^\rho=k_\rho(x,\cdot)$. For $a_1,\ldots,a_m\in\C$ and $x_1,\ldots,x_m\in\calX$,
\begin{align*}
\left\|\sum_{i=1}^m a_i\Psi_{x_i}\right\|_\calH^2 =\sum_{i,j=1}^m a_i\overline{a_j}k_\rho(x_i,x_j) =\left\|\sum_{i=1}^m a_i\phi_{x_i}^\rho \right\|_{\calH_\rho}^2.
\end{align*}
The rule $T\phi_x^\rho=\Psi_x$ therefore extends uniquely to an isometry $T:\calH_\rho\to\calH$ with range $M$. For every $f\in\calH$ and $x\in\calX$, the reproducing identity and the definition of the adjoint give
\begin{align*}
(T^*f)(x) =\langle T^*f,\phi_x^\rho\rangle_{\calH_\rho} =\langle f,T\phi_x^\rho\rangle_\calH =\langle f,\Psi_x\rangle_\calH =\int f\,d\rho_x =(Sf)(x).
\end{align*}
Thus the pointwise operator $S$ agrees with the bounded operator $T^*:\calH\to\calH_\rho$. Since $T$ is an isometry,
\begin{align*}
ST=T^*T=I_{\calH_\rho},\qquad S\calH=\calH_\rho, \qquad \ker S=M^\perp, \qquad TS=TT^*=P_M.
\end{align*}
In particular, $S|_M$ is the isometric inverse of $T$. If $Sf=g$, then $f=Tg+(f-Tg)$ with $f-Tg\in\ker S=M^\perp$. Hence $f=Tg$ is the unique solution of least norm, and
\begin{align*}
\|g\|_{\calH_\rho}=\|Tg\|_\calH =\inf_{Sf=g}\|f\|_\calH,
\end{align*}
which concludes the proof.
\end{proof}

\section{Domination of measures and bounded densities} \label{a:measure_domination}

The following elementary lemma identifies domination by a constant multiple of a measure with an essentially bounded Radon--Nikodym derivative.

\begin{lem}\label{l:measure_domination} Let $\mu_1$ and $\mu_2$ be measures on a measurable space $(Y,\Sigma)$, where $\mu_2$ is $\sigma$-finite, and let $0\le C<\infty$. Then the following are equivalent: \begin{enumerate} \item[{\rm (i)}] $\mu_1(\Delta)\le C\mu_2(\Delta)$ for every $\Delta\in\Sigma$; \item[{\rm (ii)}] $\mu_1\ll\mu_2$ and $g=d\mu_1/d\mu_2\in L^\infty(\mu_2)$ with $\|g\|_{L^\infty(\mu_2)}\le C$.
\end{enumerate}
\end{lem}
\begin{proof} Under (ii), for every $\Delta\in\Sigma$,
\begin{align*}
\mu_1(\Delta) &=\int\one_\Delta\,d\mu_1 =\int\one_\Delta g\,d\mu_2 \le C\int\one_\Delta\,d\mu_2 =C\mu_2(\Delta).
\end{align*}
This is (i).

Conversely, the case $C=0$ is immediate. Assume (i) with $C>0$. It follows directly that $\mu_1\ll\mu_2$. Moreover, $\mu_1$ is $\sigma$-finite, so the Radon--Nikodym theorem supplies a density $g = d\mu_1/d\mu_2\ge 0$. Choose mutually disjoint sets $\Delta_1,\Delta_2,\ldots\in\Sigma$ such that $\mu_2(\Delta_j) < \infty$ for each $j\in\N$ and $Y = \bigcup_{j=1}^\infty\Delta_j$. For $n\in\N$, let $\Delta_n' := \{x\in Y : g(x)\ge C + 1/n\}$. Then
\begin{align*}
\begin{aligned} C\mu_2(\Delta_n'\cap\Delta_j) &\ge \mu_1(\Delta_n'\cap\Delta_j) =\int_{\Delta_n'\cap\Delta_j}g\,d\mu_2 \ge (C+\tfrac1n)\mu_2(\Delta_n'\cap\Delta_j).
\end{aligned}
\end{align*}
Hence $\mu_2(\Delta_n'\cap\Delta_j) = 0$ for all $j\in\N$, and $\mu_2(\Delta_n') = \sum_{j=1}^\infty\mu_2(\Delta_n'\cap\Delta_j) = 0$ for all $n\in\N$. It follows that $g\le C$ $\mu_2$-almost everywhere, and therefore $\|g\|_{L^\infty(\mu_2)}\le C$.
\end{proof}

\section{Integration with respect to the measure \texorpdfstring{$\nu$}{nu}}\label{a:integration_nu} The measure $\nu$ is the second marginal of the one-step measure $\mu(dx)\rho_x(dy)$. The following kernel version of Tonelli--Fubini makes this relation explicit. For deterministic dynamics it reduces to the familiar equivalence $f\in L^1(\mu\circ F^{-1})$ if and only if $f\circ F\in L^1(\mu)$, together with $\int f\,d(\mu\circ F^{-1}) = \int f\circ F\,d\mu$.

\begin{lem}[Kernel Tonelli--Fubini formula]\label{l:nu}
For every Borel function $h:\calX\to[0,\infty]$,
\begin{align}\label{e:murhox_nonnegative}
\int h\,d\nu = \int\left(\int h(y)\,d\rho_x(y)\right)d\mu(x),
\end{align}
with equality in $[0,\infty]$. Consequently, for every complex Borel function $f$, the following are equivalent: \begin{enumerate}
\item[\rm(i)] $f\in L^1(\nu)$;
\item[\rm(ii)] $I_{|f|}\in L^1(\mu)$, cf.\ \eqref{e:I_h};
\item[\rm(iii)] $f\in\calD(\mu)$ and $\calK_0|f|\in L^1(\mu)$.
\end{enumerate}
Whenever these conditions hold,
\begin{align}\label{e:murhox}
\int f\,d\nu = \int\left(\int f(y)\,d\rho_x(y)\right)d\mu(x) = \int\calK_0f\,d\mu.
\end{align}
\end{lem}
\begin{proof}
The definition of $\nu$ proves \eqref{e:murhox_nonnegative} for indicator functions. Linearity extends it to nonnegative simple functions, and monotone convergence then gives the identity for every nonnegative Borel function. Taking $h=|f|$ proves the equivalence of (i) and (ii). The same identity shows that (ii) is equivalent to $I_{|f|}<\infty$ almost everywhere together with $\calK_0|f|=I_{|f|}$ $\mu$-almost everywhere and $\calK_0|f|\in L^1(\mu)$; this is precisely (iii). Under these conditions, the inner integral of $f$ is absolutely convergent for almost every $x$. Apply \eqref{e:murhox_nonnegative} to the positive and negative parts of the real and imaginary parts of $f$ to obtain \eqref{e:murhox}.
\end{proof}

\section{A countable norming class for total variation}\label{a:tv}
The total-variation norm is defined by a supremum over all Borel sets. On a Polish space, the same supremum may be taken over a fixed countable class. This countability is essential in Theorem~\ref{t:compact:stoch:Linf}, where it permits a single exceptional null set for all indicator tests. Fix a complete metric $d$ that induces the topology of $\calX$; such a metric may be chosen because $\calX$ is Polish.

\begin{thm}\label{t:TV_with_calE}
Fix a complete compatible metric $d$ on $\calX$, let $\{x_i : i\in I\}\subset\calX$ be countable and dense, and let $\calE_1$ consist of all balls $B_r(x_i)$ with $r\in (0,\infty)\cap\Q$ and $i\in I$. Then the collection $\calE$ of all finite unions of members of $\calE_1$ is countable, and, for every difference $\tau = \mu_1 - \mu_2$ of probability measures $\mu_1$ and $\mu_2$ on $(\calX,\calB_\calX)$,
\begin{align}\label{e:norming_class}
\|\tau\|_{TV} = 2\,\sup_{E\in\calE}|\tau(E)|.
\end{align}
\end{thm}
\begin{proof}
For $m\geq1$, let $\calE_m$ be the collection of unions of $m$ members of $\calE_1$. Each $\calE_m$ is countable, and hence so is $\calE=\bigcup_{m\geq1}\calE_m$.

Let $p=\mu_1+\mu_2$. Since $\calX$ is Polish, $p$ is a finite Radon measure. Given $\Delta\in\calB_\calX$ and $\varepsilon>0$, outer regularity supplies an open $U\supset\Delta$ with $p(U\backslash\Delta)<\varepsilon$. Therefore,
\begin{align*}
|\tau(\Delta)|\leq |\tau(U)|+p(U\backslash\Delta) <|\tau(U)|+\varepsilon.
\end{align*}
Given an open $U$, inner regularity supplies a compact $K\subset U$ with $p(U\backslash K)<\varepsilon$, and similarly $|\tau(U)|<|\tau(K)|+\varepsilon$. Consequently,
\begin{align}\label{e:TV_open_compact_reduction}
\frac12\|\tau\|_{TV} &=\sup_{\Delta\in\calB_\calX}|\tau(\Delta)| =\sup_{U\subset\calX\ \mathrm{open}}|\tau(U)| =\sup_{K\subset\calX\ \mathrm{compact}}|\tau(K)|.
\end{align}
It remains to approximate compact sets from outside by members of $\calE$. Fix a compact $K$ and $\varepsilon>0$. Outer regularity gives an open $U\supset K$ with $p(U\backslash K)<\varepsilon$. The rational balls centered at $\{x_i:i\in I\}$ form a base for the topology. Indeed, for $x\in U$, choose $r_0>0$ such that $B_{r_0}(x)\subset U$, then choose $x_i$ with $d(x,x_i)<r_0/3$ and a rational $r$ satisfying
\begin{align*}
d(x,x_i)<r<r_0-d(x,x_i).
\end{align*}
Then $x\in B_r(x_i)\subset B_{r_0}(x)\subset U$. Thus members of $\calE_1$ contained in $U$ cover $K$. A finite subcover has union $E\in\calE$ satisfying $K\subset E\subset U$, and
\begin{align*}
|\tau(K)|\leq |\tau(E)|+|\tau(E)-\tau(K)| \leq |\tau(E)|+p(E\backslash K) <|\tau(E)|+\varepsilon.
\end{align*}
Taking suprema and applying \eqref{e:TV_open_compact_reduction} proves ``$\le$'' in \eqref{e:norming_class}. The reverse inequality follows from $\calE\subset\calB_\calX$.
\end{proof}

\begin{cor}\label{c:meas_map_meas} Let $\Delta\in\calB_\calX$ and equip $M = \{\rho_x : x\in\Delta\}$ with the total variation metric. If $M$ is separable, then the map $\Theta : \Delta\to M$ defined by $\Theta(x) = \rho_x$ is Borel measurable.
\end{cor}
\begin{proof} Since $M$ is separable, it suffices to prove that $\Theta^{-1}(\ol B_r(\rho_{x_0}))\in\calB_\calX$ for every $r>0$ and $x_0\in\Delta$, where $\ol B_r(\rho_{x_0})$ is the closed ball in $M$ with center $\rho_{x_0}$ and radius $r$. Let $\calE$ be the countable family from Theorem~\ref{t:TV_with_calE}. Its norming property gives
\begin{align*}
\Theta^{-1}(\ol B_r(\rho_{x_0})) =\bigcap_{E\in\calE} \left\{x\in\Delta: |\rho_x(E)-\rho_{x_0}(E)|\le\frac r2\right\}.
\end{align*}
The last set is Borel because $\calE$ is countable and $x\mapsto\rho_x(E)$ is measurable for every $E\in\calE$.
\end{proof}

\section{Koopman invariance of the Gaussian RKHS}\label{a:GaussianRKHS}
For $\sigma>0$, consider the Gaussian kernel
\begin{align*}
k(x,y) = \exp\left(-\frac{\|x-y\|^2}{2\sigma^2}\right),\qquad x,y\in\R^n,
\end{align*}
and denote the corresponding RKHS on a set $\calX\subset\R^n$ by $\calH_\sigma(\calX)$. Let $F:\calX\to\calX$ be arbitrary.

\begin{thm}[Gaussian RKHS invariance and exact norm]\label{t:Gaussian_exact_norm}
Let $\sigma > 0$, and assume that $\calX$ has non-empty interior. Then the Gaussian RKHS $\calH = \calH_\sigma(\calX)$ is Koopman invariant for deterministic systems $x^+ = F(x)$ if and only if
\begin{align*}
F(x)=Ax+b,
\end{align*}
where $b\in\R^n$ and $A\in\R^{n\times n}$ is invertible with $\|A\|_{2\to2}\le 1$. In this case,
\begin{align}\label{e:gaussian_norm}
\|\calK_\calH\|_{\calH\to\calH}=|\det A|^{-1/2}.
\end{align}
\end{thm}
\begin{proof}
Let $\wt\calH:=\calH_\sigma(\mathbb R^n)$ and $U:=\operatorname{int}(\mathcal X)$. We first observe that the restriction operator
\[
R_{\mathcal X}: \wt\calH\to\calH, \qquad R_{\calX}f = f|_\calX,
\]
is an isometric isomorphism. Indeed, the restriction theorem for RKHSs \cite[Corollary 5.8]{pr} implies that $R_{\mathcal X}$ is surjective and that
\[
\|h\|_{\mathcal H}
=
\inf\left\{
\|f\|_{\widetilde{\mathcal H}}:
f|_{\mathcal X}=h
\right\}.
\]
Moreover, every function in $\widetilde{\mathcal H}$ is real analytic on $\mathbb R^n$, see \cite[Proposition 3.6]{SteinHushScovel06}. Since $U$ is nonempty and open, the identity theorem for real-analytic functions implies that $R_{\mathcal X}$ is injective. Consequently, every $h\in\mathcal H$ has a unique extension to  $\widetilde{\mathcal H}$, and
\[
\|R_{\mathcal X}f\|_{\mathcal H} =\|f\|_{\widetilde{\mathcal H}}.
\]

\medskip
\noindent
\textbf{Necessity.}
Assume that $\mathcal H$ is Koopman invariant, and let $R_U^\calX: \mathcal H\to\mathcal H_\sigma(U)$ denote restriction from $\calX$ to $U$. This operator is bounded, and hence so is $T : \wt\calH\to\calH_\sigma(U)$,
\[
T := R_U^\calX\mathcal K_{\mathcal H} R_{\mathcal X}.
\]
By construction, we have
\[
(Tf)(x)=f(F(x)), \qquad x\in U.
\]
Theorem~1.1 of \cite{Ikeda22}, applied to the Gaussian kernel, now implies that there exist $b\in\mathbb R^n$ and an invertible $A\in\R^{n\times n}$ with $\|A\|_{2\to2}\leq 1$ such that
\[
F(x)=Ax+b,\qquad x\in U.
\]
It remains to extend this identity from $U$ to all of $\mathcal X$. For this, fix $y\in\mathbb R^n$, and consider the kernel section
\[
\varphi_y(z) := k_\sigma(z,y).
\]
Since $\varphi_y|_{\mathcal X}\in\mathcal H$, Koopman invariance gives
\[
h_y:= \varphi_y\circ F \in\mathcal H.
\]
Let $\widetilde h_y\in\widetilde{\mathcal H}$ be its unique extension to $\mathbb R^n$. For $x\in U$ we have
\[
\widetilde h_y(x) = k_\sigma(F(x),y) = k_\sigma(Ax+b,y).
\]
Both functions
\[
x\longmapsto\widetilde h_y(x) \quad\text{and}\quad x\longmapsto k_\sigma(Ax+b,y)
\]
are real analytic on $\mathbb R^n$. By the identity theorem,
\[
\widetilde h_y(x) = k_\sigma(Ax+b,y) \qquad \text{for all }x\in\mathbb R^n.
\]
In particular,
\[
k_\sigma(F(x),y) = k_\sigma(Ax+b,y)
\]
for all $x\in\mathcal X$ and $y\in\mathbb R^n$. Choosing $y=F(x)$ yields
\[
1 = \exp\left(-\frac{\|Ax+b-F(x)\|_2^2} {2\sigma^2} \right),
\]
and therefore
\[
F(x)=Ax+b \qquad\text{for all }x\in\mathcal X.
\]

\medskip
\noindent
\textbf{Sufficiency and operator norm.}
Conversely, suppose that $F(x)=Ax+b$, $x\in\mathcal X$, where $A\in\mathrm{GL}_n(\mathbb R)$ and
$\|A\|_{2\to2}\le 1$. We extend $F$ to the affine map $\wt F(x) = Ax+b$ on $\R^n$. We first prove that composition with $\wt F$ defines a bounded operator $\calK_{\wt\calH}$ on  $\widetilde{\mathcal H}$ and determine its norm.

By \cite[Theorem 10.12]{Wend04}, $\wt\calH$ consists exactly of those $f\in L^2(\R^n)$ for which the following expression is finite:
\[
\|f\|_{\widetilde{\mathcal H}}^2 = c_\sigma \int_{\mathbb R^n} |\widehat f(\xi)|^2 \exp\left(\frac{\sigma^2\|\xi\|_2^2}{2}\right)d\xi,
\]
where $c_\sigma = (2\pi)^{-n}(2\pi\sigma^2)^{-n/2}$ and
\[
\widehat f(\xi) = \int_{\mathbb R^n} f(x)e^{-i x\cdot\xi}\,dx,
\]
initially defined on $L^1(\R^n)\cap L^2(\R^n)$ and extended to $L^2(\R^n)$ by the Plancherel theorem. Since $A$ is invertible,
\[
\widehat{f\circ \wt F}(\xi) = \frac{1}{|\det A|} \exp\left(i\langle A^{-T}\xi,b\rangle\right)\widehat f(A^{-T}\xi).
\]
Consequently, the change of variables $\xi=A^T\eta$ gives
\begin{align*}
\|f\circ \wt F\|_{\widetilde{\mathcal H}}^2
&= \frac{c_\sigma}{|\det A|} \int_{\mathbb R^n}|\widehat f(\eta)|^2\exp\left(\frac{\sigma^2\|A^T\eta\|_2^2}{2}\right)d\eta.
\end{align*}
Because $\|A\|_{2\to2}\le 1$, we have
\[
\|A^T\eta\|_2\leq\|\eta\|_2.
\]
It follows that
\[
\|f\circ\wt F\|_{\widetilde{\mathcal H}}^2 \leq \frac{1}{|\det A|} \|f\|_{\widetilde{\mathcal H}}^2.
\]
Thus composition with $\wt F$ defines a bounded operator $\calK_{\wt\calH}$ on $\widetilde{\mathcal H}$, and
\[
\|\calK_{\wt\calH}\|_{\widetilde{\mathcal H}\to\widetilde{\mathcal H}}\leq|\det A|^{-1/2}.
\]
To prove the reverse inequality, let $\varepsilon>0$ and choose a nonzero function $f_\varepsilon\in\widetilde{\mathcal H}$ such that
\[
\operatorname{supp}\widehat f_\varepsilon \subseteq B_\varepsilon(0).
\]
For $\eta\in B_\varepsilon(0)$,
\[
\|A^T\eta\|_2^2-\|\eta\|_2^2 \geq-\varepsilon^2.
\]
Therefore,
\[
\frac{\|f_\varepsilon\circ\wt F\|_{\widetilde{\mathcal H}}^2}{\|f_\varepsilon\|_{\widetilde{\mathcal H}}^2} \geq \frac{1}{|\det A|}\exp\left(-\frac{\sigma^2\varepsilon^2}{2}\right).
\]
Letting $\varepsilon\downarrow0$ proves \eqref{e:gaussian_norm}. Finally, since $F=\wt F|_{\mathcal X}$, we have
\[
\mathcal K_{\mathcal H} = R_{\mathcal X}\,\calK_{\wt\calH}\,R_{\mathcal X}^{-1}.
\]
As $R_{\mathcal X}$ is an isometric isomorphism, this identity proves Koopman invariance of $\mathcal H$ and yields
\[
\|\mathcal K_{\mathcal H}\|_{\mathcal H\to\mathcal H} = \|\calK_{\wt\calH}\|_{\widetilde{\mathcal H}\to\widetilde{\mathcal H}} = |\det A|^{-1/2},
\]
which concludes the proof.
\end{proof}

\begin{rem} Such deterministic rigidity fails for Markov dynamics: the Gaussian RKHS can be Koopman invariant for stochastic systems \cite{Philipp2025}.
\end{rem}

\section{Compact embeddings of reproducing kernel Hilbert spaces}\label{a:RKHS_compact_embeddings}
The range-space arguments in Theorems~\ref{t:RKHS_comp_det} and \ref{t:RKHS_comp} reduce compactness to an inclusion between two RKHSs. This appendix characterizes such inclusions by uniform spectral decay of finite Gram matrices. The abstract Hilbert-space lemma below supplies the passage from finite compressions to the full inclusion. For a closed subspace $M$ of a Hilbert space, let $P_M$ denote the orthogonal projection onto $M$. For a finite-rank nonnegative self-adjoint operator $T$, let $\lambda_n(T)$ be its $n$th eigenvalue in nonincreasing order, counting multiplicities, and zero for $n>\rank T$.

\begin{lem}\label{l:A_compact}
A bounded nonnegative self-adjoint operator $A$ on a Hilbert space $\calH$
is compact if and only if
\begin{align}\label{e:sup_Efinite}
\sup_{\dim E<\infty}\lambda_n(P_EAP_E)\longrightarrow0 \qquad(n\to\infty).
\end{align}
\end{lem}
\begin{proof} For $\varepsilon>0$, let $M_\varepsilon=\ran\one_{[\varepsilon,\infty)}(A)$. By the spectral theorem, $A$ is compact exactly when every $M_\varepsilon$ is finite-dimensional. Indeed, in that case the finite-rank operators $A\one_{[1/m,\infty)}(A)$ converge to $A$ in norm; conversely, a compact operator cannot be bounded below by a positive constant on an infinite-dimensional subspace.

If $r=\dim M_\varepsilon<\infty$, then $0\le A\le\varepsilon I$ on $M_\varepsilon^\perp$. For every finite-dimensional $E$, we have $\dim(E\cap M_\varepsilon^\perp)\ge\dim E-r$ and $\langle Au,u\rangle\le\varepsilon\|u\|^2$ for $u\in E\cap M_\varepsilon^\perp$. The min--max principle for $P_EAP_E|_E$ therefore gives $\lambda_{r+1}(P_EAP_E)\le\varepsilon$ (trivially if $\dim E\le r$), proving \eqref{e:sup_Efinite}. Conversely, if $M_\varepsilon$ is infinite-dimensional, choose an $n$-dimensional $E_n\subset M_\varepsilon$ for each $n$. Then $P_{E_n}AP_{E_n}|_{E_n}\ge\varepsilon I_{E_n}$ and $\ran(P_{E_n}AP_{E_n}) = E_n$, so $\lambda_n(P_{E_n}AP_{E_n})\ge\varepsilon$ for every $n$. Thus, \eqref{e:sup_Efinite} fails.
\end{proof}

\begin{thm}\label{t:RKHS_comp_emb} Let $k_1$ and $k_2$ be Hermitian positive-definite kernels on $\calX$, with RKHSs $\calH_1$ and $\calH_2$, respectively. For an ordered finite tuple $X=(x_1,\ldots,x_m)$, denote
\begin{align*}
K_{k_j,X}=(k_j(x_r,x_s))_{r,s=1}^m
\qquad\text{and}\qquad
B_X=(K_{k_2,X}^\dagger)^{1/2}K_{k_1,X} (K_{k_2,X}^\dagger)^{1/2}.
\end{align*}
Then $\calH_1$ is compactly embedded in $\calH_2$ if and only if $k_1\preceq c^2k_2$ for some $c>0$ and
\begin{align}\label{e:supla_n}
\sup_{X\text{ finite}}\lambda_n(B_X)\longrightarrow0 \qquad(n\to\infty).
\end{align}
\end{thm}

\begin{rem}
The domination condition $k_1\preceq c^2k_2$ is equivalent to
\begin{align*}
\ker K_{k_2,X}\subset\ker K_{k_1,X}, \qquad \|B_X\|_{2\to2}\le c^2
\end{align*}
for every finite $X$; see \cite[Remark B.2]{KohnPhil24}. Thus the boundedness and compactness conditions are both expressed through the same finite Gram matrices.
\end{rem}

\begin{proof}[Proof of Theorem \rmref{t:RKHS_comp_emb}]
Assume first that $\calH_1\subset\calH_2$ continuously. Let $J:\calH_1\to\calH_2$ be the inclusion. The operator $A=JJ^*$ is bounded, nonnegative, and self-adjoint, while $J$ is compact if and only if $A$ is compact. For the canonical features $\phi_{j,x}=k_j(x,\cdot)$, $J^*\phi_{2,x}=\phi_{1,x}$.

Fix $X=(x_1,\ldots,x_m)$ and let $T_X:\C^m\to\calH_2$ satisfy $T_Xe_r=\phi_{2,x_r}$. With $G_{j,X}:=\overline{K_{k_j,X}}$ and the linear-in-the-first-variable convention fixed in Subsection~\ref{s:invariance:det:RKHS}, the exact matrix identities are
\begin{align*}
T_X^*T_X=G_{2,X},\qquad T_X^*AT_X=G_{1,X}.
\end{align*}
With $E_X=\ran T_X$, the map
\begin{align*}
U_X=T_X(G_{2,X}^\dagger)^{1/2}: \ran G_{2,X}\longrightarrow E_X
\end{align*}
is unitary. The matrix $\overline{B_X}$ vanishes on $\ker G_{2,X}$, and
\begin{align*}
U_X^*(P_{E_X}A|_{E_X})U_X =\overline{B_X}|_{\ran G_{2,X}}.
\end{align*}
Complex conjugation is antiunitary and carries the nonnegative Hermitian matrix $B_X$ to $\overline{B_X}$; the two matrices therefore have the same eigenvalues. Hence the eigenvalues of $B_X$ are precisely those of $P_{E_X}AP_{E_X}|_{E_X}$, completed by zeros. If $J$ is compact, Lemma~\ref{l:A_compact} now gives \eqref{e:supla_n}.

Conversely, suppose that $k_1\preceq c^2k_2$ and that \eqref{e:supla_n} holds. Aronszajn's inclusion theorem supplies the bounded inclusion $J$, and hence $A=JJ^*$. There is nothing to prove if $A=0$. Otherwise, given $\varepsilon>0$, choose $N$ such that
\begin{align*}
\sup_{X\text{ finite}}\lambda_n(P_{E_X}AP_{E_X}) <\frac{\varepsilon}{2}\qquad(n\ge N).
\end{align*}
Let $E\subset\calH_2$ be finite-dimensional. Approximate an orthonormal basis of $E$ by linearly independent vectors in the dense feature span. For a sufficiently close approximation, their span $E'\subset E_X$ for some finite $X$ and
\begin{align*}
\|P_E-P_{E'}\|<\frac{\varepsilon}{4\|A\|}.
\end{align*}
For completeness, if a synthesis operator $V$ is close to an isometry, then the projection onto its range is
\begin{align*}
P_{\ran V}=V(V^*V)^{-1}V^*.
\end{align*}
This proves the required continuity of the projections. Since $E'\subset E_X$, compression monotonicity gives
\begin{align*}
\lambda_n(P_{E'}AP_{E'})\leq\lambda_n(P_{E_X}AP_{E_X}).
\end{align*}
The min--max principle and eigenvalue perturbation now yield, for $n\ge N$,
\begin{align*}
\lambda_n(P_EAP_E) &\le\|P_EAP_E-P_{E'}AP_{E'}\| +\lambda_n(P_{E'}AP_{E'})\\ &\le2\|A\|\|P_E-P_{E'}\| +\lambda_n(P_{E_X}AP_{E_X}) <\varepsilon.
\end{align*}
Take the supremum over finite-dimensional $E$ and apply Lemma~\ref{l:A_compact}. It follows that $A$ is compact, and thus so is $J$.
\end{proof}

\section{Compact subsets of \texorpdfstring{$L^1$}{L1}}\label{a:L1_compactness}
The exact kernel criterion in Theorem~\ref{t:compact2:Lp} rests on a compactness theorem for arbitrary subsets of $L^1$. Weak relative compactness controls concentration on small sets and escape to infinity, but norm compactness additionally requires local control of mean oscillations. The condition below is usually called the local Bocce condition; see \cite[Theorem 4.7.29]{Bogachev07} and \cite{Girardi1991}.

\begin{thm}\label{t:Lp_comp_family} 
Let $\mu$ be a $\sigma$-finite measure on the measurable space $(\calX,\Sigma)$. A set $\mathscr F\subset L^1(\mu;\C)$ is relatively compact in $L^1(\mu;\C)$ if and only if the following conditions hold:
\begin{enumerate}
\item[{\rm (a)}] $\mathscr F$ is bounded: $\sup_{f\in\mathscr F}\|f\|_1<\infty$;
\item[{\rm (b)}] $\mathscr F$ has uniformly absolutely continuous integrals:
\begin{align*}
\lim_{\delta\downarrow0} \sup_{\substack{f\in\mathscr F,\ \Delta\in\Sigma\\\mu(\Delta)\leq\delta}} \int_\Delta|f|\,d\mu=0;
\end{align*}
\item[{\rm (c)}] $\mathscr F$ has uniformly small tails: for every $\varepsilon>0$ there is an $E\in\Sigma$ with $\mu(E)<\infty$ such that
\begin{align*}
\sup_{f\in\mathscr F}\int_{\calX\backslash E}|f|\,d\mu<\varepsilon;
\end{align*}
\item[{\rm (d)}] $\mathscr F$ satisfies the following local mean-oscillation condition: for every $\varepsilon>0$ and every $\Delta\in\Sigma$ with $0<\mu(\Delta)<\infty$, there are measurable $\Delta_1,\ldots,\Delta_n\subset\Delta$ with $0<\mu(\Delta_j)<\infty$ such that, for each $f\in\mathscr F$, some $j\in\{1,\ldots,n\}$ satisfies
\begin{align}\label{e:bocce_oscillation}
\frac1{\mu(\Delta_j)^2} \int_{\Delta_j}\!\int_{\Delta_j}|f(x)-f(z)|\,d\mu(x)\,d\mu(z) <\varepsilon.
\end{align}
\end{enumerate}
\end{thm}

\begin{rem}
Conditions {\rm(a)--(c)} are equivalent to weak relative compactness in $L^1(\mu)$; see \cite[Theorem 4.7.20]{Bogachev07}. For $A\in\Sigma$ with $0<\mu(A)<\infty$, let
\begin{align*}
\langle f\rangle_A=\frac1{\mu(A)}\int_A f\,d\mu, \qquad \operatorname{osc}_A(f) =\frac1{\mu(A)}\int_A|f-\langle f\rangle_A|\,d\mu = \left\langle|f - \langle f\rangle_A|\right\rangle_A.
\end{align*}
Jensen's inequality and the triangle inequality give
\begin{align*}
\operatorname{osc}_A(f) \leq\frac1{\mu(A)^2}\int_A\!\int_A|f(x)-f(z)|\,d\mu(x)d\mu(z) \leq2\operatorname{osc}_A(f).
\end{align*}
Consequently, condition~{\rm(d)} is equivalent to the formulation in \cite[Theorem 4.7.29]{Bogachev07} and \cite{Girardi1991}.
\end{rem}

\begin{proof}[Proof of Theorem \rmref{t:Lp_comp_family}]
First assume that $\mathscr F$ is relatively compact. Boundedness gives {\rm(a)}. To prove {\rm(b)}, fix $\eta>0$ and choose a finite $\eta/2$-net $f_1,\ldots,f_N$ for $\mathscr F$. Absolute continuity of the finitely many integrals $\int|f_j|\,d\mu$ supplies $\delta>0$ such that $\mu(A)\leq\delta$ implies $\int_A|f_j|\,d\mu<\eta/2$ for every $j$. Approximation of any $f\in\mathscr F$ by one of the $f_j$ proves {\rm(b)}.

For {\rm(c)}, choose increasing $E_m\in\Sigma$ with $\mu(E_m)<\infty$ and $E_m\uparrow\calX$. Given $\eta>0$, take a finite $\eta/2$-net as above and then choose $m$ so large that $\int_{\calX\backslash E_m}|f_j|\,d\mu<\eta/2$ for every $j$. The triangle inequality gives the required uniform tail bound.

For {\rm(d)}, fix $0<\mu(\Delta)<\infty$ and $\varepsilon>0$, and let $K$ be the compact closure of $\mathscr F$. For each $f\in K$, partition $\C$ into countably many Borel sets of diameter less than $\varepsilon/2$. At least one inverse image $A_f\subset\Delta$ has positive measure, and
\begin{align*}
D_{A_f}(f):=\frac1{\mu(A_f)^2} \int_{A_f}\!\int_{A_f}|f(x)-f(z)|\,d\mu(x)d\mu(z) <\frac\varepsilon2.
\end{align*}
For every positive finite-measure $A$ and $f,g\in L^1(\mu)$,
\begin{align*}
|D_A(f)-D_A(g)| \leq\frac{2}{\mu(A)}\|f-g\|_{L^1(\mu)}.
\end{align*}
The sets $\{g:D_{A_f}(g)<\varepsilon\}$ form an open cover of $K$. A finite subcover therefore supplies the finitely many sets required in {\rm(d)}.

Conversely, assume {\rm(a)--(d)} and fix $\eta>0$. By {\rm(c)}, choose $E\in\Sigma$ with $\mu(E)<\infty$ and
\begin{align*}
\sup_{f\in\mathscr F}\|\one_{\calX\backslash E}f\|_1<\frac\eta3.
\end{align*}
If $\mu(E)=0$, the zero function alone is already an $\eta$-net for $\mathscr F$. Otherwise, the restricted family $\mathscr F_E=\{\one_Ef:f\in\mathscr F\}$ satisfies {\rm(a)}, {\rm(b)}, and the finite-measure version of condition~{\rm(d)} on $E$. The real families $\operatorname{Re}\mathscr F_E$ and $\operatorname{Im}\mathscr F_E$ satisfy the same conditions because the real and imaginary parts do not increase the pairwise oscillation. The finite-measure compactness criterion \cite[Theorem 4.7.29]{Bogachev07} therefore makes both real families relatively compact. Identifying complex $L^1$ with the product of its real and imaginary parts shows that $\mathscr F_E$ is relatively compact in $L^1(E;\C)$.

Choose a finite $\eta/3$-net $g_1,\ldots,g_M$ for $\mathscr F_E$, extending the $g_j$ by zero off $E$. For every $f\in\mathscr F$, some $j$ satisfies
\begin{align*}
\|f-g_j\|_1 \leq\|\one_Ef-g_j\|_1+\|\one_{\calX\backslash E}f\|_1 <\frac{2\eta}{3}.
\end{align*}
Thus $\mathscr F$ is totally bounded. Since $L^1(\mu;\C)$ is complete, $\mathscr F$ is relatively compact.
\end{proof}

\bibliographystyle{abbrvnat}
\bibliography{References}

\bigskip
\section*{Author affiliations}
\end{document}